\RequirePackage{fix-cm}
\documentclass[smallextended]{svjour3}       % onecolumn (second format)
\smartqed  % flush right qed marks, e.g. at end of proof
\usepackage{graphicx}
\usepackage{hyperref}

\usepackage{booktabs}       % professional-quality tables

\usepackage{hyperref} %url

\usepackage[most]{tcolorbox}

\usepackage{pifont}
\newcommand{\cmark}{\ding{51}}
\newcommand{\xmark}{\ding{55}}

\usepackage[font=footnotesize]{subfig} %subfigure
\usepackage[export]{adjustbox} % subfig align

\usepackage{mathtools} % math
\usepackage{amsmath,amsfonts,amssymb}
\usepackage{bm} %bold
\usepackage{bbm}
\DeclarePairedDelimiter{\norm}{\lVert}{\rVert}
\newcommand{\argmin}{\operatornamewithlimits{argmin}}

\usepackage{algorithm} %algo
\usepackage{algpseudocode}
\makeatletter
\renewcommand{\ALG@beginalgorithmic}{\small}
\makeatother

\usepackage{bold-extra} % bold

\usepackage[numbers]{natbib} %cite

\usepackage{listings} %code
\usepackage{xcolor} 
\definecolor{codegreen}{rgb}{0,0.6,0}
\definecolor{codegray}{rgb}{0.5,0.5,0.5}
\definecolor{codepurple}{rgb}{0.58,0,0.82}
\definecolor{backcolour}{rgb}{0.97,0.97,0.98}
\definecolor{jordyblue}{rgb}{0.55,0.72,0.89}
\definecolor{mediumblue}{rgb}{0,0,0.8}
\definecolor{fungreen}{rgb}{0.07,0.36,0.21}
\lstdefinestyle{mystyle}{
    backgroundcolor=\color{backcolour},   
    commentstyle=\color{jordyblue},
    keywordstyle=\color{mediumblue},
    numberstyle=\tiny\color{codegray},
    stringstyle=\color{fungreen},
    basicstyle=\ttfamily\footnotesize,
    breakatwhitespace=false,         
    breaklines=true,                 
    captionpos=b,                    
    keepspaces=true,                 
    numbers=left,                    
    numbersep=5pt,                  
    showspaces=false,                
    showstringspaces=false,
    showtabs=false,                  
    tabsize=2
}
\newcommand{\pyepo}{\textsl{PyEPO}}
\newcommand{\gurobi}{\textsl{Gurobi}}
\newcommand{\copt}{\textsl{COPT}}
\newcommand{\cplex}{\textsl{CPLEX}}

\newcommand{\glpk}{\textsl{GLPK}}

\newcommand{\cvxpy}{\textsl{CVXPY}}

\newcommand{\gurobipy}{\textsl{GurobiPy}}
\newcommand{\coptpy}{\textsl{COPTPy}}
\newcommand{\pyomo}{\textsl{Pyomo}}
\newcommand{\sklearn}{\textsl{Scikit-Learn}}
\newcommand{\autosklearn}{\textsl{Auto-Sklearn}}
\newcommand{\pytorch}{\textsl{PyTorch}}
\newcommand{\tensorflow}{\textsl{TensorFlow}}
\newcommand{\mxnet}{\textsl{MXNet}}

\newcommand{\optnet}{\textbf{\texttt{OptNet}}}
\newcommand{\qpth}{\textbf{\texttt{qpth}}}
\newcommand{\dqp}{\textbf{\texttt{DQP}}}
\newcommand{\qptl}{\textbf{\texttt{QPTL}}}
\newcommand{\mipaal}{\textbf{\texttt{MIPaaL}}}
\newcommand{\intopt}{\textbf{\texttt{IntOpt}}}
\newcommand{\cvxpylayers}{\textbf{\texttt{CvxpyLayers}}}
\newcommand{\spo}{\textbf{\texttt{SPO+}}}
\newcommand{\sporel}{\textbf{\texttt{SPO+ Rel}}}
\newcommand{\spows}{\textbf{\texttt{SPO+ WS}}}
\newcommand{\dbb}{\textbf{\texttt{DBB}}}
\newcommand{\dbbrel}{\textbf{\texttt{DBB Rel}}}

\newcommand{\dpo}{\textbf{\texttt{DPO}}}

\newcommand{\pfyl}{\textbf{\texttt{PFYL}}}
\newcommand{\pfylrel}{\textbf{\texttt{PFYL Rel}}}

\newcommand{\spolo}{\textbf{\texttt{SPO+ L1}}}
\newcommand{\spolt}{\textbf{\texttt{SPO+ L2}}}

\newcommand{\dbblo}{\textbf{\texttt{DBB L1}}}
\newcommand{\dbblt}{\textbf{\texttt{DBB L2}}}

\newcommand{\pfyllo}{\textbf{\texttt{PFYL L1}}}
\newcommand{\pfyllt}{\textbf{\texttt{PFYL L2}}}

\newcommand{\nce}{\textbf{\texttt{NCE}}}
\newcommand{\ltr}{\textbf{\texttt{LTR}}}

\begin{document}

\captionsetup[subfigure]{labelformat=empty} % remove subfig label

\title{PyEPO: A PyTorch-based End-to-End Predict-then-Optimize Library for Linear and Integer Programming%\thanks{Grants or other notes
}
\titlerunning{PyEPO: A Predict-then-Optimize Library for Linear and Integer Programming}        % if too long for running head

\author{Bo Tang         \and
        Elias B.~Khalil %etc.
}


\institute{\href{https://orcid.org/0000-0002-6035-5167}
           {\includegraphics[scale=0.06]{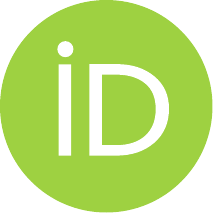}\hspace{1mm}Bo Tang} \at
              Department of Mechanical and Industrial Engineering, University of Toronto\\
              5 King's College Rd, Toronto, ON M5S 3G8\\
              \email{botang@mie.utoronto.ca}           %  \\
           \and
           \href{https://orcid.org/0000-0001-5844-9642}
	       {\includegraphics[scale=0.06]{orcid.pdf}\hspace{1mm}Elias B.~Khalil} \at
              Department of Mechanical and Industrial Engineering, University of Toronto\\
              5 King's College Rd, Toronto, ON M5S 3G8\\
              \email{khalil@mie.utoronto.ca}
}

\date{Received: date / Accepted: date}

\maketitle

\begin{abstract}
% The increasing reliance on decision-making inspires the development of optimization methods and toolboxes. 
In deterministic optimization, it is typically assumed that all problem parameters are fixed and known. In practice, however, some parameters may be a priori unknown but can be estimated from contextual information. A typical predict-then-optimize approach separates predictions and optimization into two distinct stages. Recently, end-to-end predict-then-optimize has emerged as an attractive alternative. This work introduces the \pyepo{} package, a \pytorch-based end-to-end predict-then-optimize library in Python. To the best of our knowledge, \pyepo{} (pronounced like \textit{pineapple} with a silent ``n") is the first such generic tool for linear and integer programming with predicted objective function coefficients. It includes various algorithms such as surrogate decision losses, black-box solvers, and perturbated methods. \pyepo{} offers a user-friendly interface for defining new optimization problems, applying state-of-the-art algorithms, and using custom neural network architectures. We conducted experiments comparing various methods on problems such as the Shortest Path, the Multiple Knapsack, and the Traveling Salesperson Problem, discussing empirical insights that may guide future research.~\pyepo{} and its documentation are available at \url{https://github.com/khalil-research/PyEPO}.

\end{abstract}

\section{Introduction}
\label{sec:intro}

Predictive modeling is ubiquitous in real-world decision-making. For instance, in many applications, the objective function coefficients, such as travel time in a routing problem, electricity price in a power system, and asset return in portfolio optimization, are unknown at the time of decision-making. In this work, we focus on the commonly used paradigm of prediction followed by optimization in the context of linear programs or integer linear programs. In this paradigm, it is assumed that a set of features describes an instance of the optimization problem. A regression model maps the features to the (unknown) objective function coefficients. Then, a deterministic optimization problem is solved to obtain a solution. Due to its wide applicability and simplicity compared to other frameworks for optimization under uncertain parameters, the predict-then-optimize paradigm has garnered increasing attention in recent years.

% idea: this is the way predict-then-optimize training is done
One natural approach is to proceed in two stages, first training an accurate predictive model and then solving the downstream optimization problem using predicted coefficients. The advantage of the two-stage approach is the direct utilization of existing machine learning methods and optimization solvers. However, prediction errors, such as mean squared error, cannot adequately measure the quality of decisions. While a perfect prediction would always yield an optimal decision, learning a prediction model without errors is impractical. \citet{bengio1997using,ford2015beware}, and~\citet{elmachtoub2021smart} reported that training a predictive model based on prediction error leads to worse decisions compared to directly considering the decision error. Consequently, the state-of-the-art alternative is to integrate optimization into prediction, taking into account its impact on the decision, the so-called end-to-end learning framework.

% idea: these solvers can provide a solution, but improving the predictive model based on it requires additional computations and integration which is not provided
End-to-end predict-then-optimize requires embedding an optimization solver into the model training loop. Classical solution approaches for linear and integer linear models, including graph algorithms, linear programming, integer programming, constraint programming, etc., are well-established and efficient in practice. In addition, commercial solvers such as \gurobi{} \cite{gurobi},\cplex{} \cite{cplex2009v12} and \copt{} \cite{copt} are highly optimized, allowing users to easily express business or academic problems as optimization models without a deep understanding of the underlying theory and algorithms. However, embedding a solver for end-to-end learning necessitates additional computation and integration (e.g., gradient calculation), which current software does not provide.

% idea: although these tools make it easy to train ML models, they cannot seamlessly integrate training loss functions that require solving a discrete optimization problem
On the other hand, the field of machine learning has witnessed tremendous growth in recent decades. Breakthroughs in deep learning have led to remarkable improvements in several complex tasks. As a result, neural networks now pervade disparate applications that span computer vision, natural language, and planning, among others. Python-based machine learning frameworks such as \sklearn{} \cite{scikit-learn}, \tensorflow{} \cite{abadi2016tensorflow}, \pytorch{} \cite{NEURIPS2019_9015}, \mxnet{} \cite{chen2015mxnet}, etc., have been developed and extensively used for research and production needs. Although deep learning has proven highly effective in regression and classification, it lacks the capability to handle constrained optimization such as integer linear programming.

% idea: here are some prominent attempts at bridging this gap
Since \citet{amos2017optnet} first introduced a neural network layer for generic mathematical optimization, there have been several prominent attempts to bridge the gap between optimization solvers and deep learning frameworks. A critical component is typically a differentiable block for optimization tasks. With a differentiable optimizer or a surrogate decision loss function, neural network packages enable the computation of gradients for optimization operations and then update predictive model parameters based on a loss function that depends on decision quality.  

% idea: here are the issues with existing codebases 
Although research code implementing a number of predict-then-optimize training algorithms has been made available for particular classes of optimization problems and/or predictive models ~\cite{gould2016differentiating, donti2017task, wilder2019melding, ferber2020mipaal, mandi2020interior, elmachtoub2021smart, poganvcic2019differentiation, mandi2020interior, NIPS2017_192fc044, agrawal2019differentiating, agrawal2019differentiable}, there is a pressing need for a generic end-to-end learning framework, particularly for linear and integer programming. In this paper, we propose the open-source software package \pyepo{} that aims to customize and train end-to-end predict-then-optimize for optimization problems with linear objective functions, such as linear and integer programming. Our contributions are as follows:

\begin{enumerate}
    \item We implement \spo{} (``Smart Predict-then-Optimize+") loss \cite{elmachtoub2021smart}, \dbb{} (differentiable black-box) solver \cite{poganvcic2019differentiation}, \dpo{} (differentiable perturbed optimizer), \pfyl{} (perturbed Fenchel-Young loss), among other typical end-to-end methods for linear and integer programming.

    \item We build \pyepo{} based on \pytorch{}. As one of the most popular deep learning frameworks, \pytorch{} facilitates easy integration and usage of any deep neural network.

    \item We provide interfaces to the Python-based optimization modeling frameworks \gurobipy{}, \coptpy{} and \pyomo{}. These high-level modeling languages allow non-specialists to formulate optimization models with \pyepo{}.

    \item We enable parallel computing for the forward and backward pass in \pyepo{}. Optimizations during training are carried out in parallel, allowing users to leverage multiple processors to reduce training time.

    \item We present new benchmark datasets for end-to-end predict-then-optimize, allowing us to compare the performance of different approaches.

    \item We conduct and analyze a comprehensive set of experiments for end-to-end predict-then-optimize. The performance of different methods across various datasets demonstrate the competitiveness of end-to-end learning and highlight the surprising effect of relaxations and regularization.
\end{enumerate}

\pyepo{} focuses on the development of software that offers interfaces for various methodologies. It also provides new benchmark datasets and conducts experiments. In addition to \pyepo{}, there are comprehensive surveys \cite{sadana2023survey, mandi2023decision} that complement this paper. \citet{sadana2023survey} provides a theoretical and algorithmic review and taxonomy for more general contextual optimization, including predict-then-optimize approaches. Meanwhile, the survey by \citet{mandi2023decision} delivers an extensive analysis and experiments on end-to-end predict-then-optimize methods. Notably, the experiments on \citet{mandi2023decision} differ from \pyepo{}: they test more methods on several different datasets but do not explore the implications of relaxation or coefficient regularization.

\section{Related work}
\label{sec:liter}

% A variety of approaches built on an end-to-end learning framework have emerged in recent years. 
In early work on the topic, \citet{bengio1997using} introduced a differentiable portfolio optimizer and suggested that direct optimization with financial criteria performs better in neural networks than the mean squared error of the predicted values. \citet{kao2009directed} trained a linear regressor with a convex combination of prediction error and decision error, but only considered unconstrained quadratic programming. \citet{domke2012generic} investigated generic gradient descent methods to minimize unconstrained energy functions. 

More recently, interest has shifted towards constrained optimization problems. For instance, \citet{gould2016differentiating} allows differentiating bilevel optimization problems with and without constraints. Subsequent studies have comprehensively explored differentiable constraint optimization, covering various types of optimization problems such as quadratic programming, linear programming, and integer linear programming. A comparison of gradient-based methodologies for end-to-end constrained optimization is presented in Table \ref{tab:comp} and Table \ref{tab:cost}. Additionally, \citet{elmachtoub2020decision} proposed a tree-based approach that does not rely on gradients.

%\subsection{Solution Prediction without Differentiable Optimization}

%Since end-to-end learning for constrained optimization is promising, approaches to directly predict optimization solutions without the use of a differentiable optimization solver, such as Hopfield Networks \citep{hopfield1982neural} and Pointer Networks \citep{vinyals2015pointer} have emerged. Although these machine learning architectures predict approximate solutions with good computational efficiency, they are difficult to guarantee feasibility and generalize to more complex constraints. Therefore, the state-of-art alternative is end-to-end predict-then-optimize, which integrates optimization into prediction and takes into account the impact on the decision during the training procedure.

\subsection{KKT-based Implicit Differentiation}
\label{subsec:kkt}

% Gradient-based methods begin by asking whether the gradient of an optimal solution w.r.t. the parameters of the predictive model can be computed efficiently.
Gradient-based end-to-end learning requires well-defined first-order derivatives of optimization. The Karush-Kuhn-Tucker (KKT) conditions become an attractive option because they make the optimization problem with hard constraints differentiable.

\citet{amos2017optnet} proposed \optnet{}, which derives implicit differentiation of constrained quadratic programs from KKT conditions. Building on \optnet{}, \citet{donti2017task} investigated a general end-to-end framework, \dqp{}, for learning with constrained quadratic programming, which simultaneously obtains optimal solutions and their gradients. Although linear programming is a special case of quadratic programming, \dqp{} has no ability to tackle linear objective functions because the optimal solution exhibits piecewise constancy, leading to zero gradients almost everywhere.

Furthermore, \citet{wilder2019melding} added a small quadratic objective term to \dqp{} to ensure nonzero gradients for linear programming, resulting in \qptl{}. \citet{wilder2019melding} also discussed relaxation and rounding for the approximation of binary problems. \citet{ferber2020mipaal} followed up on \qptl{} with \mipaal{}, a cutting-plane approach to support (mixed) integer programming. With the cutting-plane method, \mipaal{} generates (potentially exponentially many) valid cuts to convert a discrete problem into an equivalent continuous problem, which is theoretically sound for combinatorial optimization, but extremely time-consuming. In addition, \citet{mandi2020interior} introduced \intopt{} to compute gradients for linear programming with the log-barrier term instead of the quadratic term of \qptl{}. Except for \mipaal{} \cite{ferber2020mipaal}, end-to-end learning approaches for (mixed) integer programming use linear relaxation during training but evaluate with optimal integer solutions at test time.

In addition to \dqp{} and its extension, \citet{agrawal2019differentiable} introduced a generic framework, \cvxpylayers{}, which differentiates through KKT conditions of the conic program. The central idea in \cvxpylayers{} is to canonicalize disciplined convex programming as conic programming so that it is applicable to a wider range of convex optimization problems. Similarly to \dqp{}, \cvxpylayers{} encounters difficulties when differentiating through linear programs. Thus, a squared solution is used as a heuristic.

However, these KKT-based methods require a specific quadratic or conic programming solver, and linear programming necessitates additional objective terms. Therefore, the efficiency and accuracy of the above implementations is not comparable to commercial MILP solvers such as \gurobi{} \cite{gurobi} and \cplex{} \cite{cplex2009v12}. Moreover, these methods are designed for continuous problems and do not naturally support discrete optimization.

\begin{table}[ht]
    \begin{center}
    \small
    \resizebox{1.0\textwidth}{!}{
    \begin{tabular}{l|l|lllll}
    \hline
        Method                                       &In PyEPO & w/ Unk Constr & Disc Var & Lin Obj & Quad Obj \\ \hline
        \dqp{} \cite{donti2017task} \href{https://github.com/locuslab/e2e-model-learning}{[code]}                   & \xmark & \cmark & \xmark & \xmark & \cmark \\
        \qptl{} \cite{wilder2019melding} \href{https://github.com/bwilder0/aaai_melding_code}{[code]}              & \xmark & \xmark & \xmark & \cmark & \cmark \\
        \mipaal{} \cite{ferber2020mipaal}               & \xmark & \xmark & \cmark & \cmark & \xmark \\
        \intopt{} \cite{mandi2020interior} \href{https://github.com/JayMan91/NeurIPSIntopt}{[code]}             & \xmark & \xmark & \xmark & \cmark & \xmark \\
        \cvxpylayers{} \cite{agrawal2019differentiable} \href{https://github.com/cvxgrp/cvxpylayers}{[code]}& \xmark & \cmark & \xmark & \cmark & \cmark \\ \hline
        \spo{} \cite{elmachtoub2021smart} \href{https://github.com/paulgrigas/SmartPredictThenOptimize}{[code]}              & \cmark & \xmark & \cmark & \cmark & \xmark \\
        \sporel{} \cite{mandi2020smart} \href{https://github.com/JayMan91/aaai_predit_then_optimize}{[code]}                & \cmark & \xmark & \cmark & \cmark & \xmark \\
        \spows{} \cite{mandi2020smart}  \href{https://github.com/JayMan91/aaai_predit_then_optimize}{[code]}                & \xmark & \xmark & \cmark & \cmark & \xmark \\
        \dbb{} \cite{poganvcic2019differentiation}  \href{https://github.com/martius-lab/blackbox-backprop}{[code]}    & \cmark & \xmark & \cmark & \cmark & \xmark \\
        \dpo{} \cite{berthet2020learning} \href{https://github.com/google-research/google-research/tree/master/perturbations}{[code]}             & \cmark & \xmark & \cmark & \cmark & \xmark \\
        \pfyl{} \cite{berthet2020learning} \href{https://github.com/google-research/google-research/tree/master/perturbations}{[code]}             & \cmark & \xmark & \cmark & \cmark & \xmark \\ 
        \nce{} \cite{mulamba2020contrastive} \href{https://github.com/CryoCardiogram/ijcai-cache-loss-pno}{[code]}             & \cmark & \xmark & \cmark & \cmark & \cmark \\ 
        \ltr{} \cite{mandi2022decision} \href{https://github.com/JayMan91/ltr-predopt}{[code]}             & \cmark & \xmark & \cmark & \cmark & \cmark \\ \hline
    \end{tabular}}
    \caption{Methodology Comparison}\label{tab:comp}
    {\parbox{4.5in}{
        \small This is a comparison diagram for different methodologies. \\
        The first set of methods uses KKT conditions, and the second part is based on differentiable approximations. \\
        "In PyEPO" denotes whether the method is available in \pyepo{}. \\
        % "Grad" means whether the method computes gradients. \\
        "w/ Unk Constr" denotes whether unknown parameters occur in constraints. \\
        "Disc Var" denotes whether the method supports integer variables. \\
        "Lin Obj" denotes whether the method supports a linear objective function. \\
         "Quad Obj" denotes whether the method supports a quadratic objective function.}
    }
    \end{center}
\end{table}

\subsection{Surrogate Loss or Gradients}
\label{subsec:Appro}

Since KKT conditions may not be ideal for linear objective functions, researchers have also explored the design of surrogate loss or gradients. For example, \citet{elmachtoub2021smart} proposed regret (SPO loss in their paper) to measure decision error. However, the gradient of regret is still zero almost everywhere and is undefined otherwise for the linear objective function. To address this, they introduced \spo{}, a convex and sub-differentiable loss, to ensure a nonzero subgradient for training. As in previous approaches, \spo{} solves an optimization problem in each forward pass. Unlike KKT-based methods, \spo{} is limited to the linear objective function. Because optimization is the computational bottleneck, \citet{mandi2020smart} applied \spo{} to combinatorial problems using relaxation and warm starting techniques. They reported that relaxation reduces solving time at the cost of performance.

\citet{poganvcic2019differentiation} computed a subgradient from continuous interpolation of solutions, an approach they referred to as the ``differentiable black-box solver", \dbb{}. The interpolation approximation is nonconvex but avoids vanishing gradients. Compared to \spo{}, \dbb{} requires an extra optimization for the backward pass, and the loss of \dbb{} is flexible (e.g., the Hamming distance in their paper).

In addition, \citet{berthet2020learning} applied perturbed optimizer \dpo{} by adding random noise to the objective function coefficients so that the piecewise constant solutions are smoothed by probability and further constructed the Fenchel-Young loss \pfyl{} by duality. As a sampling-based method, the number of optimization problems that need to be solved for each iteration is the number of samples. 

Moreover, \citet{mulamba2020contrastive} adopted the noise-contrastive estimate (\nce{}) to compute a surrogate loss. This approach treats a subset of suboptimal feasible solutions as negative samples and maximizes the gap between the optimal solution and these negative samples. Inspired by \nce{}, \citet{mandi2022decision} converted the predict-then-optimize task into a "Learning to Rank" (\ltr{}) problem \cite{liu2009learning}, ranking the subset of feasible solutions based on the objective value. Both \nce{} and \ltr{} offer the flexibility of being unrestricted in terms of the type of optimization problem, which presents a significant advantage. However, to our knowledge, these methods have not been tested on non-linear problems.

As the key algorithms of \pyepo{}, \spo{}, \dbb{}, \dpo{}, and \pfyl{} are further discussed in Section \ref{sec:prel}.

\subsection{Code for Predict-then-Optimize}

\subsubsection{Research Code}
%\noindent\textbf{Research code.}
Table~\ref{tab:comp} lists and links to the codebases that will be discussed next.
\citet{amos2017optnet} developed a \pytorch-based solver \qpth{} for \optnet{}, which efficiently solves quadratic programs and computes their gradients. The solver was based on a primal-dual interior-point method \cite{mattingley2012cvxgen} and can handle batches of quadratic programs on a GPU. Using this solver, \citet{donti2017task} provided an open-source repository to reproduce the \dqp{} experiments; the repository was specifically designed for inventory, power scheduling, and battery storage problems. \citet{wilder2019melding} provided code for \qptl{} for budget allocation, bipartite matching, and diverse recommendation. \mipaal{} \cite{ferber2020mipaal} relied on \qpth{} and used \cplex{} to generate cutting planes, but open-source code is not available. \citet{mandi2020interior} released \intopt{} code for knapsack, shortest path, and power scheduling. 

\citet{berthet2020learning} contributed the \tensorflow-based \dpo{} and \pfyl{} implementation, which provided universal functions for end-to-end predict-then-optimize but required users to create additional helper methods for tensor operations. \citet{elmachtoub2021smart} provided an implementation of \spo{} in Julia, which contained the shortest path and portfolio optimization problems, while \citet{mandi2020smart} implemented \spo{} with Python for the knapsack and power scheduling problems. The code for \dbb{} \cite{poganvcic2019differentiation} is applied to the shortest path, the travel salesperson, the ranking, perfect matching, and graph matching are available. Additionally, \nce{} \cite{mulamba2020contrastive} covers the knapsack, energy-cost scheduling, and bipartite matching, while \ltr{} \cite{mandi2022decision} covers the shortest path, energy-cost scheduling, and bipartite matching problems.

Except for \dpo{} and \pfyl{}, the above contributions provided solutions to specific optimization problems, and \citet{berthet2020learning} did not package \dpo{} and \pfyl{} as a generic library. In conclusion, they were confined to research-grade code aimed at reproducing experimental results.

%\cvxpylayers{} \cite{agrawal2019differentiable} is the first generic end-to-end predict-then-optimize learning framework\footnote{\url{https://github.com/cvxgrp/cvxpylayers}}. In contrast to the above codes, it requires modeling with a domain-specific language \cvxpy{}, which is embedded into a differentiable layer in a straightforward way. The emergence of cvxpylayers{} provides a more powerful tool for academia and industry. However, the solver of \cvxpylayers{} cannot compete with commercial solvers on efficiency, especially for linear and integer linear programming. Since end-to-end training requires repeated optimization in each iteration, the inefficiency of the solver becomes a bottleneck. In addition, the nature of \cvxpylayers{} means that it cannot support training with integer variables, which limits applicability to many real-world decision-making problems.

\subsubsection{Software Package}
%\noindent\textbf{Software packages.} 
\cvxpylayers{} \cite{agrawal2019differentiable} is the first generic end-to-end predict-then-optimize learning framework. Unlike the aforementioned codes, it requires modeling with a domain-specific language \cvxpy{}, which is embedded into a differentiable layer in a straightforward manner. The emergence of \cvxpylayers{} provides a more powerful tool for academia and industry. However, the solver of \cvxpylayers{} cannot compete with commercial solvers on efficiency, especially for linear programming. Since end-to-end training requires repeated optimization in each iteration, the inefficiency of the solver becomes a bottleneck. In addition, the nature of \cvxpylayers{} means that it cannot support training with integer variables, which limits its applicability to many real-world decision-making problems.

\begin{table}[ht]
    \begin{center}
    \small
    \resizebox{1.0\textwidth}{!}{
    \begin{tabular}{l|l}
        \toprule
        Method                                                & Computation per Gradient                                                                              \\ \hline
        \dqp{} \cite{donti2017task}                     & GPU-based interior-point method for quadratic programming                                 \\
        \qptl{} \cite{wilder2019melding}                & GPU-based interior-point method for quadratic programming                                 \\
        \mipaal{} \cite{ferber2020mipaal}               & Cutting-plane method + GPU-based interior-point method for quadratic programming    \\
        \intopt{} \cite{mandi2020interior}              & GPU-based interior-point method for quadratic programming                                 \\
        \cvxpylayers{} \cite{agrawal2019differentiable} & GPU-based interior-point method for conic programming                                     \\
        \spo{} \cite{elmachtoub2021smart}               & Linear/integer programming                                                                \\
        \sporel{} \cite{mandi2020smart}                 & Linear programming                                                                        \\
        \spows{} \cite{mandi2020smart}                  & Integer programming with warm starting                                                   \\
        \dbb{} \cite{poganvcic2019differentiation}      & Two Linear/integer programming                                                          \\ 
        \dpo{} \cite{berthet2020learning}              & Sampling with multiple linear/integer programming solves with random noise             \\
        \pfyl{} \cite{berthet2020learning}              & Sampling with multiple linear/integer programming solves with random noise             \\ 
        \nce{} \cite{mulamba2020contrastive}              & Linear/integer programming            \\
        \ltr{} \cite{mandi2022decision}              & Linear/integer programming           \\ \hline
        \end{tabular}}
        \end{center}
    \caption{Computational cost per gradient calculation for different methodologies.} \label{tab:cost}
\end{table}

\section{Preliminaries}
\label{sec:prel}

\subsection{Definitions and Notation}
\label{subsec:def}

For the sake of convenience, we define the following linear programming problem without loss of generality, where the decision variables are $\bm{w} \in \mathbb{R}^d$ and all $w_i \geq 0$, the objective function coefficients (cost vector) are $\bm{c} \in \mathbb{R}^d$, the constraint coefficients are $\bm{A} \in \mathbb{R}^{k \times d}$, and the right-hand sides of the constraints are $\bm{b} \in \mathbb{R}^k$:
\begin{equation}
\begin{aligned}
\min_{\bm{w}} \quad & \bm{c}^\intercal \bm{w} \\
\textrm{s.t.} \quad & \bm{A}\bm{w} \leq \bm{b} \\
& \bm{w} \geq \bm{0} \\
\end{aligned}
\end{equation}

When some variables $w_i$ are restricted to integers, we obtain a (mixed) integer program:

\begin{equation}
\begin{aligned}
\min_{\bm{w}} \quad & \bm{c}^\intercal \bm{w} \\
\textrm{s.t.} \quad & \bm{A}\bm{w} \leq \bm{b} \\
& \bm{w} \geq \bm{0} \\
& w_i \in \mathbb{Z} \quad \forall i \in I \subseteq \{1, 2, ..., d\}\\
\end{aligned}
\end{equation}

For linear and integer programming, let $S$ be the feasible region, $z^*(\bm{c})$ be the optimal objective value with respect to the objective function coefficients $\bm{c}$, and $\bm{w}^*(\bm{c}) \in W^*(\bm{c})$ be a particular optimal solution derived from some solver. We define the optimal solution set $W^*(\bm{c})$ because there may be multiple optima.

As mentioned above, some coefficients are unknown and need to be predicted prior to optimization. Here, we assume that only the objective function coefficients $\bm{c}^i$ are unknown, but they correlate with a feature vector $\bm{x}^i \in \mathbb{R}^p$. Given a training dataset $\mathcal{D} = \{(\bm{x}^1, \bm{c}^1), (\bm{x}^2, \bm{c}^2), ..., (\bm{x}^n, \bm{c}^n)\}$ or $\mathcal{D} = \{(\bm{x}^1, \bm{w}^*(\bm{c}^1)), (\bm{x}^2, \bm{w}^*(\bm{c}^2)), ..., (\bm{x}^n, \bm{w}^*(\bm{c}^n)\})$, a machine learning predictor $g(\cdot)$ can be trained to minimize a loss function $\mathcal{L}(\cdot)$, where $\bm{\theta}$ are predictor parameters and $\hat{\bm{c}}^i = g(\bm{x}^i; \bm{\theta})$ is the prediction of the objective function coefficients $\bm{c}^i$.

\subsection{The Two-Stage Method}
\label{subsec:2s}

\begin{figure}[htbp]
    \centering
    \includegraphics[width=1\textwidth]{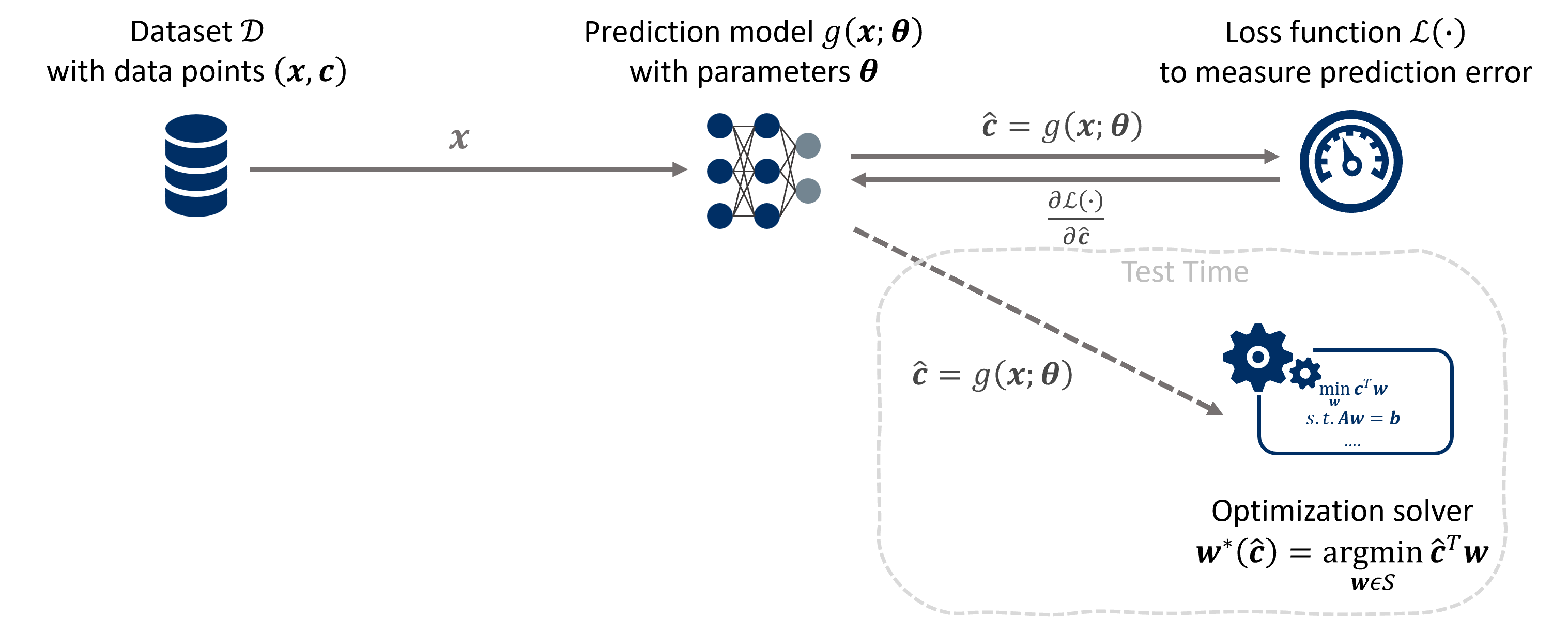}
    \caption{Illustration of the two-stage predict-then-optimize framework: A labeled dataset $\mathcal{D}$ of $(\bm{x},\bm{c})$ pairs is used to fit a machine learning predictor that minimizes prediction error. At test time (grey box), the predictor is used to estimate the parameters of an optimization problem, which is then tackled with an optimization solver. The two stages are thus completely separate.}
    \label{fig:2s}
\end{figure}

As Figure \ref{fig:2s} shows, the two-stage approach trains a predictor $g(\cdot)$ by minimizing a loss function with respect to the true coefficients $\bm{c}$, such as mean squared error (MSE), $\mathcal{L}_{\text{MSE}}(\hat{\bm{c}}, \bm{c}) = \frac{1}{n} \norm{\hat{\bm{c}} - \bm{c}}^2$. Following training and given an instance with feature vector $\bm{x}$, the predictor produces coefficients $\hat{\bm{c}}=g(\bm{x}; \bm{\theta})$, which is then used to solve the optimization problem. It decomposes the predict-then-optimize problem into a traditional regression and then an optimization. However, it introduces a mismatch between the decision error and the prediction error.

\subsection{Gradient-based End-to-end Predict-then-Optimize}
\label{subsec:e2e}

\begin{figure}[htbp]
    \centering
    \includegraphics[width=1\textwidth]{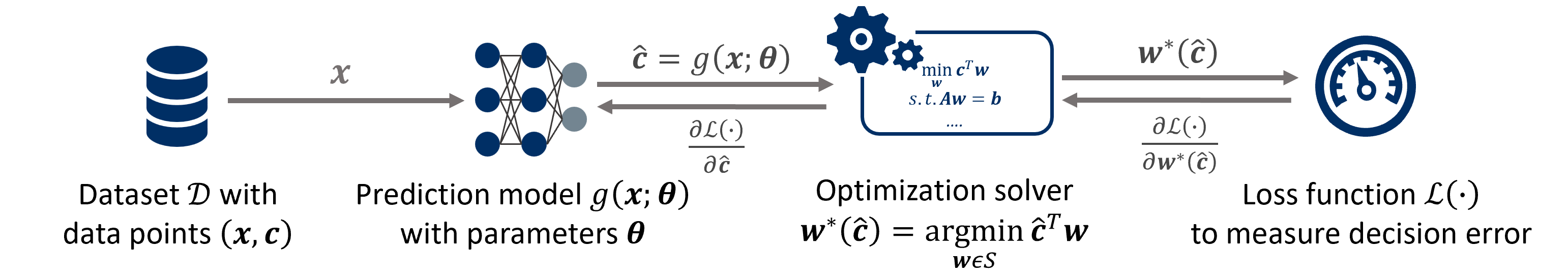}    
    \caption{Illustration of the end-to-end predict-then-optimize framework: A labeled dataset $\mathcal{D}$ of $(\bm{x},\bm{c})$ (or $(\bm{x},\bm{w}^*(\bm{c}))$) pairs is used to fit a machine learning predictor that directly minimizes decision error. The critical component is an optimization solver, which is embedded into a differentiable predictor (e.g., a neural network). At test time, this approach is similar to the two-stage approach from Figure~\ref{fig:2s}.}
    \label{fig:e2e}
\end{figure}

The main drawbacks of the two-stage approach are that it does not account for decision error during training and that it always requires the true coefficients as labels for supervised learning. In contrast, the end-to-end predict-then-optimize method depicted in Figure \ref{fig:e2e} aims to minimize the decision error and has the potential to learn with only optimal solutions. In line with deep learning terminology, we will use the term ``backward pass" to refer to the gradient computation and parameter updating via the backpropagation algorithm. In order to incorporate optimization into the prediction, we can derive the derivative of the optimization task and then apply the gradient descent algorithm, Algorithm~\ref{alg:grad}, to update the predictor parameters.

\begin{algorithm}
  \caption{End-to-end Gradient Descent}\label{alg:grad}
  \small
  \begin{algorithmic}[1]
    \Require coefficient matrix $\bm{A}$, right-hand side $\bm{b}$, data $\mathcal{D}$
    \State Initialize predictor parameters $\bm{\theta}$ for predictor $g(\bm{x}; \bm{\theta})$
    \For{epochs}
      \For{each batch of training data $(\bm{x}, \bm{c})$ or $(\bm{x}, \bm{w}^* (\bm{c}))$}
        \State Sample batch of the $\bm{c}$ or $\bm{w}^* (\bm{c})$ with the corresponding features $\bm{x}$
%        \State Predict coefficients $\hat{\bm{c}}$ using predictor $g(\bm{x}; \bm{\theta})$
        \State Predict coefficients using predictor $\hat{\bm{c}} := g(\bm{x}; \bm{\theta})$
        \State Forward pass to compute optimal solution $\bm{w}^*(\hat{\bm{c}}) :=  {\argmin}_{\bm{w} \in S} \hat{\bm{c}}^\intercal \bm{w}$ 
        \State Forward pass to compute decision loss $\mathcal{L}(\hat{\bm{c}}, \cdot)$ 
        \State Backward pass from loss $\mathcal{L}(\hat{\bm{c}}, \cdot)$  to update parameters $\bm{\theta}$ with gradient
      \EndFor
    \EndFor
  \end{algorithmic}
\end{algorithm}

For an appropriately defined loss function that penalizes decision error, the chain rule can be used to calculate the gradient of the loss with respect to the predictor parameters as follows:

\begin{align}
\frac{\partial \mathcal{L}(\hat{\bm{c}}, \cdot)}{\partial \bm{\theta}}
& =
\frac{\partial \mathcal{L}(\hat{\bm{c}}, \cdot)}{\partial \hat{\bm{c}}}
\frac{\partial \hat{\bm{c}}}{\partial \bm{\theta}} \label{eq:rule1} \\ 
& = 
\frac{\partial \mathcal{L}(\hat{\bm{c}}, \cdot)}{\partial \bm{w}^*(\hat{\bm{c}})}
\frac{\partial \bm{w}^*(\hat{\bm{c}})}{\partial \hat{\bm{c}}}
\frac{\partial \hat{\bm{c}}}{\partial \bm{\theta}} \label{eq:rule2} \\ 
\textrm{Note: } & \frac{\partial \hat{\bm{c}}}{\partial \bm{\theta}} =
\frac{\partial g(\bm{x}; \bm{\theta})}{\partial \bm{\theta}}
\end{align}

The last term $\frac{\partial \hat{\bm{c}}}{\partial \bm{\theta}}$ is the gradient of the predictions with respect to the model parameters, which is trivial to calculate in modern deep learning frameworks. The challenging part is to compute the a surrogate decision loss function $\frac{\partial \mathcal{L}(\hat{\bm{c}}, \cdot)}{\partial \hat{\bm{c}}}$ or constrained optimization gradient $\frac{\partial \bm{w}^*(\hat{\bm{c}})}{\partial \hat{\bm{c}}}$. Since the mapping from coefficients $\bm{c}$ to solution vector $\bm{w}^*$ is piecewise constant for linear and integer programming, the predictor parameters cannot be updated with gradient descent. Thus, both \spo{} and \pfyl{} are surrogate loss functions that adhere to Equation~\ref{eq:rule1}, which derive the gradient of decision loss $\frac{\partial \mathcal{L}(\hat{\bm{c}}, \cdot)}{\partial \hat{\bm{c}}}$, while \dbb{} and \dpo{}, follow Equation~\ref{eq:rule2}, approximate $\frac{\partial \bm{w}^*(\hat{\bm{c}})}{\partial \hat{\bm{c}}}$. Furthermore, \dbb{} and \dpo{} provide flexibility to customize loss functions, and the use of loss functions in our experiments is detailed in Table~\ref{tab:model} and Table~\ref{tab:cnn}.

\subsubsection{Decision Loss}
\label{subsec:loss}

The concept of regret (also known as SPO loss \cite{elmachtoub2021smart}) has been introduced to quantify the error in decision-making. It is defined as the gap of the objective value between the true optimal solution $\bm{w}^*(\bm{c})$ and the optimal solution acquired using the predicted coefficients $\bm{w}^*(\hat{\bm{c}})$. Formally, it is represented as: 
\begin{equation}
    \mathcal{L}_{\text{Regret}}(\hat{\bm{c}}, \bm{c}) = \bm{c}^\intercal \bm{w}^*(\hat{\bm{c}}) - z^*(\bm{c}).
    \label{eq:regret}
\end{equation}

\begin{figure}[htb]
    \centering
    \subfloat[]{\includegraphics[width=0.48\textwidth]{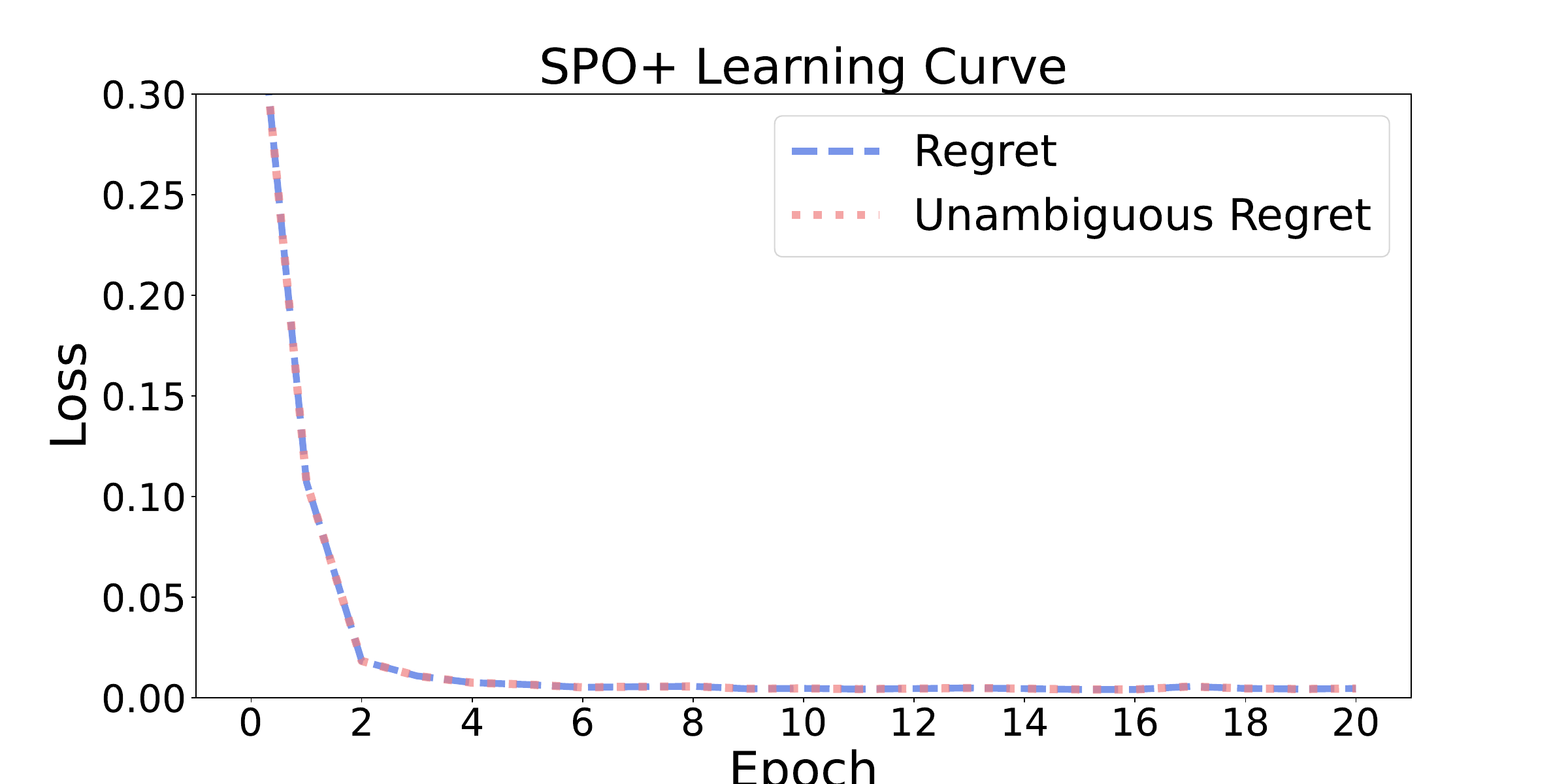}} 
    \subfloat[]{\includegraphics[width=0.48\textwidth]{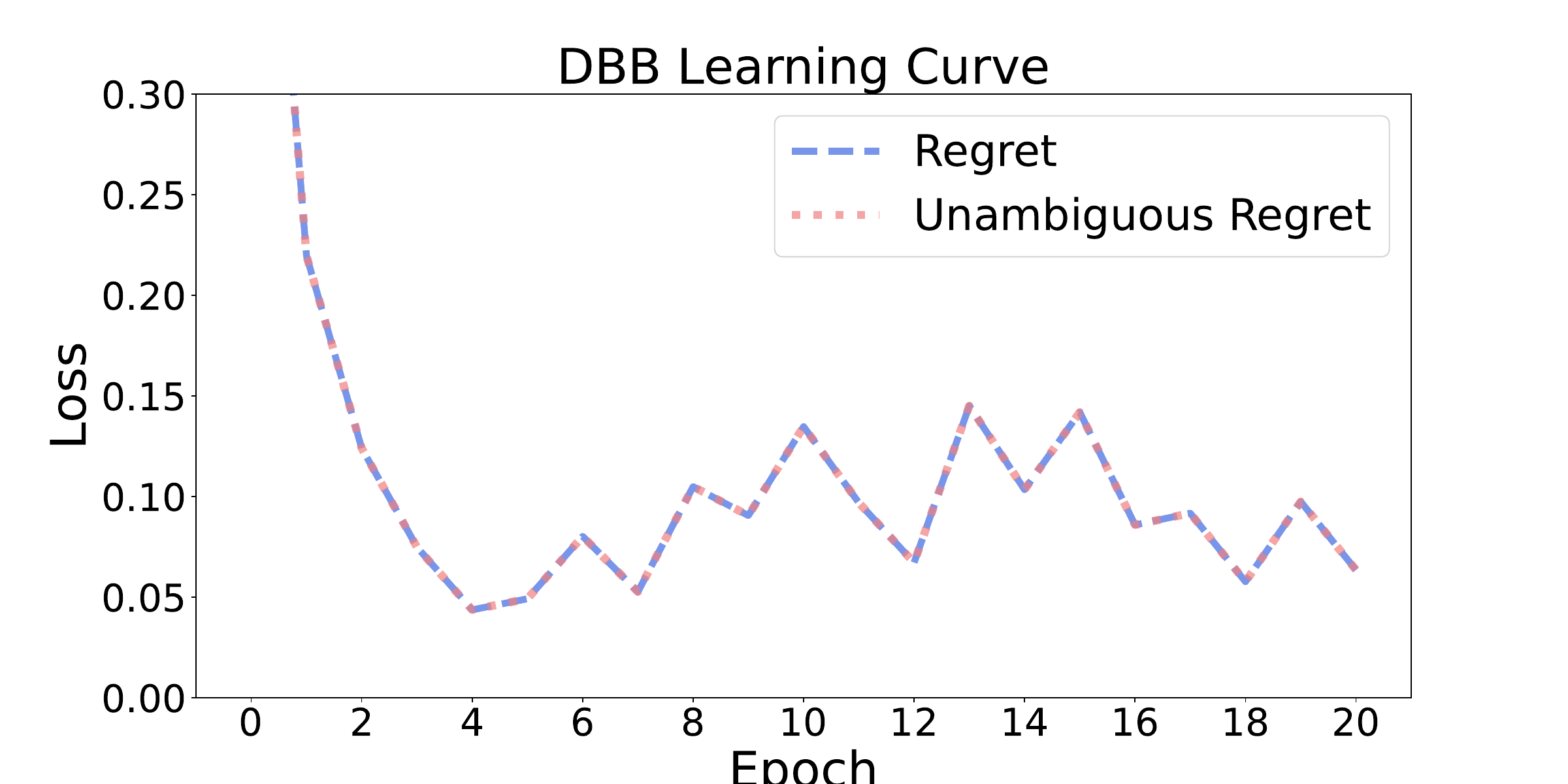}} \\
    \caption{As shown for the learning curves of the training of \spo{} (Left) and \dbb{} (Right) on the shortest path, regret and unambiguous regret in the various tasks overlap almost exactly.}
    \label{fig:rcurve}
\end{figure}

Given the coefficients of the objective function $\hat{\bm{c}}$, there may be multiple optimal solutions to ${\min}_{\bm{w} \in S} \hat{\bm{c}}^\intercal \bm{w}$. To address this, \citet{elmachtoub2021smart} devised the ``unambiguous" regret (also called an unambiguous SPO loss): $\mathcal{L}_{\text{URegret}}(\hat{\bm{c}}, \bm{c}) = {\max}_{\bm{w} \in W^*(\bm{c})} \bm{w}^\intercal \bm{c} - z^*(\bm{c})$. This loss considers the worst case among all optimal solutions with respect to the predicted coefficients. 

An example of the worst-case scenario is when the predicted coefficients $\hat{\bm{c}}$ converge to all zeros, making all feasible solutions optimal. Consequently, there is a small likelihood that the solver always selects the particular optimal solution $\bm{w}^*(\hat{\bm{c}})$ equivalent to the true optimal solution $\bm{w}^*(\bm{c})$, leading to zero loss with zero prediction. 

\pyepo{} provides an evaluation module (Section \ref{subsec:metr}) that includes both regret and unambiguous regret. However, as Figure \ref{fig:rcurve} shows, regret and unambiguous regret are almost identical throughout the training procedures. Therefore, although unambiguous regret is more theoretically rigorous, it is not necessary to consider it in practice.

In addition to regret, the decision error can also be defined as the difference between the true solution and its prediction, such as the Hamming distance of the solutions \citep{poganvcic2019differentiation} and the squared error of the solutions \citep{berthet2020learning}. Moreover, \citet{dalle2022learning} also considered treating the objective value $\bm{c}^\intercal \bm{w}^*_{\hat{\bm{c}}}$ itself as a loss.

\subsection{Methodologies}
\label{subsec:mthd}

\subsubsection{Smart Predict-then-Optimize~\cite{elmachtoub2021smart}}

To make the decision error differentiable, \citet{elmachtoub2021smart} proposed \spo{}, a convex upper bound on the regret:

\begin{equation}
\begin{aligned}
\mathcal{L}_{\text{SPO+}}(\hat{\bm{c}}, \bm{c}) = - \underset{\bm{w} \in S}{\min} \{(2 \hat{\bm{c}} - \bm{c})^\intercal \bm{w}\} + 2 \hat{\bm{c}}^\intercal \bm{w}^* (\bm{c}) - z^*(\bm{c}).
\end{aligned}
\end{equation}

One proposed subgradient for this loss is as follows:

\begin{equation}
\begin{aligned}
2 \bm{w}^*(\bm{c}) - 2 \bm{w}^* (2 \hat{\bm{c}} - \bm{c}) \in \frac{\partial \mathcal{L}_{\text{SPO+}}(\hat{\bm{c}}, \bm{c})}{\partial \hat{\bm{c}}}
\end{aligned}
\end{equation}

Thus, we can use Algorithm \ref{alg:grad} to directly minimize $\mathcal{L}_{\text{SPO+}}(\hat{\bm{c}}, \bm{c})$ with gradient descent. This algorithm with \spo{} requires solving ${\min}_{\bm{w} \in S} (2 \hat{\bm{c}} - \bm{c})^\intercal \bm{w}$ for each training iteration. 

To accelerate \spo{} training, \citet{mandi2020smart} employed relaxations (\sporel{}) and warm starting (\spows{}) to speed up optimization. The idea of \sporel{} is to use the continuous relaxation of the integer program during training. This simplification greatly reduces training time at the expense of model performance. Compared to \spo{}, the improvement in \sporel{} in training efficiency is not negligible. For example, linear programming can be solved in polynomial time, while integer programming is worst-case exponential. In Section \ref{sec:eval}, we will further discuss this performance-efficiency trade-off. For \spows{}, \citet{mandi2020smart} suggested using previous solutions as a starting point for the integer programming solver, which potentially improves efficiency by narrowing down the search space.

\subsubsection{Differentiable Black-Box Solver~\cite{poganvcic2019differentiation}}

\dbb{} was developed by \citet{poganvcic2019differentiation} to estimate gradients from finite difference, replacing the zero gradients in $\frac{\partial \bm{w}^*(\hat{\bm{c}})}{\partial \hat{\bm{c}}}$. For the predicted coefficients $\hat{\bm{c}}$, \citet{poganvcic2019differentiation} performed a linear interpolation of the loss function $\mathcal{L} (\hat{\bm{c}},\cdot)$ between $\hat{\bm{c}}$ and $\hat{\bm{c}} + \lambda \frac{\partial \mathcal{L} (\hat{\bm{c}},\cdot)}{\partial \bm{w}^*(\hat{\bm{c}})}$. Thus, the substitute function becomes piecewise affine. Therefore, when computing $\bm{w}^*(\hat{\bm{c}})$, a useful nonzero gradient is obtained at the cost of faithfulness. The forward and backward passes are shown in Algorithm \ref{algo:dbbf} and Algorithm \ref{algo:dbbb}. The hyperparameter $\lambda > 0$ controls the degree of interpolation. However, compared to \spo{}, the approximation function of \dbb{} is nonconvex, so the convergence to a global optimum is compromised, even when the predictor is convex in its parameters.

\begin{minipage}{0.46\textwidth}
\begin{algorithm}[H]
    \centering
    \caption{DBB Forward Pass}\label{algo:dbbf}
    \begin{algorithmic}[1]
        \Require $\hat{\bm{c}}$
        \State Solve $\bm{w}^*(\hat{\bm{c}})$
        \State Save $\hat{\bm{c}}$ and $\bm{w}^*(\hat{\bm{c}})$ for backward pass \\
        \Return $\bm{w}^*(\hat{\bm{c}})$
    \end{algorithmic}
\end{algorithm}
\end{minipage}
\hfill
\begin{minipage}{0.46\textwidth}
\begin{algorithm}[H]
    \centering
    \caption{DBB Backward Pass}\label{algo:dbbb}
    \begin{algorithmic}[1]
        \Require $\frac{\partial \mathcal{L}(\hat{\bm{c}}, \cdot)}{\partial \bm{w}^*(\hat{\bm{c}})}$, $\lambda$
        \State Load $\hat{\bm{c}}$ and $\bm{w}^*(\hat{\bm{c}})$ from forward pass
        \State $\bm{c}^\prime := \hat{\bm{c}} + \lambda \frac{\partial \mathcal{L}(\hat{\bm{c}}, \cdot)}{\partial \bm{w}^* (\hat{\bm{c}})}$
        \State Solve $\bm{w}^*(\bm{c}^\prime)$ \\
        \Return $\frac{\partial \mathcal{L}(\hat{\bm{c}}, \cdot)}{\partial \hat{\bm{c}}} := \frac{1}{\lambda} (\bm{w}^*(\bm{c}^\prime) - \bm{w}^*(\hat{\bm{c}}))$ 
    \end{algorithmic}
\end{algorithm}
\end{minipage}
\vspace{.1in}

Similarly to \spo{}, \dbb{} requires solving the optimization problem in each training iteration. Thus, utilizing a relaxation approach may also work for \dbb{}. However, \citet{poganvcic2019differentiation} did not consider this option. Given the potential efficiency gains that continuous relaxation can bring, we also conducted experiments for \dbbrel{} in section \ref{sec:eval}. The same applies to \dpo{} and \pfyl{} below.

\subsubsection{Differentiable Perturbed Optimizer~\cite{berthet2020learning}}

\begin{figure}[htbp]
    \centering
    \includegraphics[width=0.8\textwidth]{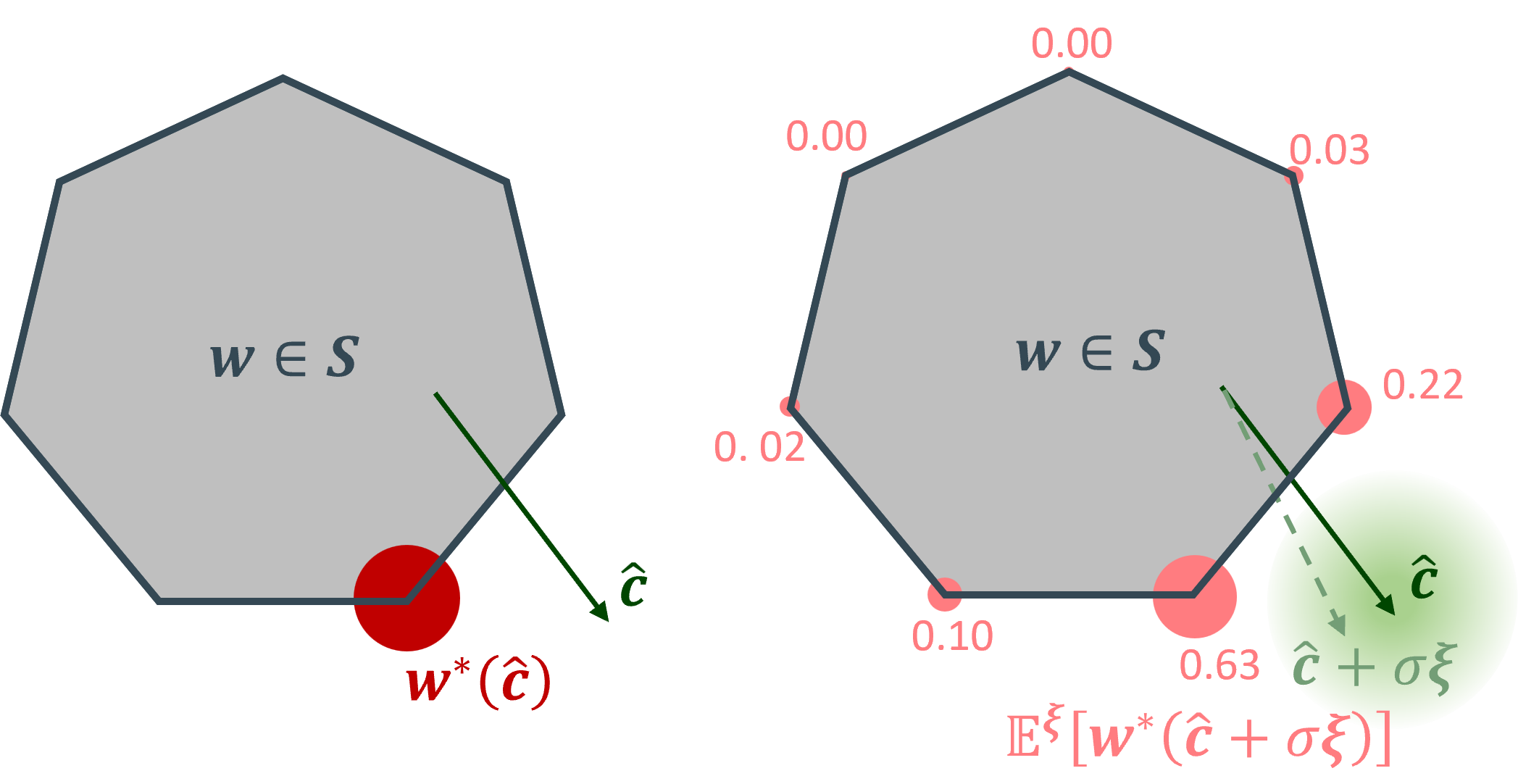}    
    \caption{Illustration of optimal solution (left) and expected optimal solution under perturbation (right): While the optimal solution $\bm{w}^*(\hat{\bm{c}})$ remains static or shifts to a neighboring extreme point as the predicted coefficient $\hat{\bm{c}}$ changes, the perturbed expectation undergoes a smooth variation.}
    \label{fig:ptb}
\end{figure}

\dpo{} \citep{berthet2020learning} uses Monte-Carlo samples to estimate solutions, in which the predicted coefficients $\hat{\bm{c}}$ are perturbed by Gaussian noise $\bm{\xi} \sim \mathcal{N}(\bm{0}, \bm{I})$. Based on the $K$ samples of perturbed coefficients $\hat{\bm{c}} + \sigma \bm{\xi}$, it returns the estimate of the optimal solution expectation $$\mathbb{E}_{\bm{\xi}} [\bm{w}^*(\hat{\bm{c}} + \sigma \bm{\xi})] \approx \frac{1}{K} 
 \sum_{\kappa}^K {\bm{w}^*(\hat{\bm{c}} + \sigma \bm{\xi}_{\kappa})}.$$ \dpo{} outputs the expectation of the optimal solution $\mathbb{E}_{\bm{\xi}}[\bm{w}^*(\hat{\bm{c}} + \sigma \bm{\xi})]$ under a perturbed coefficients with random noise $\hat{\bm{c}} + \sigma \bm{\xi}$, rather than the optimal solution $\bm{w}^*(\hat{\bm{c}})$ under a fixed coefficients $\hat{\bm{c}}.$ As visualized in Fig~\ref{fig:ptb}, the expectation can be regarded as a convex combination of different feasible solutions.

Thus, unlike the piecewise constant function $\bm{w}^*(\hat{\bm{c}})$, $\mathbb{E}_{\bm{\xi}}[\bm{w}^*(\hat{\bm{c}} + \sigma \bm{\xi})]$ varies the proportions in response to the change of $\hat{\bm{c}}$, providing a nonzero Jacobian matrix of $\hat{\bm{c}}$: $$\frac{\partial \mathbb{E}_{\bm{\xi}}[\bm{w}^*(\hat{\bm{c}} + \sigma \bm{\xi})]}{\partial \hat{\bm{c}}} \approx \frac{1}{K} \sum_{\kappa}^K {\bm{w}^*(\hat{\bm{c}} + \sigma \bm{\xi}_{\kappa})} \bm{\xi}_{\kappa}^\intercal.$$  The forward and backward pass are as follows:

\begin{minipage}{0.46\textwidth}
\begin{algorithm}[H]
    \centering
    \caption{DPO Forward Pass}\label{algo:dpof}
    \begin{algorithmic}[1]
        \Require $\hat{\bm{c}}$, $K$, $\sigma$
        \For{sample $\kappa \in \{1, ..., K\}$}
            \State Generate Gaussian noise $\bm{\xi}_{\kappa}$
            \State Solve: $\bm{w}_{\kappa}^{\bm{\xi}} := \bm{w}^*(\hat{\bm{c}} + \sigma \bm{\xi}_{\kappa})$
            \State Save $\bm{w}_{\kappa}^{\bm{\xi}}$ and $\bm{\xi}_{\kappa}$ for backward pass
        \EndFor \\
        \Return $\frac{1}{K} \sum_{\kappa}^K {\bm{w}_{\kappa}^{\bm{\xi}}}$
    \end{algorithmic}
\end{algorithm}
\end{minipage}
\hfill
\begin{minipage}{0.46\textwidth}
\begin{algorithm}[H]
    \centering
    \caption{DPO Backward Pass}\label{algo:dpob}
    \begin{algorithmic}[1]
        \Require $\frac{\partial \mathcal{L}(\hat{\bm{c}},\cdot)}{\partial \mathbb{E}_{\bm{\xi}}[\bm{w}^*]}$, $K$
        \State Load $\bm{w}_{\kappa}^{\bm{\xi}}$ and $\bm{\xi}_{\kappa}$ from forward pass
        \State Compute $\frac{\partial \mathbb{E}_{\bm{\xi}}[\bm{w}^*]}{\partial \hat{\bm{c}}} := \frac{1}{K} \sum_{\kappa}^K {\bm{w}_{\kappa}^{\bm{\xi}}} \bm{\xi}_{\kappa}^\intercal$ 
        \State Compute $\frac{\mathcal{L}(\hat{\bm{c}},\cdot)}{\partial \hat{\bm{c}}} := \frac{\partial \mathcal{L}(\hat{\bm{c}},\cdot)}{\partial \mathbb{E}_{\bm{\xi}}[\bm{w}^*]} \frac{\partial \mathbb{E}_{\bm{\xi}}[\bm{w}^*]}{\partial \hat{\bm{c}}}$\\
        \Return $\frac{\mathcal{L}(\hat{\bm{c}},\cdot)}{\partial \hat{\bm{c}}}$
    \end{algorithmic}
\end{algorithm}
\end{minipage}

\subsection{Perturbed Fenchel-Young Loss~\cite{berthet2020learning}}

Instead of using an arbitrary loss for \dpo{}, \citet{berthet2020learning} further constructed the Fenchel-Young loss \citep{blondel2020learning} to directly compute the decision error $\mathcal{L}_{\text{FY}}(\hat{\bm{c}}, \bm{w}^*(\bm{c}))$ and its gradient $\frac{\partial \mathcal{L}_{\text{FY}}(\hat{\bm{c}}, \bm{w}^*(\bm{c}))}{\partial \hat{\bm{c}}}$. Compared to \dpo{}, \pfyl{} avoids the inefficient calculation of the Jacobian matrix $\nabla \bm{w}^*(\hat{\bm{c}})$ and includes a theoretically sound loss function. 

The loss of \pfyl{} is based on the Fenchel duality: The expectation of the perturbed minimizer is defined as $F(\bm{c}) =  \mathbb{E}_{\bm{\xi}}[\underset{\bm{w} \in S}{\min} {\{(\bm{c}+\sigma \bm{\xi})^\intercal \bm{w}\}}]$, and the dual of $F(\bm{c})$, denoted by $\Omega (\bm{w}^*({\bm{c}})),$ is utilized to define the Fenchel-Young loss to reduce duality gap:
$$\mathcal{L}_{\text{FY}}(\hat{\bm{c}}, \bm{w}^*({\bm{c}})) =  \hat{\bm{c}}^\intercal \bm{w}^*({\bm{c}}) - F(\hat{\bm{c}}) - \Omega (\bm{w}^*({\bm{c}})),$$
where the gradient of the loss is 
$$\frac{\partial \mathcal{L}_{\text{FY}}(\hat{\bm{c}}, \bm{w}^*({\bm{c}}))}{\partial \hat{\bm{c}}} = \bm{w}^*({\bm{c}}) - \frac{\partial F(\hat{\bm{c}})}{\partial \hat{\bm{c}}} = \bm{w}^*({\bm{c}}) - \mathbb{E}_{\bm{\xi}} [\bm{w}^* (\hat{\bm{c}}+\sigma \bm{\xi})].$$
Similar to \dpo{}, we can estimate the well-defined gradient as
$$\frac{\partial \mathcal{L}_{\text{FY}}(\hat{\bm{c}}, \bm{w}^*({\bm{c}}))}{\partial \hat{\bm{c}}} \approx \bm{w}^*({\bm{c}}) - \frac{1}{K} 
 \sum_{\kappa}^K \bm{w}^* (\hat{\bm{c}}+\sigma \bm{\xi})$$

\section{Implementation and Modeling}
\label{sec:tool}

The core module of \pyepo{} is an ``autograd" function inherited from \pytorch{} \cite{paszke2017automatic}. Such functions perform a forward pass to yield optimal solutions or decision losses and a backward pass to obtain nonzero gradients so that the prediction model can learn from the decision error. Thus, our implementation extends \pytorch{}, which facilitates the deployment of end-to-end predict-then-optimize tasks using any neural network that can be built with \pytorch{}.

\subsection{Modeling Language}

We choose \gurobipy{} \cite{gurobi}, \coptpy{} \cite{copt} or \pyomo{} \cite{hart2017pyomo} to build optimization models. All of these are algebraic modeling languages (AMLs) written in Python. 

\gurobipy{} and \coptpy{} are Python interfaces to \gurobi{} and \copt{}, which combine the expressiveness of a modeling language with the flexibility of a programming language. As an official interface of \gurobi{}, \gurobipy{} has a simple algebraic syntax and natively supports all features of \gurobi{} for modeling, and \coptpy{} offers similar functionality. Recognizing that users might not have access to a \gurobi{} or \copt{} license, we also implemented a \pyomo{} interface as an alternative. \pyomo{} is an open-source optimization modeling language that supports a variety of solvers, including \gurobi{} and \glpk{}. 

Each of these tools provides a natural way to express mathematical programming models. Therefore, users without specialized optimization knowledge can easily build and maintain optimization models through high-level algebraic representations. 

Moreover, \pyepo{} also allows users to construct optimization models from scratch coding with any algorithm and solver, facilitating fast and flexible model customization for both research and production purposes.

\begin{figure}[htb]
    \centering
    \subfloat[]{\includegraphics[width=0.48\textwidth]{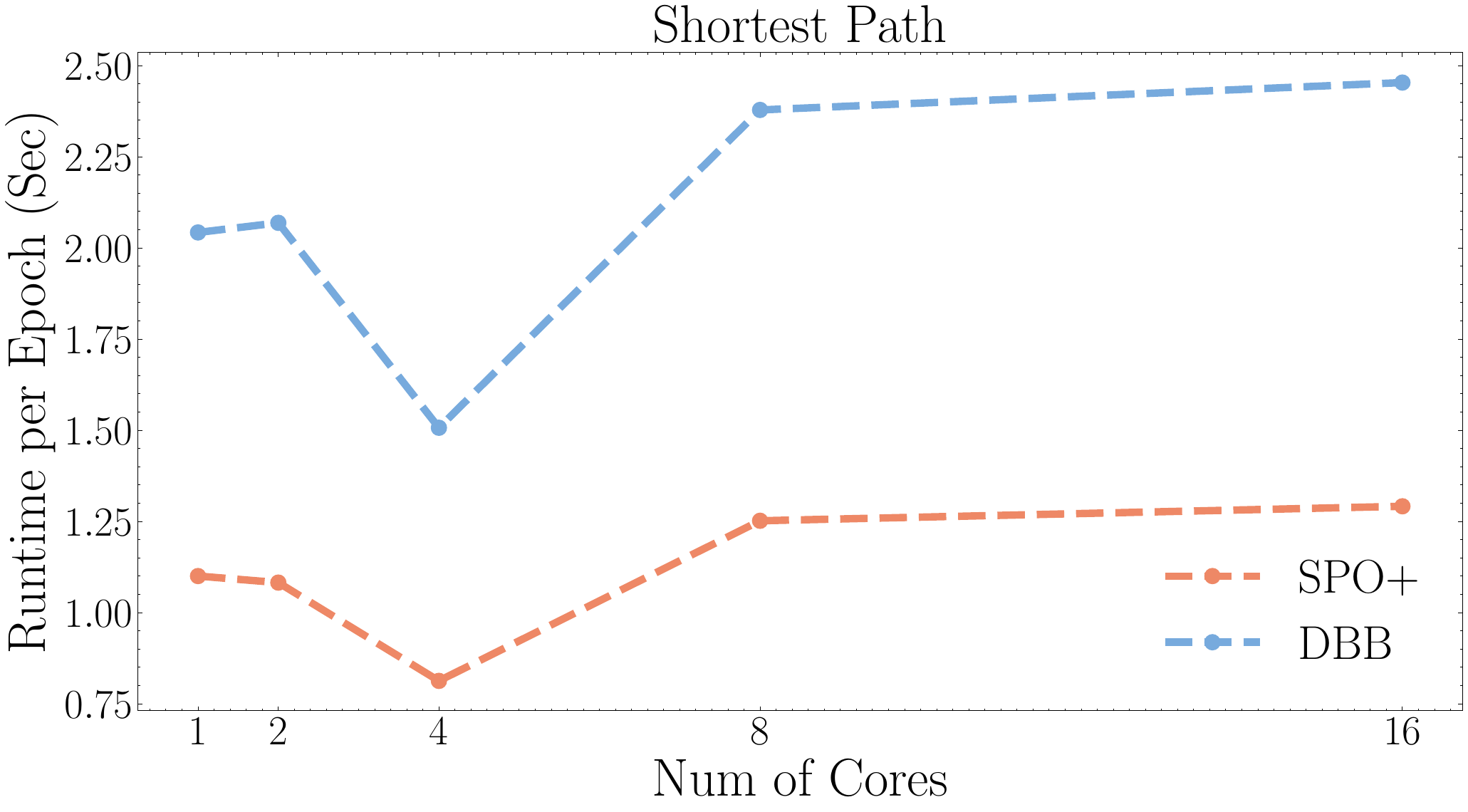}} 
    \subfloat[]{\includegraphics[width=0.48\textwidth]{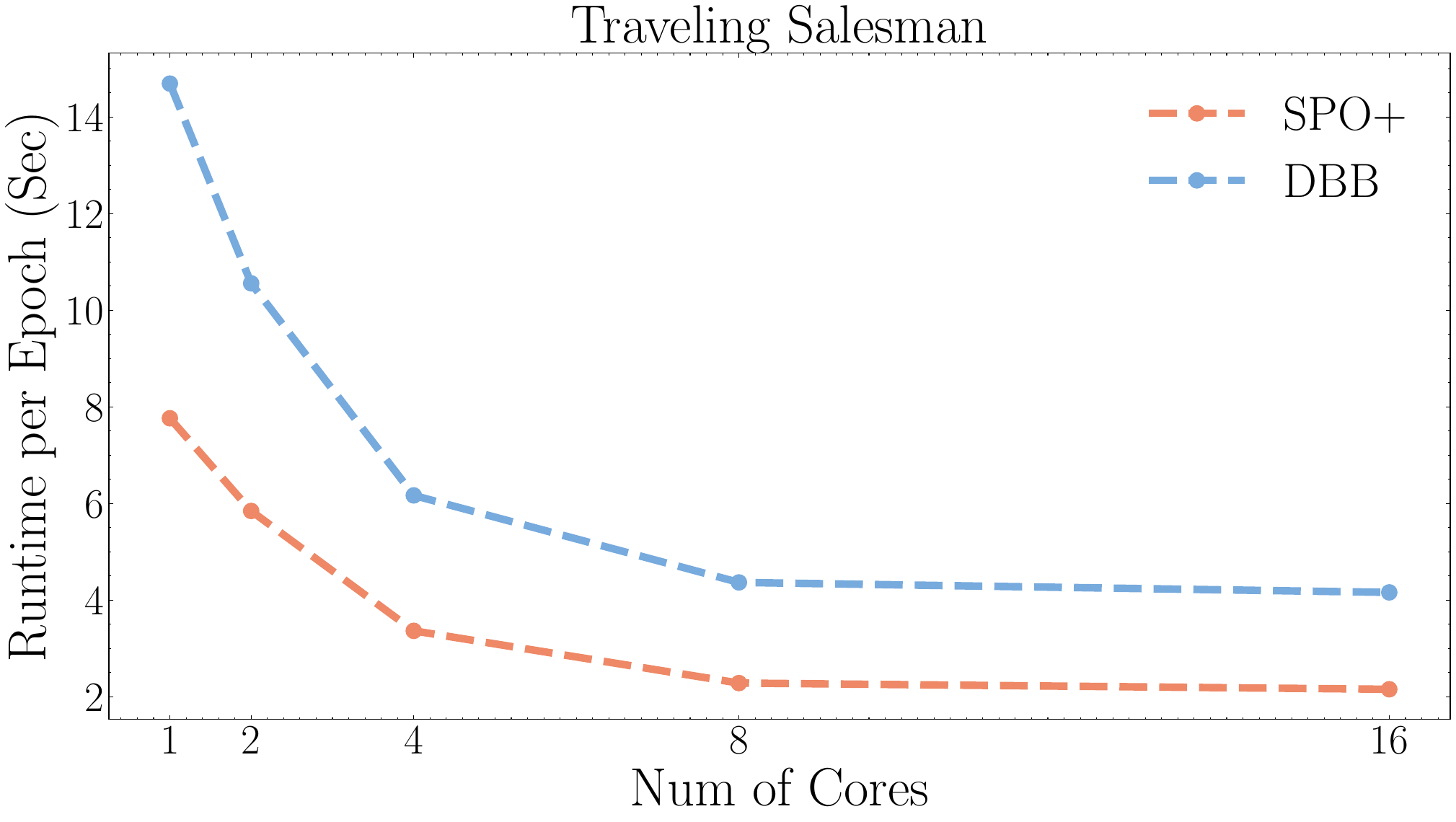}} \\
    \caption{Parallel Efficiency: Although creating a new process incurs additional overhead, parallel computing with an appropriate number of processors can effectively reduce the training time for \spo{} and \dbb{}.}
    \label{fig:parallel}
\end{figure}

\subsection{Parallel Computing}

In addition, \pyepo{} supports parallel computing. For methods such as \spo{}, \dbb{}, \dpo{}, and \pfyl{}, the computational cost of repeated optimization is the major challenge, especially for integer programs. These methods require solving a batch of optimization problems per iteration to obtain the solution and the gradient. Parallel computing helps mitigate this issue by distributing the computational load across multiple processors.

Figure \ref{fig:parallel} illustrates the average running time per epoch for a mini-batch gradient descent algorithm with a batch size of $32$, plotted against the number of cores. The decrease in running time per epoch is sublinear with respect to the number of cores. This may be due to the overhead associated with initializing additional cores, which can dominate computation costs. For example, in Figure \ref{fig:parallel}, for the shortest path, the easily solvable polynomial problem, the running time actually increases when the number of cores exceeds $4$. In contrast, the more complicated $\mathcal{NP}$-complete problem, TSP, continues to benefit slightly from additional cores. Overall, we believe this feature to be crucial for large-scale predict-then-optimize tasks.

\subsection{Optimization Model}
\label{subsec:optmodel}

The first step in using \pyepo{} is to create an optimization model that inherits from the \texttt{optModel} class. Since \pyepo{} tackles predict-then-optimize with unknown objective function coefficients, it is necessary to instantiate an optimization model, \texttt{optModel}, with fixed constraints and variable coefficients. Such an optimization model would accept different objective function coefficients and find the corresponding optimal solutions under identical constraints. The construction of \texttt{optModel} is separated from the autograd functions (see Section \ref{subsec:func}), and then the instance of \texttt{optModel} will be passed as an argument to the autograd functions.

\subsubsection{Optimization Model from Scratch}

In \pyepo{}, the \texttt{optModel} functions as a black-box, which means that it does not rely on a specific algorithm or solver. This design grants users the flexibility to customize their tasks. To construct an \texttt{optModel} from scratch, users must override abstract methods \texttt{\_getModel} to build a model and retrieve its variables, \texttt{setObj} to set the objective function with given coefficients, and \texttt{solve} to find an optimal solution. In addition, \texttt{optModel} provides an attribute \texttt{modelSense} to indicate whether the problem is minimization or maximization. The following shortest path example uses the Python library $NetworkX$ \cite{hagberg2008exploring} and its built-in Dijkstra's algorithm:

\begin{lstlisting}[language=Python]
import networkx as nx
from pyepo import EPO
from pyepo.model.opt import optModel

class myShortestPathModel(optModel):

    def __init__(self):
        self.modelSense = EPO.MINIMIZE
        self.grid = (5,5) # graph size
        self.arcs = self._getArcs() # list of arcs
        super().__init__()

    def _getModel(self):
        """
        A method to build a model
        """
        # build graph as model
        model = nx.Graph() 
        # add arcs as variables
        model.add_edges_from(self.arcs, cost=0)
        var = model.edges
        return model, var

    def setObj(self, c):
        """
        A method to set objective function
        """
        for i, e in enumerate(self.arcs):
            self._model.edges[e]["cost"] = c[i]

    def solve(self):
        """
        A method to solve model
        """
        # dijkstra
        s = 0 # source node
        t = self.grid[0] * self.grid[1] - 1 # target node
        path = nx.shortest_path(self._model, weight="cost",
                                source=s, target=t)
        # convert the path into active edges
        edges = []
        u = 0
        for v in path[1:]:
            edges.append((u,v))
            u = v
        # init sol & obj
        sol = [0] * self.num_cost
        obj = 0
        # convert active edges into solution and obj
        for i, e in enumerate(self.arcs):
            if e in edges:
                # active edge
                sol[i] = 1 
                # cost of active edge
                obj += self._model.edges[e]["cost"] 
        return sol, obj

    def _getArcs(self):
        """
        A helper method to get a list of arcs for the grid network
        """
        arcs = []
        h, w = self.grid
        for i in range(h):
            # edges on rows
            for j in range(w - 1):
                v = i * w + j
                arcs.append((v, v + 1))
            # edges in columns
            if i == h - 1:
                continue
            for j in range(w):
                v = i * w + j
                arcs.append((v, v + w))
        return arcs
        
optmodel = myShortestPathModel()
\end{lstlisting}

\subsubsection{Optimization Model with Gurobi}

On the other hand, we provide \texttt{optGrbModel} to create an optimization model with \gurobipy{}. Unlike \texttt{optModel}, \texttt{optGrbModel} is more lightweight but less flexible for users. To illustrate, consider the following optimization model \eqref{equa:example}, where $c_i$ represents an unknown objective function coefficient:

\begin{equation}
\label{equa:example}
\begin{aligned}
\max_{x} & \sum_{i=0}^4 c_i x_i \\
s.t. \quad & 3 x_0 + 4 x_1 + 3 x_2 + 6 x_3 + 4 x_4 \leq 12 \\
& 4 x_0 + 5 x_1 + 2 x_2 + 3 x_3 + 5 x_4 \leq 10 \\
& 5 x_0 + 4 x_1 + 6 x_2 + 2 x_3 + 3 x_4 \leq 15 \\
& \forall x_i \in \{0, 1\}
\end{aligned}
\end{equation}

Inheriting from \texttt{optGrbModel} provides a convenient way to use \gurobi{} with \pyepo{}. The only implementation required is to override \texttt{\_getModel} and return a \gurobi{} model and the corresponding decision variables. Furthermore, there is no need to manually assign a value to the attribute \texttt{modelSense}. An example of the model \eqref{equa:example} is as follows:

\begin{lstlisting}[language=Python]
import gurobipy as gp
from gurobipy import GRB
from pyepo.model.grb import optGrbModel

class myModel(optGrbModel):

    def _getModel(self):
        # ceate a model
        m = gp.Model()
        # varibles
        x = m.addVars(5, name="x", vtype=GRB.BINARY)
        # sense
        m.modelSense = GRB.MAXIMIZE
        # constraints
        m.addConstr(3*x[0]+4*x[1]+3*x[2]+6*x[3]+4*x[4]<=12)
        m.addConstr(4*x[0]+5*x[1]+2*x[2]+3*x[3]+5*x[4]<=10)
        m.addConstr(5*x[0]+4*x[1]+6*x[2]+2*x[3]+3*x[4]<=15)
        return m, x
        
optmodel = myModel()
\end{lstlisting}

The same applies to \texttt{optCoptModel} with \copt{}.

\subsubsection{Optimization Model with Pyomo}

Similarly, \texttt{optOmoModel} allows modeling mathematical programs with \pyomo{}. Unlike \texttt{optGrbModel}, \texttt{optOmoModel} requires explicitly setting the \texttt{modelSense}. Since \pyomo{} supports multiple solvers, instantiating an \texttt{optOmoModel} requires a parameter \texttt{solver} to specify the solver. The following is an implementation of the model \eqref{equa:example}:

\begin{lstlisting}[language=Python]
from pyomo import environ as pe
from pyepo import EPO
from pyepo.model.omo import optOmoModel

class myModel(optOmoModel):

    def _getModel(self):
        # sense
        self.modelSense = EPO.MAXIMIZE
        # ceate a model
        m = pe.ConcreteModel()
        # varibles
        x = pe.Var([0,1,2,3,4], domain=pe.Binary)
        m.x = x
        # constraints
        m.cons = pe.ConstraintList()
        m.cons.add(3*x[0]+4*x[1]+3*x[2]+6*x[3]+4*x[4]<=12)
        m.cons.add(4*x[0]+5*x[1]+2*x[2]+3*x[3]+5*x[4]<=10)
        m.cons.add(5*x[0]+4*x[1]+6*x[2]+2*x[3]+3*x[4]<=15)
        return m, x

optmodel = myModel(solver="glpk")
\end{lstlisting}

\subsection{Autograd Functions}
\label{subsec:func}

Training neural networks with modern deep learning libraries such as \tensorflow{} \cite{abadi2016tensorflow} or \pytorch{} \cite{NEURIPS2019_9015} requires gradient calculations for backpropagation. This is achieved through the numerical technique of automatic differentiation \cite{paszke2017automatic}. For example, \pytorch{} provides autograd functions.

The autograd functions are the core modules of \pyepo{} that solve and backpropagate the optimization problems with predicted coefficients. These functions can be integrated with different neural network architectures to achieve end-to-end predict-then-optimize for various tasks. 
%In \pyepo{}, the functions include \texttt{SPOPlus} \cite{elmachtoub2021smart}, \texttt{blackboxOpt} \cite{poganvcic2019differentiation}, \texttt{perturbedOpt} \cite{berthet2020learning}, \texttt{perturbedFenchelYoung} \cite{berthet2020learning}, etc.

\subsubsection{Function SPOPlus}

The function \texttt{SPOPlus} directly calculates the SPO+ loss, which measures the decision error of an optimization solved with predicted coefficients. This optimization is represented as an instance of \texttt{optModel} and passed into \texttt{SPOPlus} as an argument. As shown below, \texttt{SPOPlus} also requires \texttt{processes} to specify the number of processes for parallel computation.

\begin{lstlisting}[language=Python]
from pyepo.func import SPOPlus
# init SPO+ Pytorch function
spop = SPOPlus(optmodel, processes=8)
\end{lstlisting}
The parameters for the forward pass of \texttt{SPOPlus} are as follows:

\begin{itemize}
  \item \texttt{pred\_cost}: a batch of predicted coefficient vectors, one vector per instance;
  \item \texttt{true\_cost}: a batch of true coefficient vectors, one vector per instance;
  \item \texttt{true\_sol}: a batch of true optimal solutions, one vector per instance;
  \item \texttt{true\_obj}: a batch of true optimal objective values, one value per instance.
\end{itemize}
The following code block is the \texttt{SPOPlus} forward pass:

\begin{lstlisting}[language=Python]
# calculate SPO+ loss
loss = spop(pred_cost, true_cost, true_sol, true_obj)
\end{lstlisting}

\subsubsection{Function blackboxOpt}

\texttt{SPOPlus} directly computes a loss, whereas \texttt{blackboxOpt} provides a solution as a layer. Thus, \texttt{blackboxOpt} makes it possible to use various loss functions. Compared to \texttt{SPOPlus}, \texttt{blackboxOpt} requires an additional parameter \texttt{lambd}, which is the interpolation degree $\lambda$ for finite difference. According to \citet{poganvcic2019differentiation}, the values of $\lambda$ should be between $10$ and $20$.

\begin{lstlisting}[language=Python]
from pyepo.func import blackboxOpt
# init DBB solver
dbb = blackboxOpt(optmodel, lambd=10, processes=8)
# set some loss
l1 = nn.L1Loss()
\end{lstlisting}

Since \texttt{blackboxOpt} functions as a differentiable optimizer, there is only one parameter \texttt{pred\_cost} for the forward pass. As shown in the code below, the output is the optimal solution for the given predicted coefficients:

\begin{lstlisting}[language=Python]
# find the optimal solution
pred_sol = dbb(pred_cost)
# calculate additional loss
loss = l1(pred_sol, true_sol)
\end{lstlisting}

\subsubsection{Function perturbedOpt}

Same as \texttt{blackboxOpt}, \texttt{perturbedOpt} is a differentiable optimizer and provides an ``expected solution". The hyperparameters for \texttt{perturbedOpt} include the number $K$ of Monte-Carlo sample, ``\texttt{n\_sample}" and the amplitude parameter $\sigma$ of the perturbation, ``\texttt{sigma}".

\begin{lstlisting}[language=Python]
import pyepo
# init DPO solver
dpo = pyepo.func.perturbedOpt(optmodel, n_samples=3, sigma=1.0, processes=8)
# set some loss
l1 = nn.L1Loss()
\end{lstlisting}

Given predicted coefficients $\hat{\bm{c}}$, \texttt{perturbedOpt} outputs an expected solution by averaging the solutions of $K$ randomly perturbed optimization problems:

\begin{lstlisting}[language=Python]
# find the expected optimal solution
exp_sol = dpo(pred_cost)
# calculate loss
loss = l1(exp_sol, true_sol)
\end{lstlisting}

\subsubsection{Function perturbedFenchelYoung}

\texttt{perturbedFenchelYoung} uses a predicted coefficients $\hat{\bm{c}}$ and  a true solution $\bm{w}^*(\bm{c})$ to compute the Perturbed Fenchel-Young loss $\mathcal{L}_{\text{FY}} (\hat{\bm{c}}, \bm{w}^*(\bm{c}))$; it requires the same hyperparameters as \texttt{perturbedOpt}.

\begin{lstlisting}[language=Python]
import pyepo
# init Fenchel-Young loss
pfyl = pyepo.func.perturbedFenchelYoung(optmodel, n_samples=3, sigma=1.0, processes=8)
\end{lstlisting}
The below code block illustrates the calculation of Fenchel-Young loss:

\begin{lstlisting}[language=Python]
# calculate Fenchel-Young loss
loss = pfyl(pred_cost, true_sol)
\end{lstlisting}

\subsection{The \texttt{optDataset} Class for Managing Data}
\label{subsec:optdata}

The use of decision losses, such as \spo{} and \pfyl{}, requires access to true optimal solutions. To facilitate convenient training and testing in \pyepo{}, an auxiliary component \texttt{optDataset} has been introduced, though it is not strictly indispensable. \texttt{optDataset} stores the features and their associated coefficients of the objective function and solves optimization problems to obtain optimal solutions and corresponding objective values.

\texttt{optDataset} is extended from \pytorch{} Dataset. In order to obtain optimal solutions, \texttt{optDataset} requires the corresponding \texttt{optModel} (see in Section \ref{subsec:optmodel}). The parameters for the \texttt{optDataset} are as follows:

\begin{itemize}
  \item \texttt{model}: an instance of \texttt{optModel};
  \item \texttt{feats}: data features;
  \item \texttt{costs}: corresponding objective function coefficients;
\end{itemize}

Then, \pytorch{} DataLoader receives an \texttt{optDataset}, wraps the data samples, and acts as a sampler that provides an iterable over the given dataset. The batch size is required as the number of training samples that will be used in each update of the model parameters.

\begin{lstlisting}[language=Python]
import pyepo
from torch.utils.data import DataLoader

# build dataset
dataset = pyepo.data.dataset.optDataset(optmodel, feats, costs)
# get data loader
dataloader = DataLoader(dataset, batch_size=32, shuffle=True)
\end{lstlisting}

By iterating over the DataLoader, we can obtain a batch of features, true objective function coefficients, optimal solutions, and the corresponding objective values:

\begin{lstlisting}[language=Python]
for x, c, w, z in dataloader:
    # a batch of features 
    print(x)
    # a batch of true coefficients
    print(c)
    # a batch of true optimal solutions 
    print(w)
    # a batch of true optimal objective values
    print(z)
\end{lstlisting}

\subsection{End-to-End Training}

The core capability of \pyepo{} is to build an optimization model and then embed it into a \pytorch{} neural network for end-to-end training.  Here, we illustrate this with a simple linear regression model in PyTorch as an example:

\begin{lstlisting}[language=Python]
from torch import nn

# construct linear model
class linearRegression(nn.Module):

    def __init__(self):
        super(linearRegression, self).__init__()
        # fully-connected layer without activation
        self.linear = nn.Linear(num_feat, len_cost)
        
    def forward(self, x):
        out = self.linear(x)
        return out

# init model
predmodel = linearRegression()
\end{lstlisting}

We can then train the prediction model with SPO+ loss to predict unknown coefficients, make decisions, and compute decision errors. The prediction model is trained using a stochastic gradient descent (SGD) optimizer. Leveraging the automatic differentiation capabilities of \pytorch{}, the loss gradients with respect to the model parameters $\bm{\theta}$ are calculated and used to update the model parameters during training.

\begin{lstlisting}[language=Python]
import torch

# set SGD optimizer
optimizer = torch.optim.SGD(predmodel.parameters(), lr=1e-3)

# training
for epoch in range(num_epochs):
    # iterare features, coefficients, solutions, and objective values
    for x, c, w, z in dataloader:
        # forward pass
        cp = predmodel(x) # predict coefficients
        loss = spop(cp, c, w, z) # calculate SPO+ loss
        # backward pass
        optimizer.zero_grad() # reset gradients to 0
        loss.backward() # compute gradients
        optimizer.step() # update model parameters
\end{lstlisting}

\subsection{Metrics}
\label{subsec:metr}

\pyepo{} provides evaluation functions to measure model performance, in particular the two metrics mentioned in Section \ref{subsec:loss}: regret and unambiguous regret. Additionally, we define the normalized (unambiguous) regret by $$\frac{\sum_{i=1}^{n_\text{test}}{\mathcal{L}_{\text{Regret}}(\hat{\bm{c}}^i, \bm{c}^i)}}{\sum_{i=1}^{n_\text{test}}{|z^*(\bm{c}^i)|}}.$$ Both require the following parameters:

\begin{itemize}
  \item \texttt{predmodel}: a regression neural network for coefficients prediction.
  \item \texttt{optModel}: a \pyepo{} optimization model.
  \item \texttt{dataloader}: a \pyepo{} data loader.
\end{itemize}

Assume that we have trained a prediction model \texttt{predmodel} for an optimization problem \texttt{optModel}. To evaluate the performance of \texttt{predmodel} on a dataset, \texttt{dataloader}, containing instances of \texttt{optModel}, the following suffices:

\begin{lstlisting}[language=Python]
from pyepo.metric import regret, unambRegret
# compute normalized regret
regret = regret(predmodel, optmodel, dataloader)
# compute normalized unambiguous regret
uregret = unambRegret(predmodel, optmodel, dataloader)
\end{lstlisting}

\section{Benchmark Datasets}
\label{sec:data}

\subsection{Benchmark Datasets from \pyepo{}}
\label{subsec:syndata}

This section describes our new datasets designed for the end-to-end predict-then-optimize task. In general, we generate data sets in a way similar to \citet{elmachtoub2021smart}. The synthetic dataset $\mathcal{D}$ includes features $\bm{x}$ and objective function coefficients $\bm{c}$ (e.g. $\mathcal{D} = \{(\bm{x}^1, \bm{c}^1), (\bm{x}^2, \bm{c}^2), ..., (\bm{x}^n, \bm{c}^n)\}$). The feature vector $\bm{x}^i \in \mathbb{R}^p$ follows a standard multivariate Gaussian distribution $\mathcal{N}(0, \bm{I}_p)$ and the corresponding coefficient $\bm{c}^i \in \mathbb{R}^d$ comes from a (possibly nonlinear) polynomial function of $\bm{x}^i$ with additional random noise. $\epsilon^i_j \sim  \bm{U}(1-\bar{\epsilon}, 1+\bar{\epsilon})$ is the multiplicative noise term for $c^i_j$, the $j^{th}$ element of coefficients $\bm{c}^i$. 

Our dataset includes three of the most classical optimization problems: the shortest path problem, the multi-dimensional knapsack problem, and the traveling salesperson problem. \pyepo{} provides functions for generating these data with adjustable parameters: data size $n$, number of features $p$, coefficients vector dimension $d$, polynomial degree $deg$, and noise half-width $\bar{\epsilon}$.

\subsubsection{Shortest Path}
\label{subsec:sp}

Following \citet{elmachtoub2021smart}, we consider a $h \times w$ grid network, and the goal is to find the shortest path \cite{ortega2014shortest} from the northwest to the southeast. We generate a random matrix $\mathcal{B} \in \mathbb{R}^{d \times p}$, where $\mathcal{B}_{ij}$ follows the Bernoulli distribution with probability $0.5$. Then, the coefficients $\bm{c}^i$ is almost the same as in \cite{elmachtoub2021smart}, and is generated from 
\begin{equation}
    \Bigg[\frac{1}{{3.5}^{deg}} \Bigg(\frac{1}{\sqrt{p}}(\mathcal{B} \bm{x}^i)_j + 3\Bigg)^{deg} + 1\Bigg] \cdot \epsilon^i_j.
    \label{eq:randomcost}
\end{equation}

The following code generates data for the shortest path on the grid network:

\begin{lstlisting}[language=Python]
from pyepo.data.shortestpath import genData
x, c = genData(n, p, grid=(h,w), deg=deg, noise_width=e)
\end{lstlisting}

\subsubsection{Multi-Dimensional Knapsack}
\label{subsec:ks}

The multi-dimensional knapsack problem \cite{martello1990knapsack} is one of the most well-known integer programming models. It maximizes the value of the selected items under multiple resource constraints. Due to its computational complexity, solving this problem can be challenging as the number of variables (items) and constraints (resources) grows.

Assuming that the uncertain coefficients exist only in the objective function, the weights of the items and the capacity of resources remain fixed throughout the data. We denote the number of resources by $k$ and the number of items by $d$. The weights $\mathcal{W} \in \mathbb{R}^{k \times m}$ are sampled from $3$ to $8$ with a precision of one decimal place. Using the same $\mathcal{B} \in \mathbb{R}^{d \times p}$ as in Section \ref{subsec:sp}, coefficient $c_{ij}$ is calculated according to Equation~\eqref{eq:randomcost}.%from $$\bigg\lceil \frac{5}{{3.5}^{deg}} \bigg[\frac{1}{\sqrt{p}} ((\mathcal{B} \bm{x_i})_j + 3)^{deg} + 1\bigg] \cdot \epsilon_{ij} \bigg\rceil.$$

To generate $k$-dimensional knapsack data, a user simply executes the following:

\begin{lstlisting}[language=Python]
from pyepo.data.knapsack import genData
W, x, c = genData(n, p, num_item=d, dim=dim, deg=deg, noise_width=e)
\end{lstlisting}

\subsubsection{Traveling Salesperson}
\label{subsec:tsp}

As one of the most renowned combinatorial optimization problems, the traveling salesperson problem (TSP) seeks to determine the shortest possible tour that visits every node exactly once. Here, we introduce the symmetric TSP (the distance between two nodes is the same in both directions) with the number of nodes to be visited $v$.

\pyepo{} generates costs from a distance matrix, where the distance is composed of two parts: one derived from the Euclidean distance, and the other from a polynomial function of the features. For the Euclidean distance, two-dimensional coordinates are created from the mixture of the Gaussian distribution $\mathcal{N}(0, I)$ and the uniform distribution $\bm{U}(-2, 2)$. For the function of features, the polynomial kernel function is $$\frac{1}{{3}^{deg - 1}} \bigg(\frac{1}{\sqrt{p}}(\mathcal{B} \bm{x}^i)_j + 3\bigg)^{deg} \cdot \epsilon^i_j,$$ where the elements of $\mathcal{B}$ come from the multiplication of Bernoulli $\bm{B}(0.5)$ and uniform $\bm{U}(-2, 2)$.

An example of a TSP data generation is as follows:

\begin{lstlisting}[language=Python]
from pyepo.data.tsp import genData
x, c = genData(n, p, num_node=v, deg=deg, noise_width=e)
\end{lstlisting}

\subsection{Warcraft Terrain Map Images}
\label{subsec:wcdata}

\begin{figure}[htbp]
    \centering
    \subfloat[]{\includegraphics[width=0.28\textwidth]{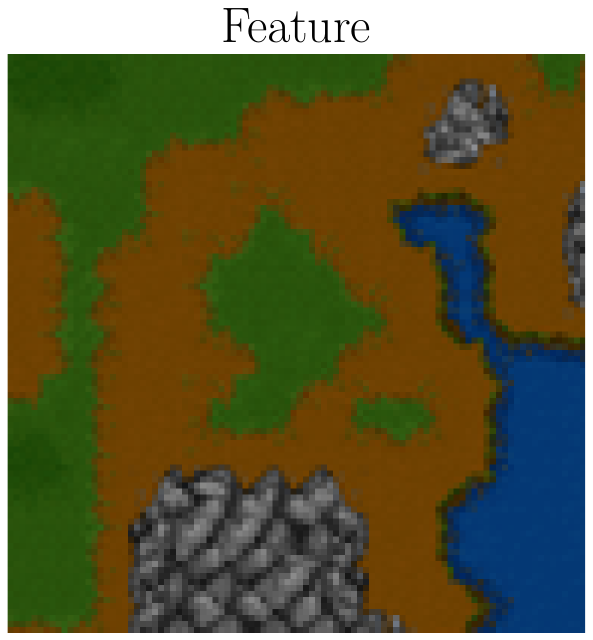}}
    \subfloat[]{\includegraphics[width=0.28\textwidth]{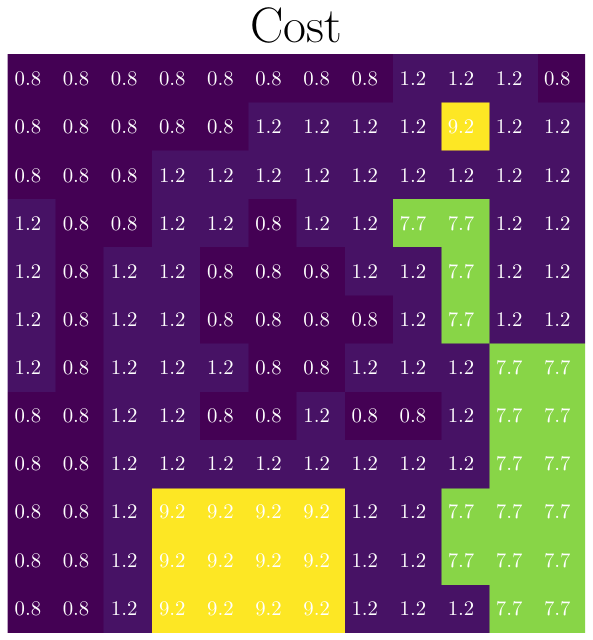}}
    \subfloat[]{\includegraphics[width=0.28\textwidth]{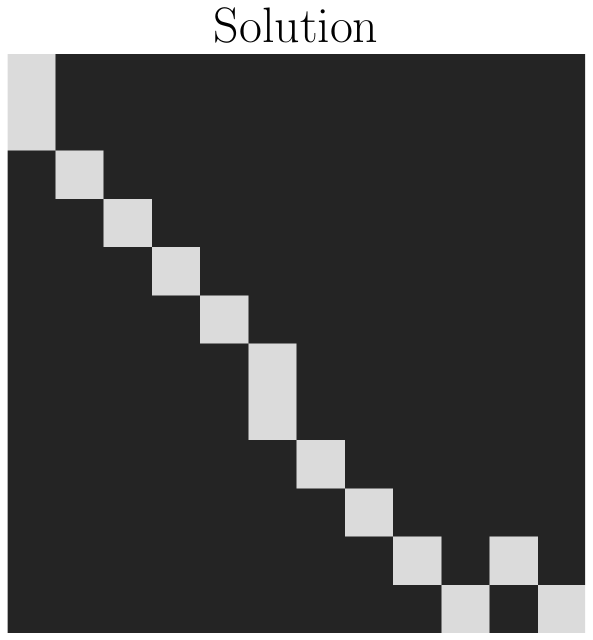}} \\
    \caption{Warcraft terrain shortest path dataset: (left) Each input feature is a $k \times k$ terrain map image as a grid of tiles; (Middle) the respective weights is a matrix indicating traveling costs; (right) the corresponding binary matrix represents the shortest path from top left to bottom right.}
    \label{fig:warcraft}
\end{figure}

\citet{poganvcic2019differentiation} proposed a dataset of maps from the popular video game Warcraft (see Figure \ref{fig:warcraft}), designed for learning the shortest path from RGB terrain images. This dataset stands out as a significant benchmark because it employs image input, a modality that was not explored in other work in this area. In accordance with the experiment of ~\citet{poganvcic2019differentiation} and \citet{berthet2020learning}, we utilize $96 \times 96$ RGB images for determining the shortest path on $12 \times 12$ grid networks.

\section{Empirical Evaluation for \pyepo{} Datasets}
\label{sec:eval}

This section presents experimental results for the benchmark datasets described in Section \ref{subsec:syndata}. The experiments aimed to investigate training time and normalized regret, as defined in Section \ref{subsec:metr}, on a test set of size $n_\text{test}=1000$. As Table \ref{tab:model} shows, the methods we compare include the two-stage approach with different predictors as well as end-to-end methods such as \spo{}/\dbb{}/\pfyl{} with a linear prediction model $g(\bm{x}; \bm{\theta})$. Notably, the result of \dpo{} was excluded due to its consistently subpar performance.

Unlike direct decision loss functions \spo{} and \pfyl{}, \dbb{} and \dpo{} allow the use of customized loss functions, offering flexibility for various problems. In the original paper, \citet{poganvcic2019differentiation} used the Hamming distance between the true optimum and the solution from prediction, whereas  \citet{berthet2020learning} applied the square difference between solutions. However, in our experiments, compared to regret, \dbb{} using the Hamming distance is only sensible for the shortest path problem but leads to much worse decisions in knapsack and TSP. For the sake of consistency, we only use regret~\eqref{eq:regret} as a loss for \dbb{} and \dpo{}. In addition, $l_1$/$l_2$ regularization means that the MAE / MSE of the predicted coefficients $\hat{\bm{c}}$ is added to the loss function. 

\begin{table}[htbp]
    \begin{center}
    \resizebox{1.0\textwidth}{!}{
    \begin{tabular}{l|l}
        \hline
        \textbf{Method}     & \textbf{Description} \\ \hline
        2-stage LR          & Two-stage method where the predictor is a linear regression \\ \hline
        2-stage RF          & Two-stage method where the predictor is a random forest with default parameters  \\ \hline
        2-stage Auto        & Two-stage method where the predictor is \autosklearn{} \cite{feurer-neurips15a} with $10$ minutes time limit and uses MSE as metric \\ \hline
        \spo{}              & Linear model with SPO+ loss \cite{elmachtoub2021smart} \\ \hline
        \pfyl{}             & Linear model with perturbed Fenchel-Young loss \cite{berthet2020learning} \\ \hline
        \dbb{}              & Linear model with differentiable black-box optimizer \cite{poganvcic2019differentiation}  and regret loss \\ \hline
        \sporel{}           & Linear model with SPO+ loss \cite{elmachtoub2021smart}, using linear relaxation for training \\ \hline
        \pfylrel{}          & Linear model with perturbed Fenchel-Young loss \cite{berthet2020learning}, using linear relaxation for training \\ \hline
        \dbbrel{}           & Linear model with differentiable black-box optimizer \cite{poganvcic2019differentiation} and regret loss, using linear relaxation for training \\ \hline
        \spolo{}            & Linear model with SPO+ loss \cite{elmachtoub2021smart}, using $l_1$ regularization for coefficients \\ \hline
        \spolt{}            & Linear model with SPO+ loss \cite{elmachtoub2021smart}, using $l_2$ regularization for coefficients \\ \hline
        \pfyllo{}           & Linear model with perturbed Fenchel-Young loss \cite{berthet2020learning}, using $l_1$ regularization for coefficients \\ \hline
        \pfyllt{}           & Linear model with perturbed Fenchel-Young loss \cite{berthet2020learning}, using $l_2$ regularization for coefficients \\ \hline
        \dbblo{}            & Linear model with differentiable black-box optimizer \cite{poganvcic2019differentiation} and regret loss , using $l_1$ regularization for coefficients \\ \hline
        \dbblt{}            & Linear model with differentiable black-box optimizer \cite{poganvcic2019differentiation} and regret loss , using $l_2$ regularization for coefficients \\ \hline
        %\spohl{}            & Fully connected neural network with $L$ hidden layers of width $h_1$, ..., $h_L$ and SPO+ loss \cite{elmachtoub2021smart} \\  \hline
        %\pfylhl{}            & Fully connected neural network with $L$ hidden layers of width $h_1$, ..., $h_L$ and perturbed Fenchel-Young loss \cite{berthet2020learning} \\  \hline
        %\dbbhl{}           & Fully connected neural network with $L$ hidden layers of width $h_1$, ..., $h_L$ and differentiable black-box optimizer \cite{poganvcic2019differentiation} \\ \hline
    \end{tabular}}
    \end{center}
    \caption{Methods compared in the experiments.}\label{tab:model}
\end{table}

All numerical experiments were conducted in Python v3.7.9 with 32 Intel E5-2683 v4 Broadwell CPUs and 32GB memory. Specifically, we used \pytorch{} \cite{NEURIPS2019_9015} v1.10.0 to train end-to-end models, and \sklearn{} \cite{scikit-learn} v0.24.2 and \autosklearn{} \cite{feurer-neurips15a} v0.14.6 as predictors for the two-stage method. \gurobi{} \cite{gurobi} v9.1.2 was the optimization solver used throughout. The version of \pyepo{} is v0.2.4\footnote{\pyepo{} is an actively maintained open-source library. To ensure the reproducibility of the results, a static branch has been established and can be accessed at \url{https://github.com/khalil-research/PyEPO/blob/MPC}.}.

\subsection{Performance Comparison between Different Methods}

We compare the performance between two-stage methods, \spo{}, \pfyl{}, and \dbb{} between varying training data size $n \in \{100, 1000, 5000\}$, polynomial degree $deg \in \{1, 2, 4, 6\}$, and noise half-width $\bar{\epsilon} \in \{0.0, 0.5\}$.  Each experiment was repeated $10$ times, using different $\bm{x}$, $\mathcal{B}$, and $\epsilon$ to generate $10$ distinct training/validation/test datasets. We use boxplots to summarize the statistical outcomes.

We also carried out small-scale experiments on the validation set to select gradient descent hyperparameters, namely batch size, learning rate, and momentum for the shortest path, knapsack, and TSP in \spo{}, \pfyl{} and \dbb{}. The hyperparameter tuning used a limited random search in the space of possible configurations, thus not guaranteeing the best performance in the results. For \pfyl{}, we arbitrarily set the number of samples $K = 1$ and the perturbation amplitude $\sigma = 1.0$.

\begin{table}[htbp]
    \begin{center}
    \begin{tabular}{l|l|l|l}
        \hline
        \textbf{Problem}  & \textbf{Parameters} & \textbf{Feature Size} & \textbf{Label Size} \\ \hline
        Shortest Path     & \begin{tabular}[c]{@{}l@{}}Height of the grid is $5$\\ Width of the grid is $5$\end{tabular}                        & 5                     & 40                      \\ \hline
        Knapsack          & \begin{tabular}[c]{@{}l@{}}Dimension of resource is   $2$\\ Number of items is $32$\\ Capacity is $20$\end{tabular} & 5                     & 32                      \\ \hline
        Traveling Salesperson & Number of nodes is $20$                                                                                             & 10                    & 190                     \\ \hline
    \end{tabular}
    \end{center}
    \caption{Problem Parameters for Performance Comparison}\label{tab:prob}
\end{table}

We generate synthetic datasets using the parameters outlined in Table \ref{tab:prob}, so the label sizes (the number of unknown objective function coefficients) $d$ are $40$, $32$, and $190$, respectively. For the TSP, we use the Dantzig–Fulkerson–Johnson (DFJ) formulation \cite{dantzig1954solution} because it is faster to solve than the alternatives.

\begin{figure}[htbp]
    \centering
    \subfloat[]{\includegraphics[width=0.48\textwidth]{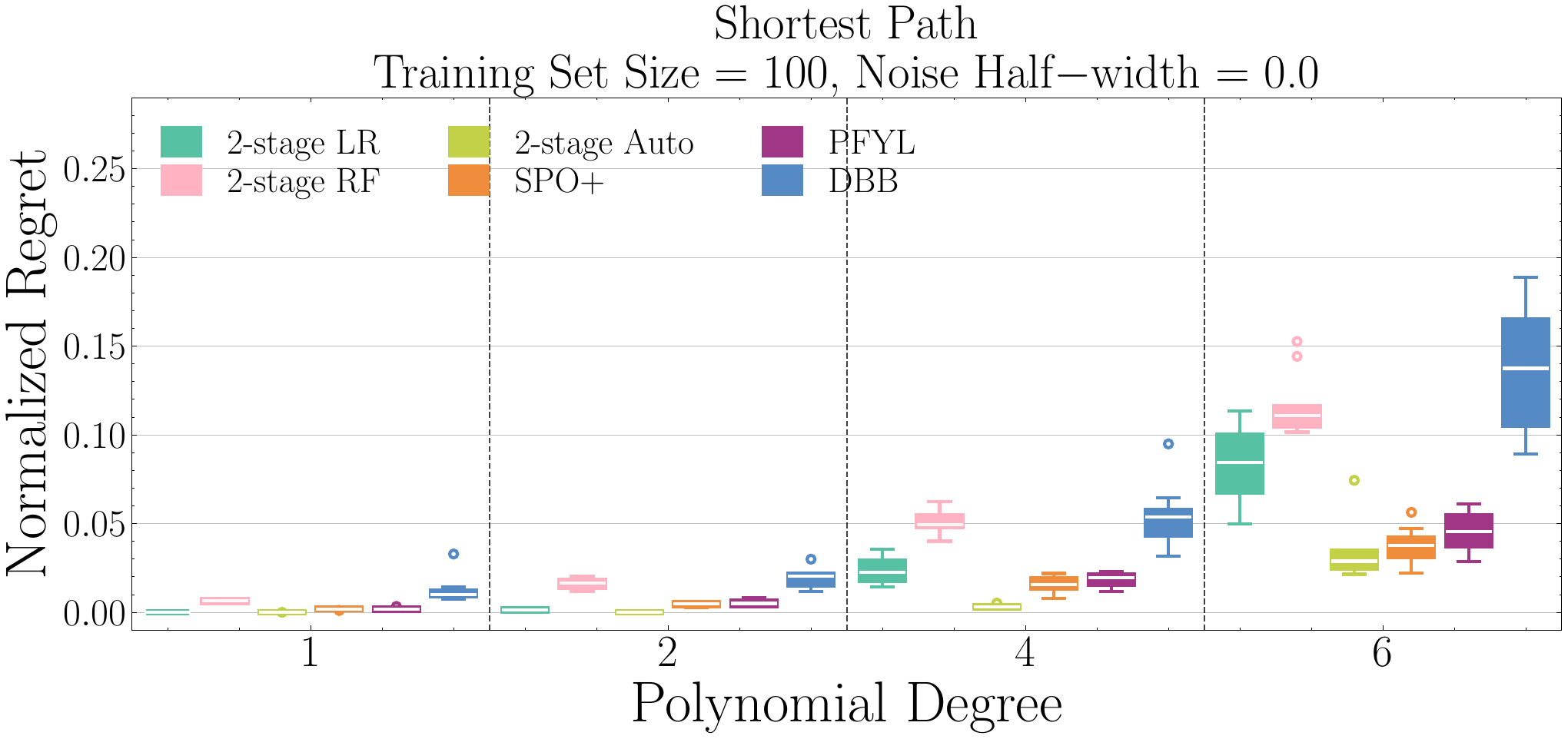}} 
    \subfloat[]{\includegraphics[width=0.48\textwidth]{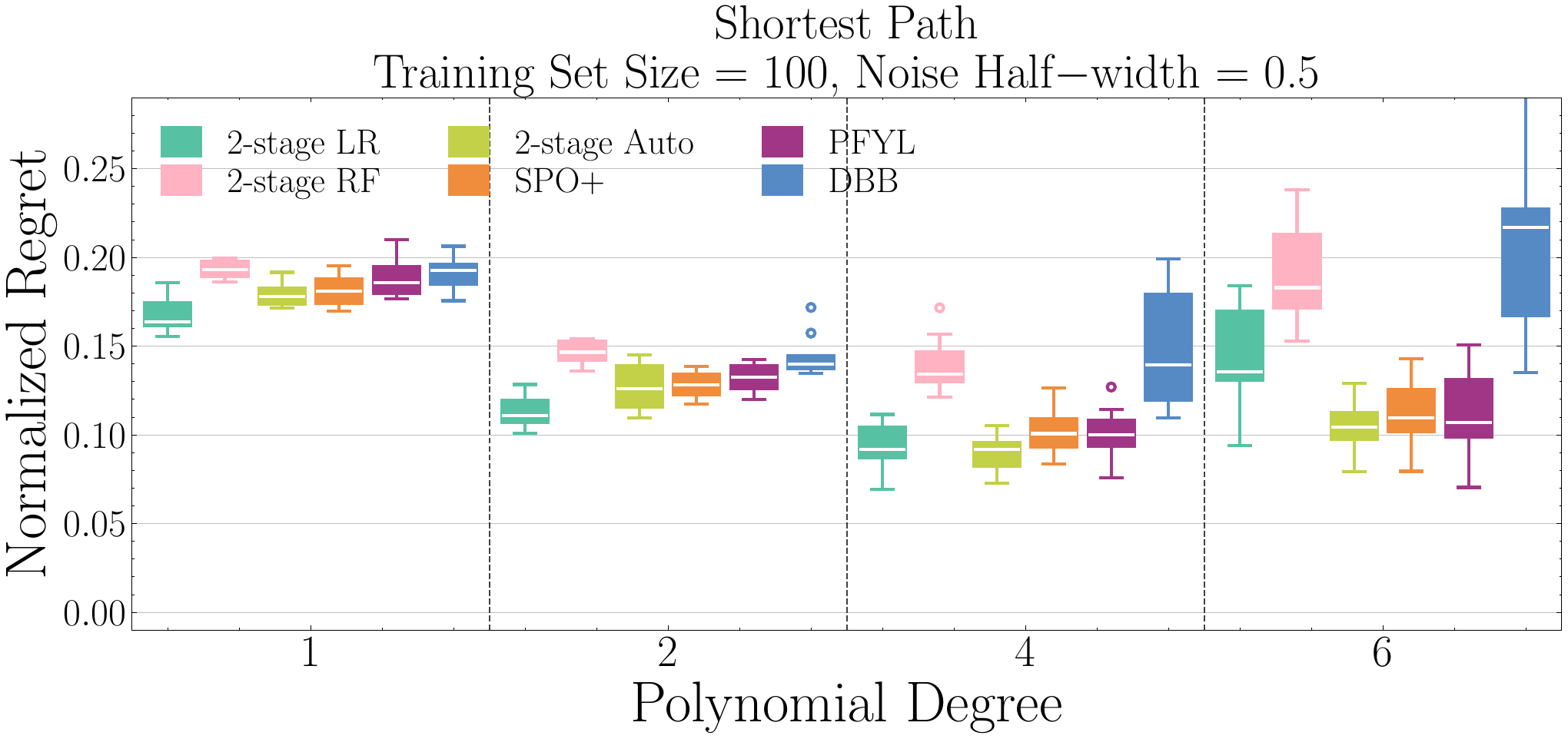}} \vspace{-2\baselineskip}
    \subfloat[]{\includegraphics[width=0.48\textwidth]{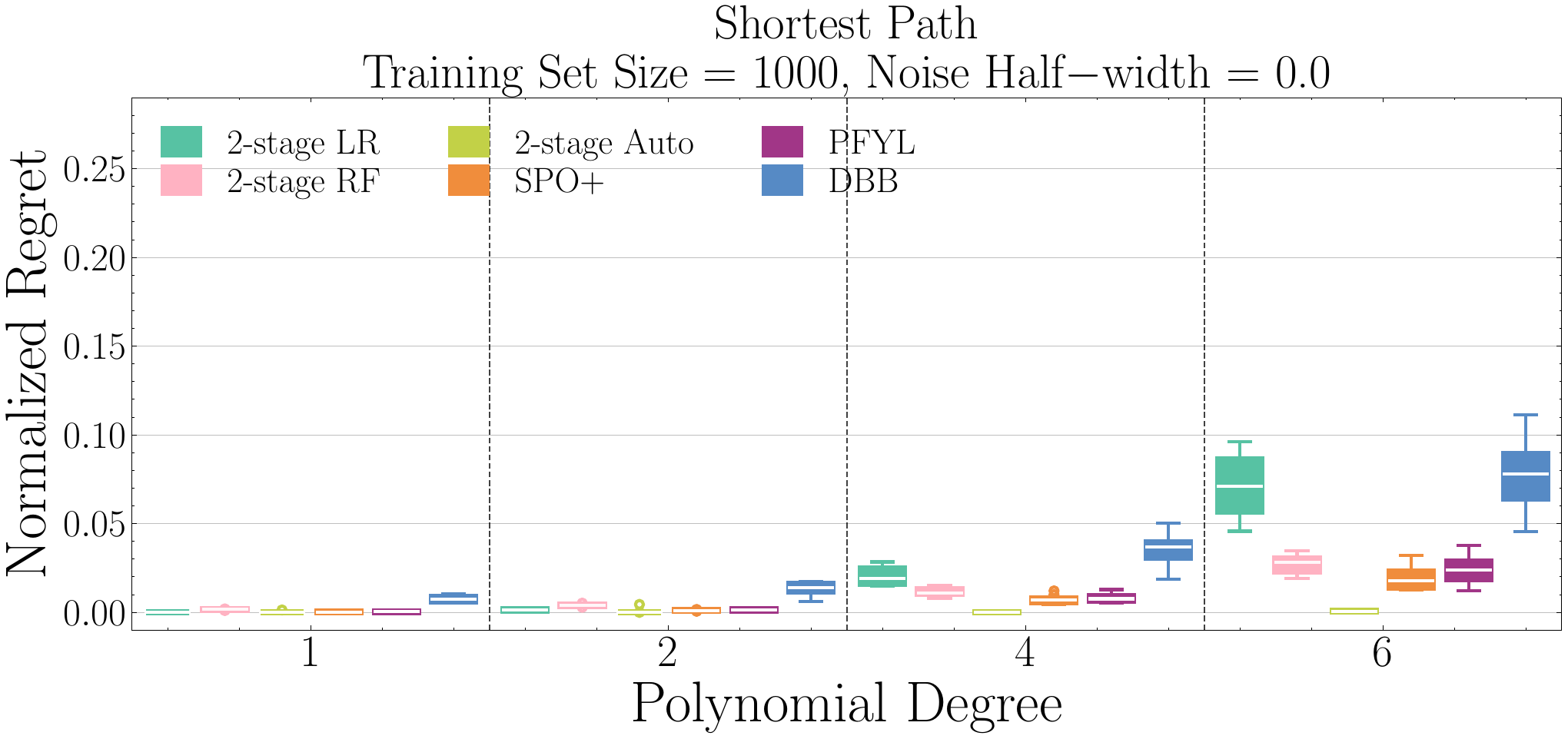}} 
    \subfloat[]{\includegraphics[width=0.48\textwidth]{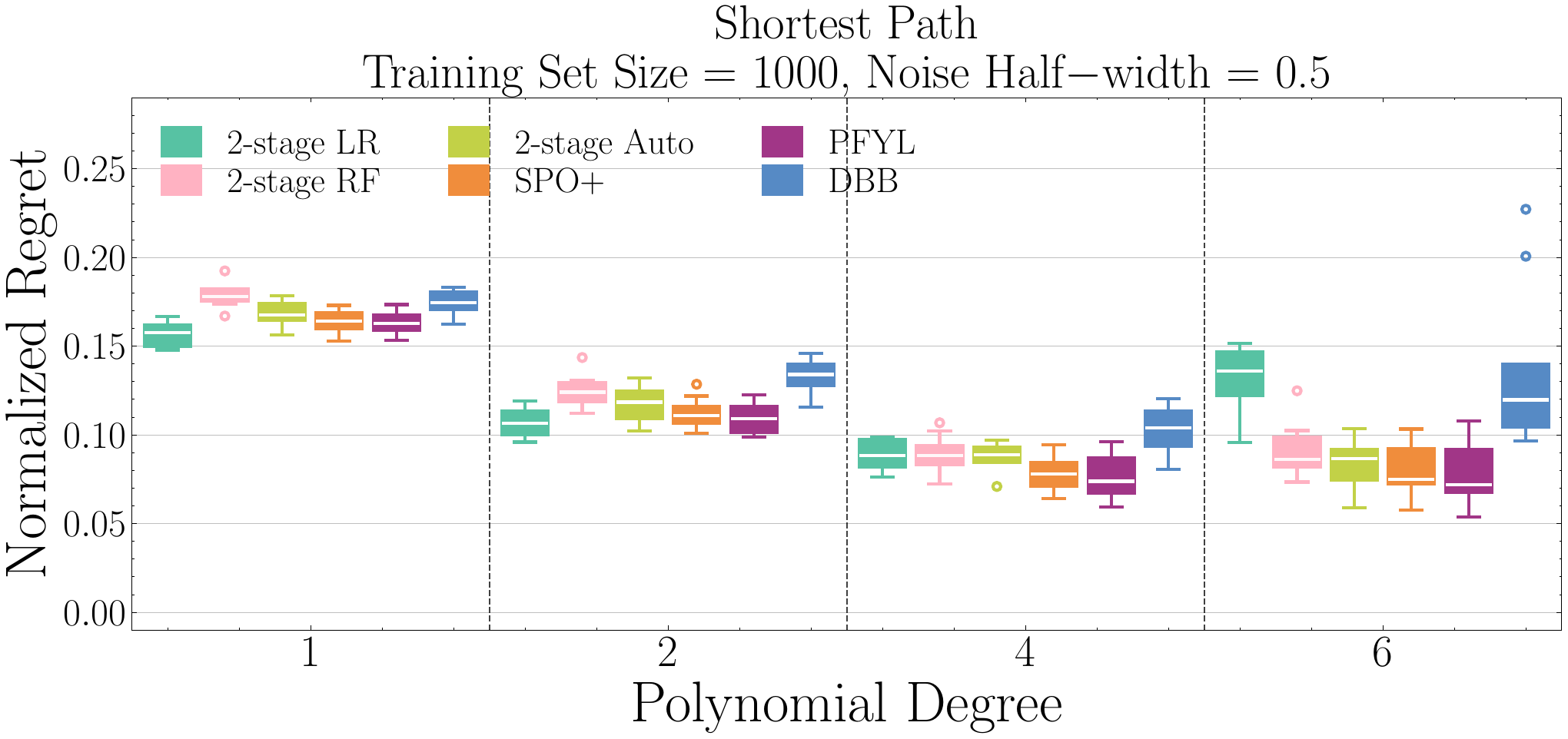}} \vspace{-2\baselineskip}
    \subfloat[]{\includegraphics[width=0.48\textwidth]{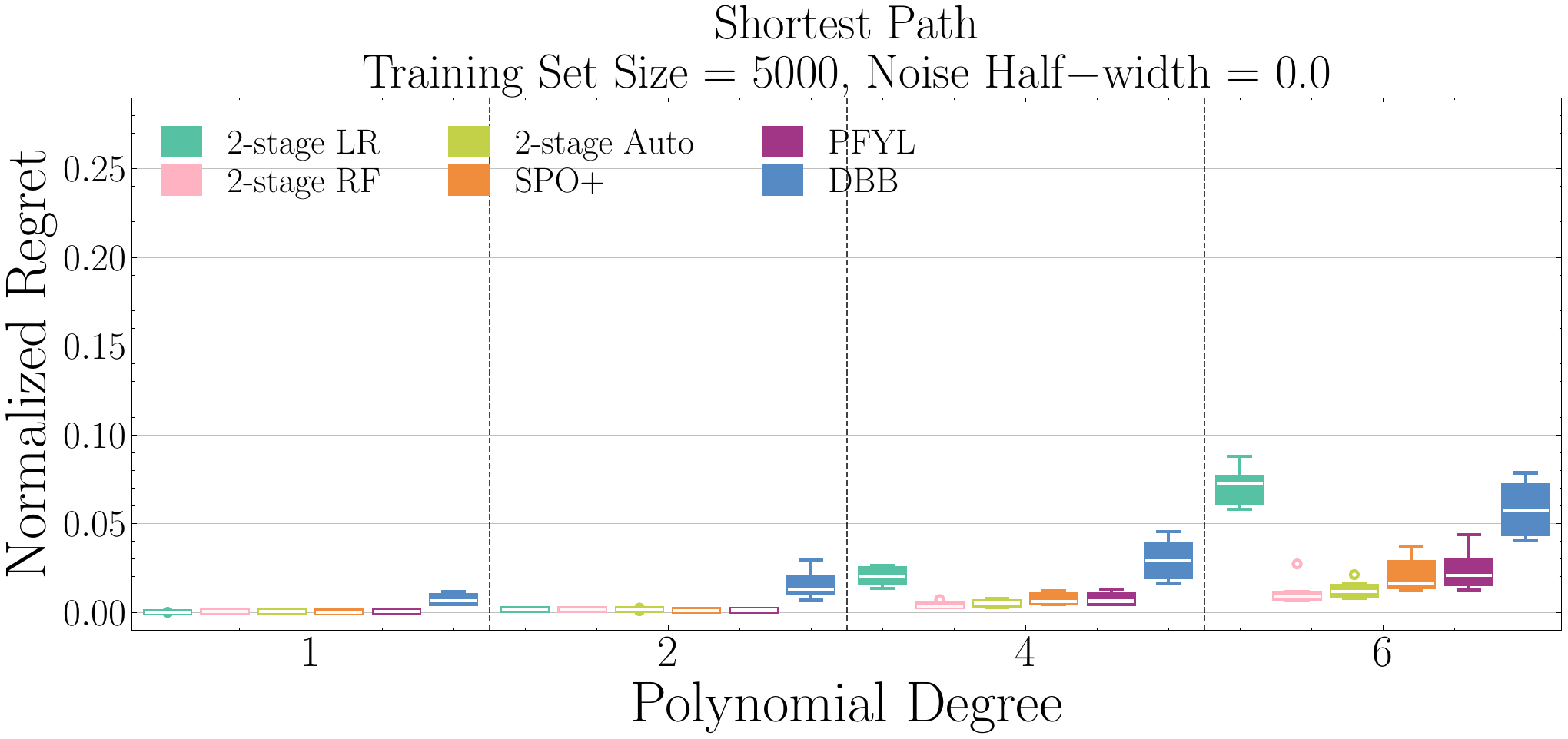}} 
    \subfloat[]{\includegraphics[width=0.48\textwidth]{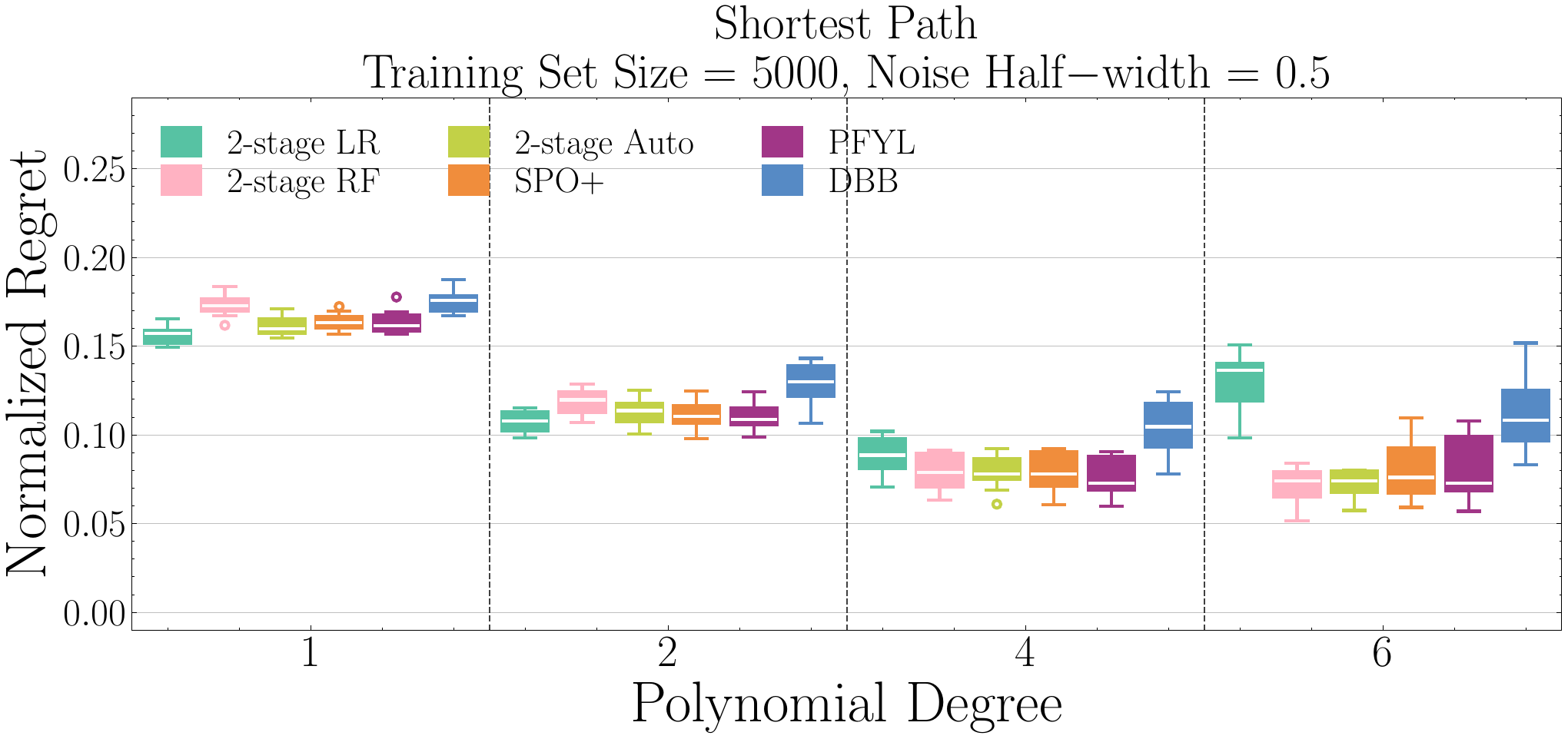}} \\
    \caption{Normalized regret for the shortest path problem on the test set: The size of the grid network is $5 \times 5$. The methods in the experiment include two-stage approaches with linear regression, random forest, and \autosklearn{} and end-to-end learning such as \spo{}, \pfyl{}, and \dbb{}. The normalized regret is visualized under different sample sizes, noise half-width, and polynomial degrees. For normalized regret, lower is better.}
    \label{fig:sp}
\end{figure}

\begin{figure}[htbp]
    \centering
    \subfloat[]{\includegraphics[width=0.48\textwidth]{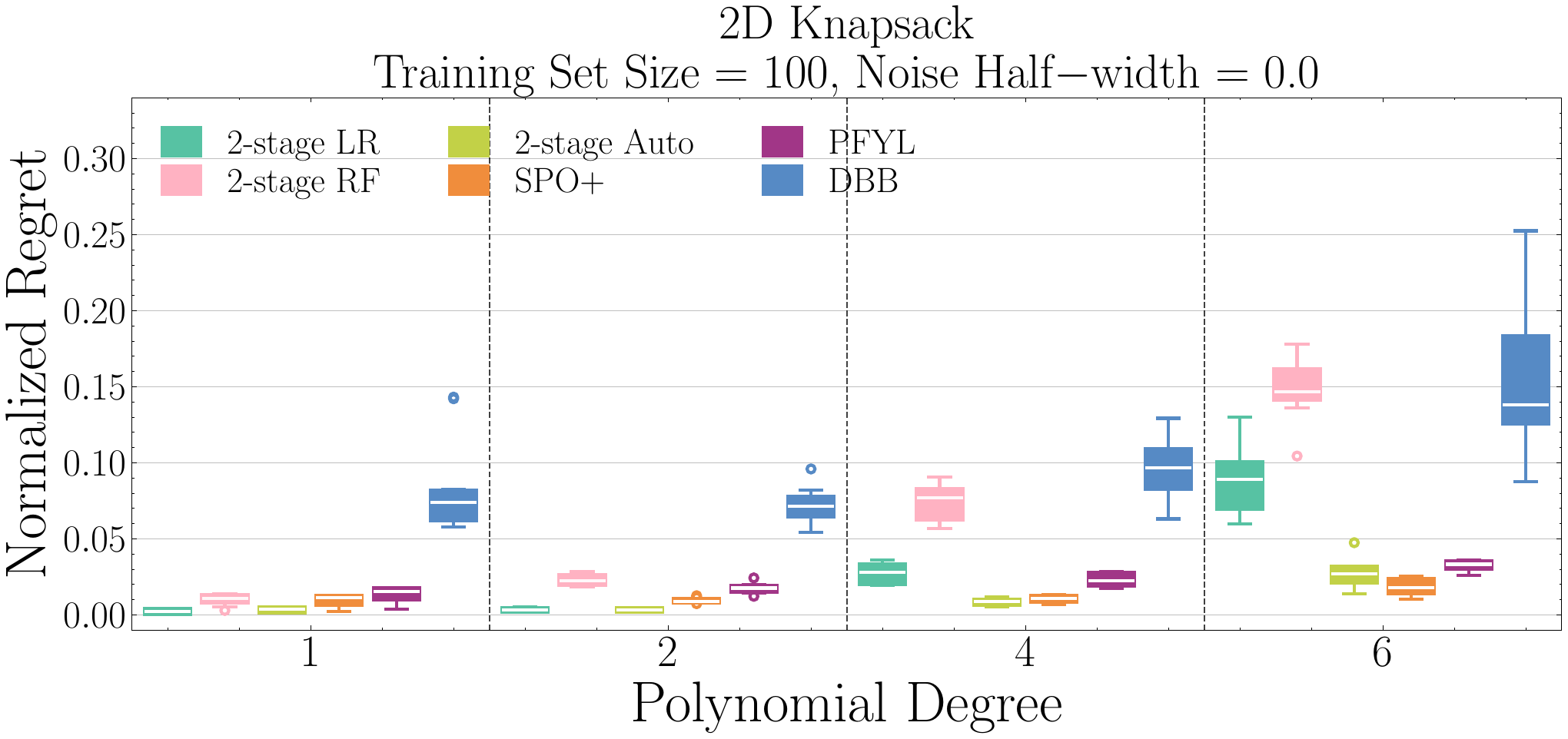}} 
    \subfloat[]{\includegraphics[width=0.48\textwidth]{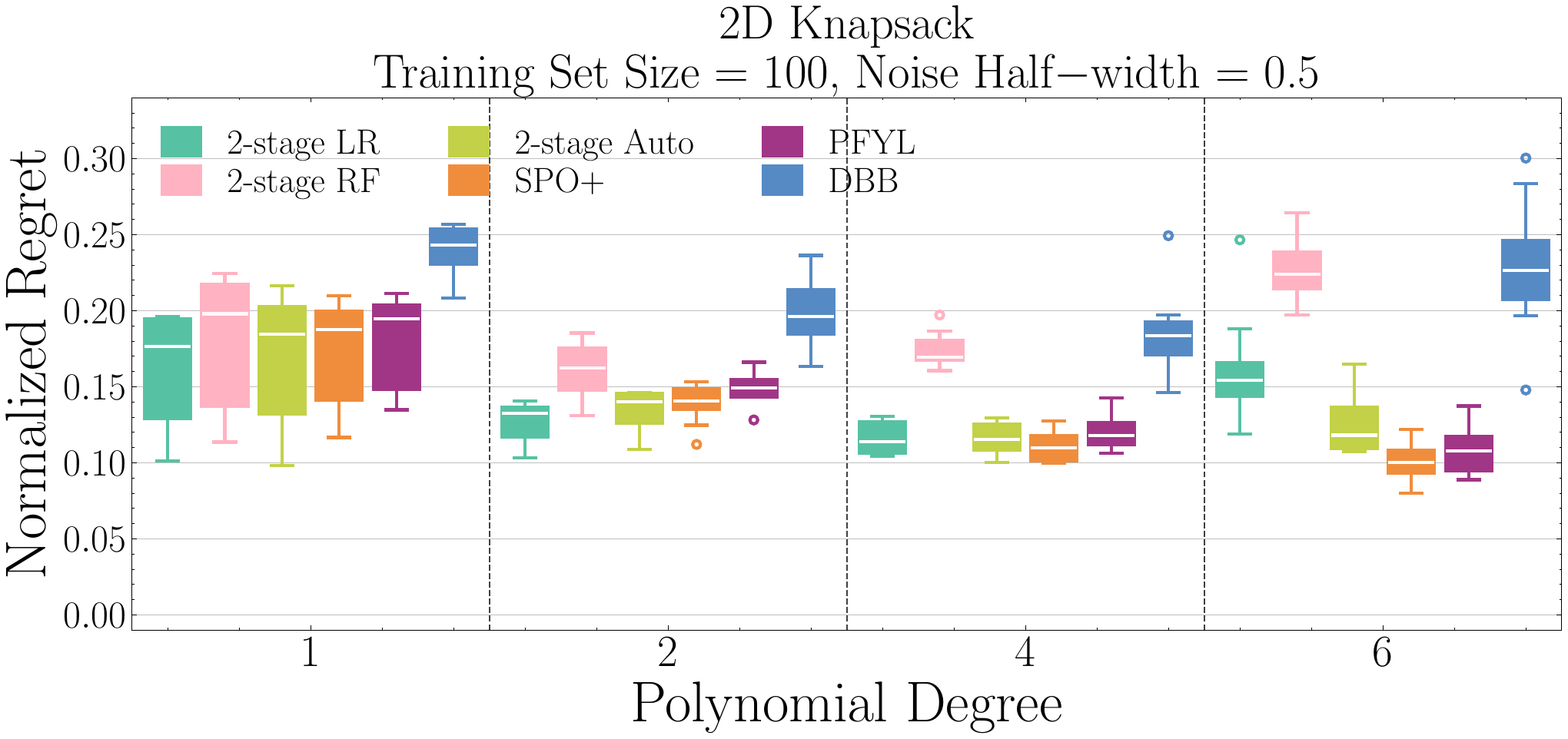}} \vspace{-2\baselineskip}
    \subfloat[]{\includegraphics[width=0.48\textwidth]{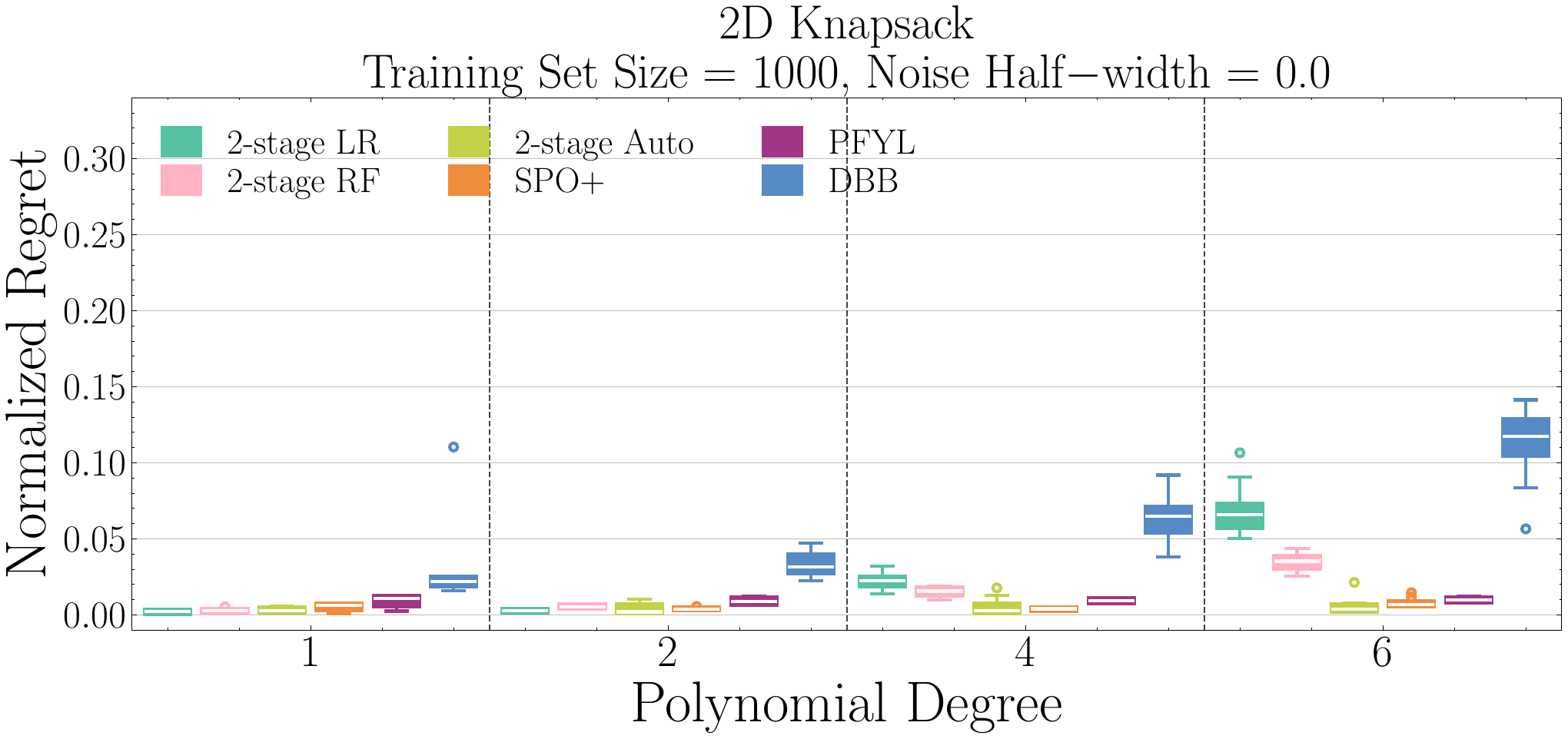}} 
    \subfloat[]{\includegraphics[width=0.48\textwidth]{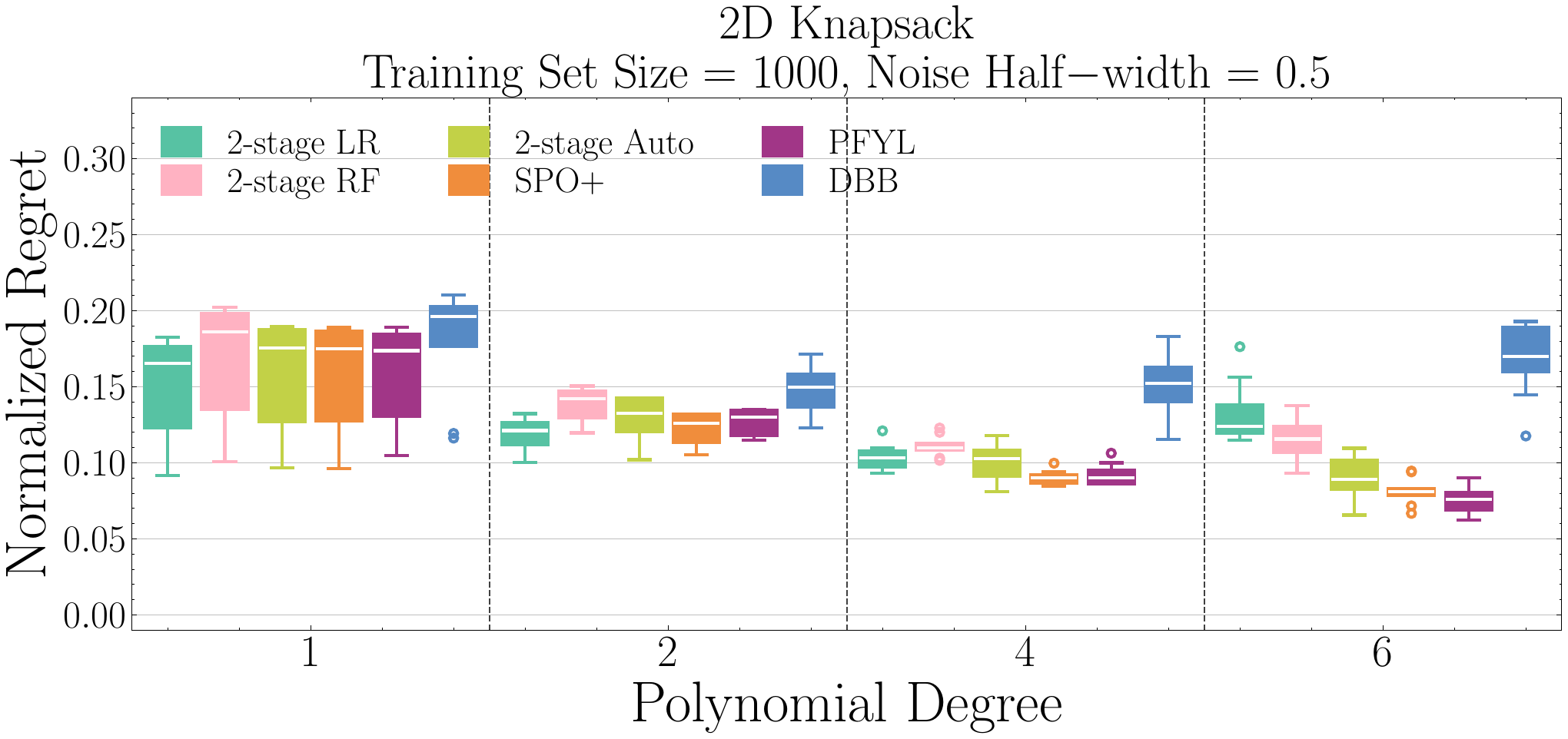}} \vspace{-2\baselineskip}
    \subfloat[]{\includegraphics[width=0.48\textwidth]{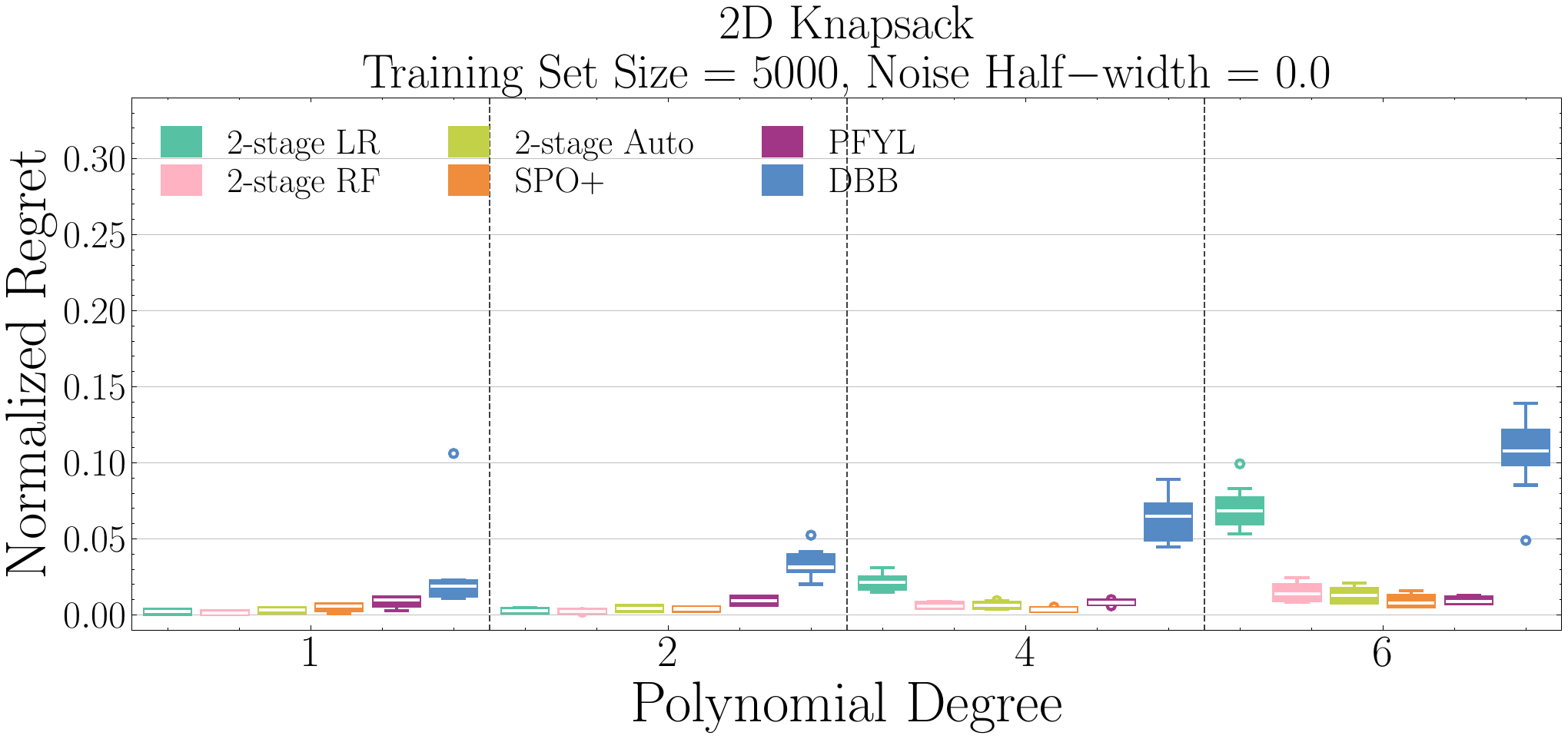}} 
    \subfloat[]{\includegraphics[width=0.48\textwidth]{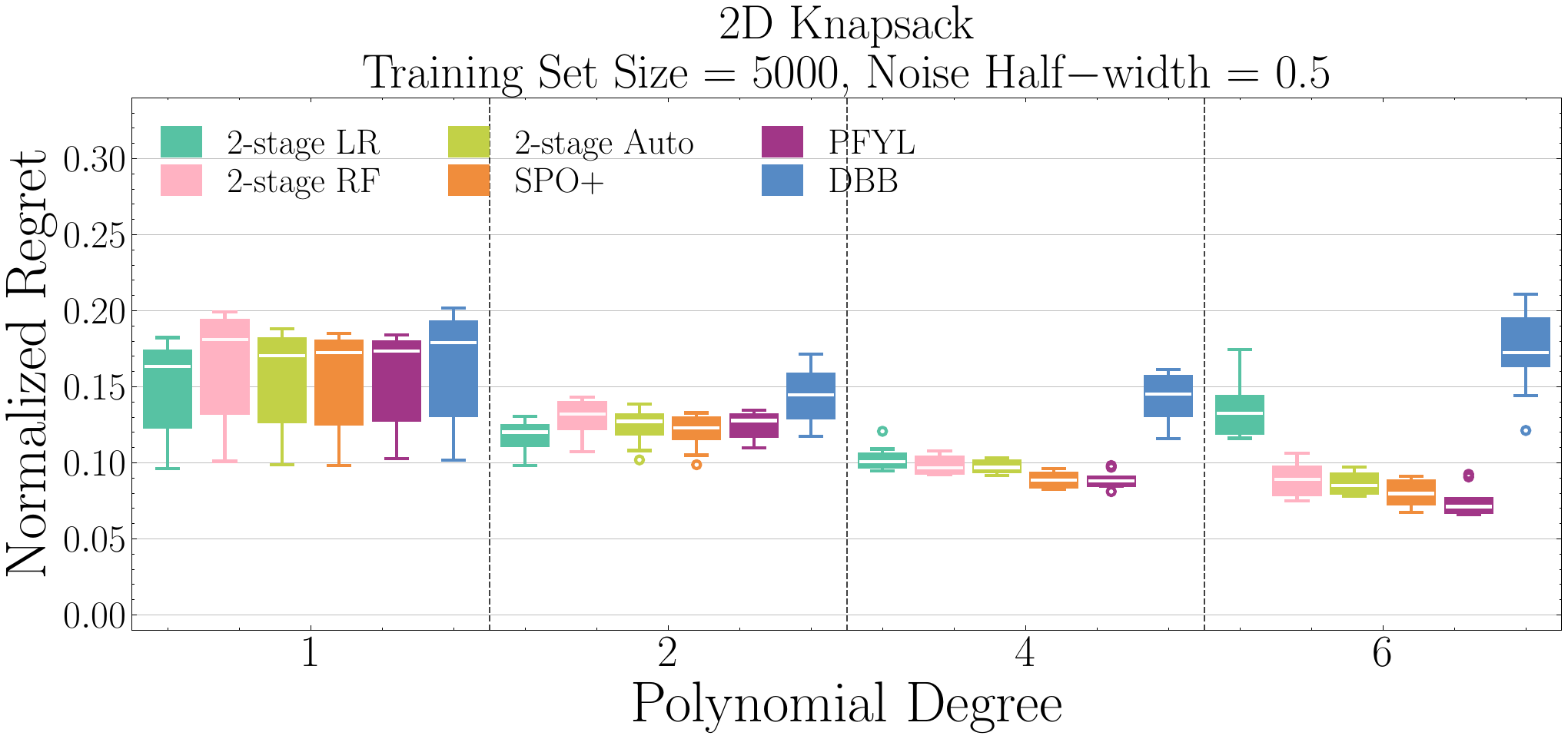}} \\
    \caption{Normalized regret for the 2D knapsack problem on the test set: There are $32$ items, and the capacity of the two resources is $20$. The methods in the experiment include two-stage approaches with linear regression, random forest, and \autosklearn{} and end-to-end learning such as \spo{}, \pfyl{} and \dbb{}. The normalized regret is visualized under different sample sizes, noise half-width, and polynomial degrees. For normalized regret, lower is better.}
    \label{fig:ks}
\end{figure}

\begin{figure}[htbp]
    \centering
    \subfloat[]{\includegraphics[width=0.48\textwidth]{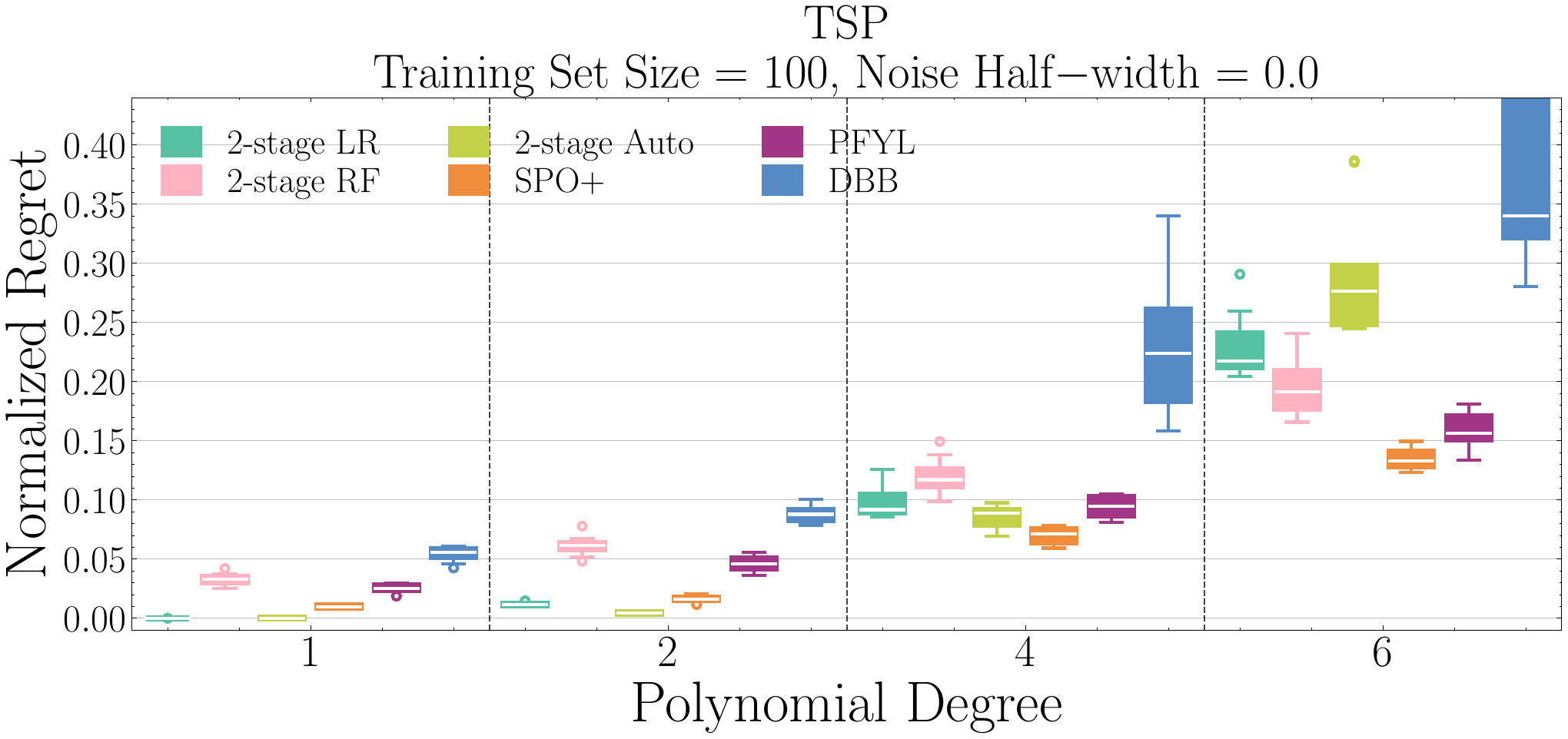}} 
    \subfloat[]{\includegraphics[width=0.48\textwidth]{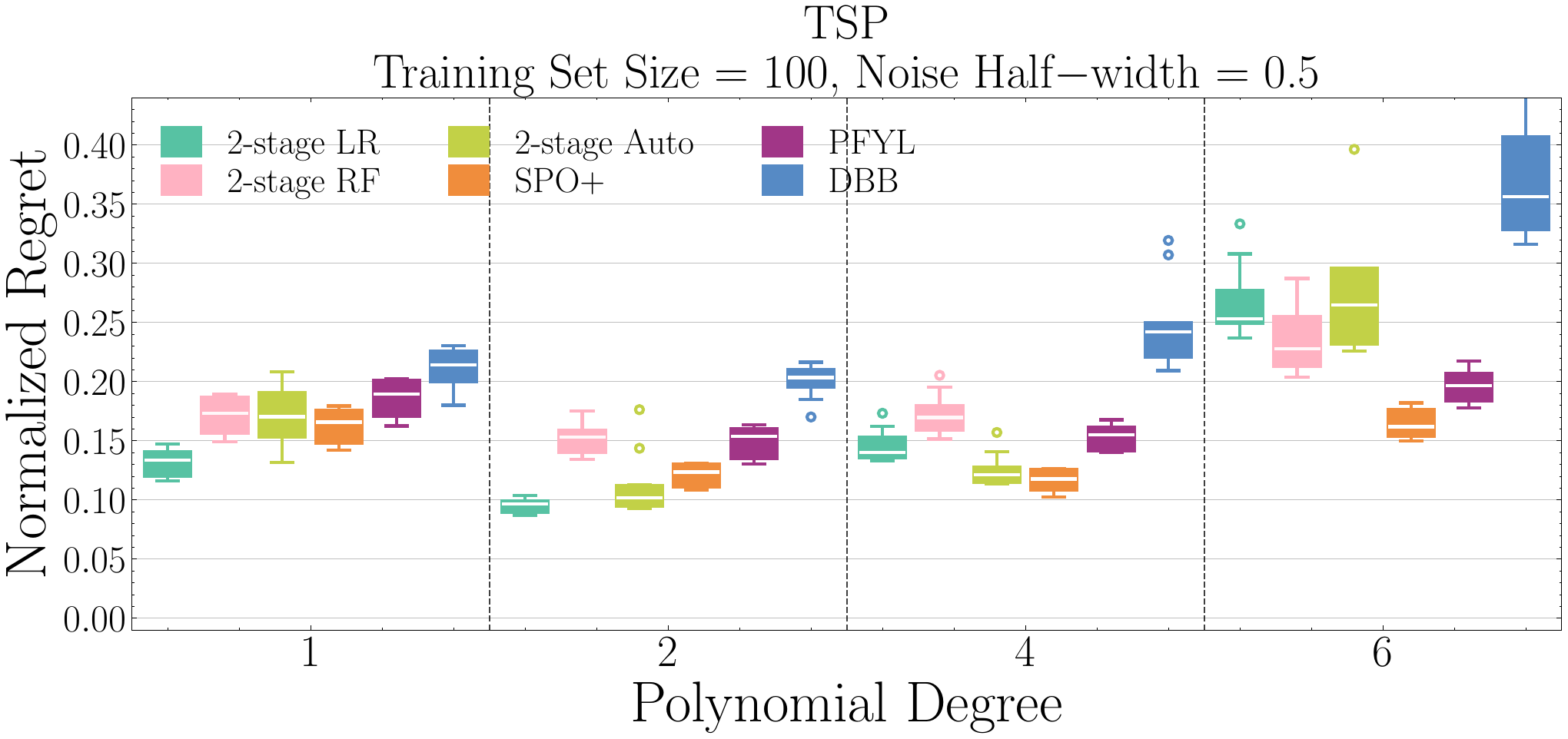}} \vspace{-2\baselineskip}
    \subfloat[]{\includegraphics[width=0.48\textwidth]{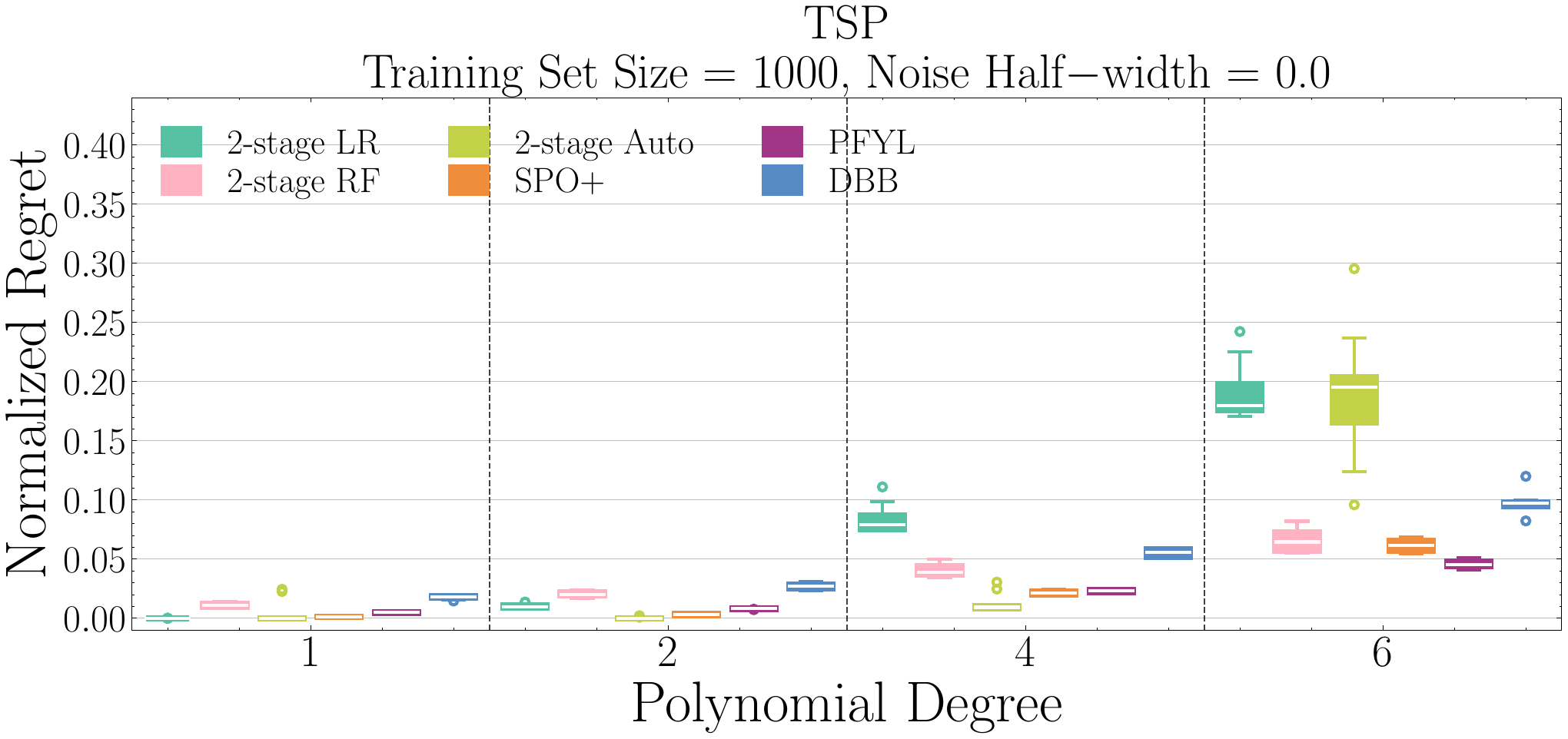}} 
    \subfloat[]{\includegraphics[width=0.48\textwidth]{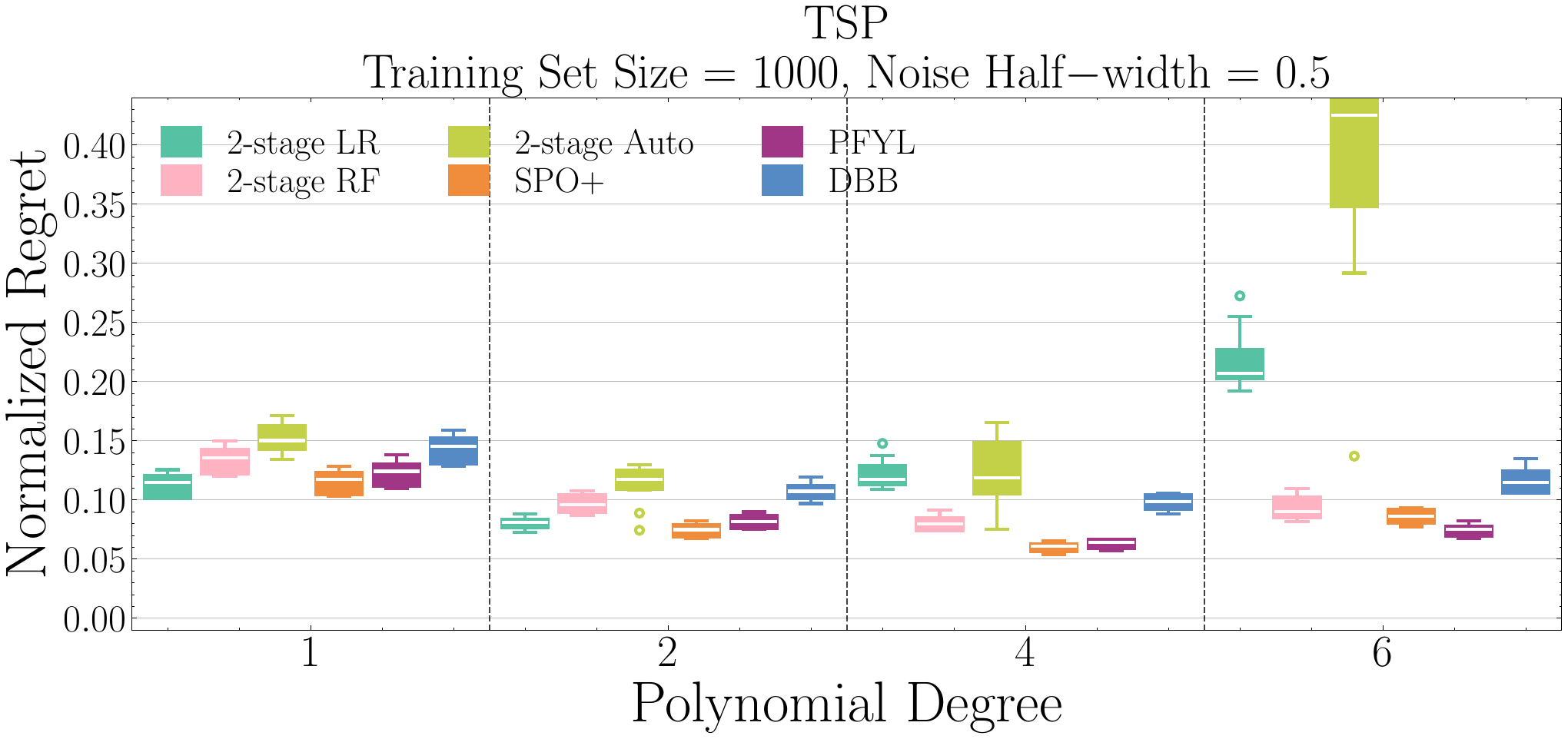}} \vspace{-2\baselineskip}
    \subfloat[]{\includegraphics[width=0.48\textwidth]{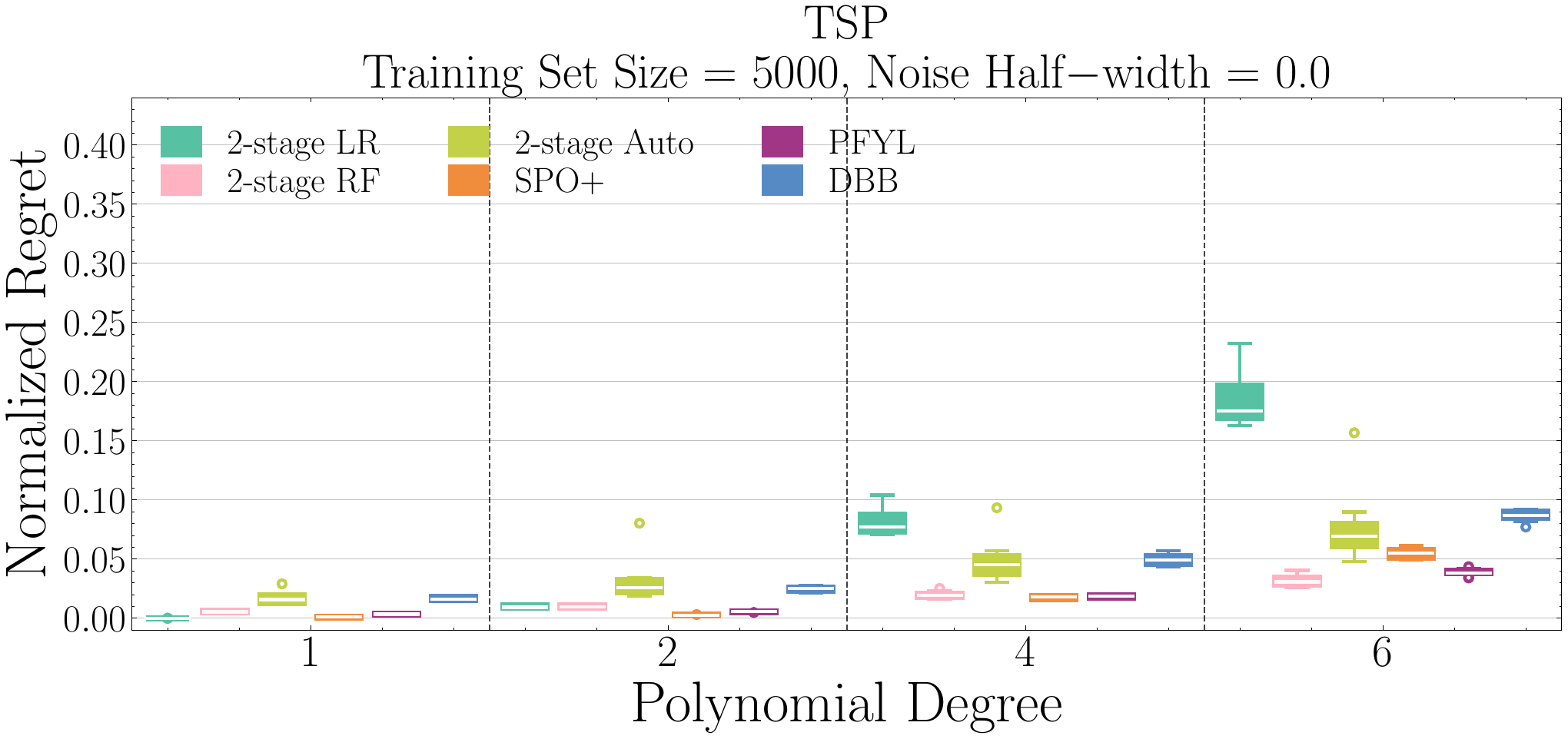}} 
    \subfloat[]{\includegraphics[width=0.48\textwidth]{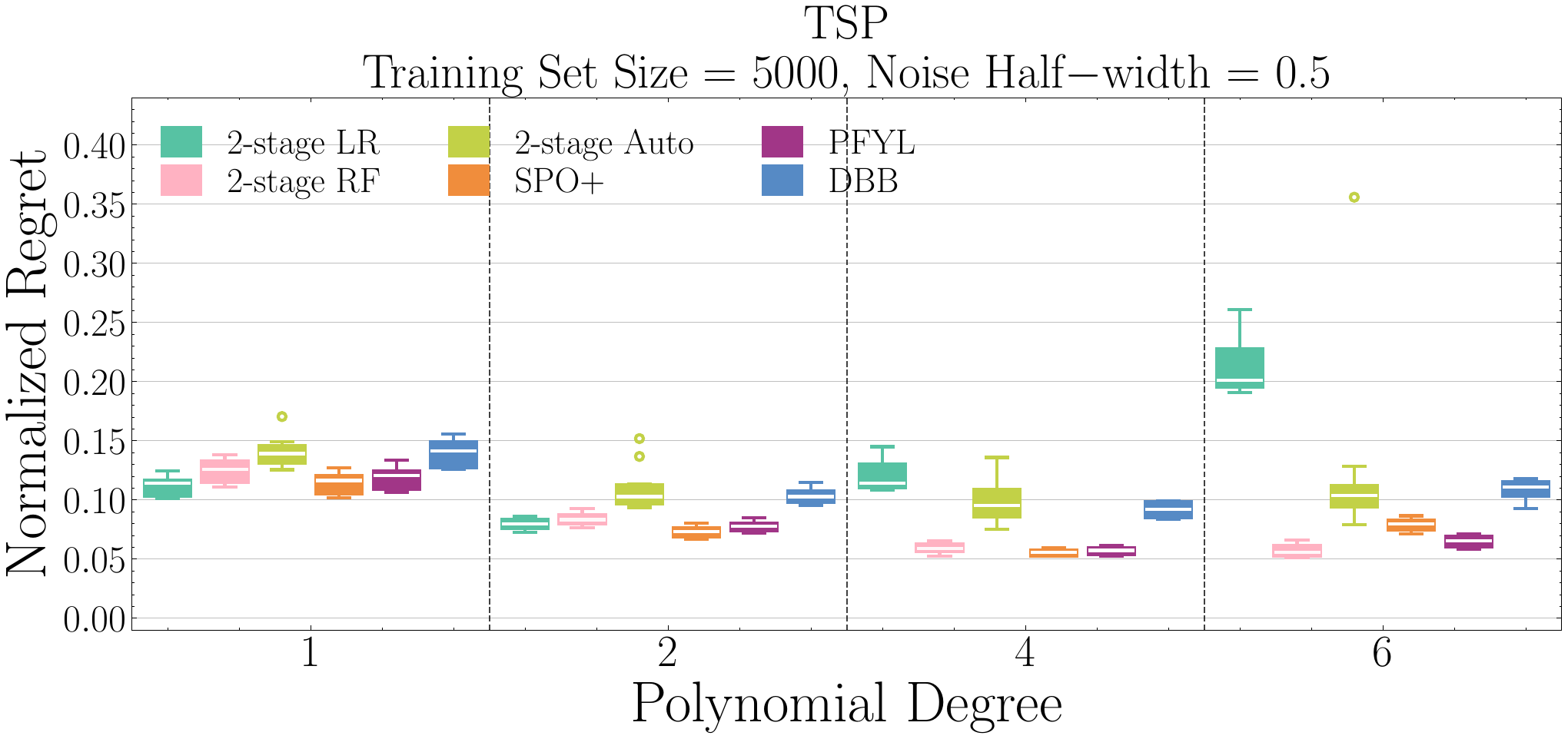}} \\
    \caption{Normalized regret for the TSP problem on the test set: There are $20$ nodes to visit. The methods in the experiment include two-stage approaches with linear regression, random forest, and \autosklearn{} and end-to-end learning such as \spo{}, \pfyl{}, and \dbb{}. The normalized regret is visualized under different sample sizes, noise half-width, and polynomial degrees. For normalized regret, lower is better.}
    \label{fig:tsp}
\end{figure}

Figures \ref{fig:sp}, \ref{fig:ks}, and \ref{fig:tsp} summarize the performance comparison for the shortest path problem, the 2D knapsack problem, and the traveling salesperson problem. These figures should be interpreted as follows: The left column is for noise-free coefficients (easier), while the right column includes noise. Each row of figures is for a training set size in increasing order. Within each figure, the degree of the polynomial that generates the coefficients from the feature vector increases from left to right. For each such polynomial degree, the different methods' boxplots are shown, summarizing the test set normalized regret results across the $10$ different experiments; lower is better.

The two-stage linear regression (2-stage LR) performs well at lower polynomial degrees but loses its advantage at higher polynomial degrees. The two-stage random forest (2-stage RF) is robust at high polynomial degrees but requires a large amount of training data. With $5000$ data samples, the random forest achieves the best performance in many cases. 
The two-stage method with automated hyperparameter tuning using the \autosklearn{} tool~\cite{feurer-neurips15a} (discussed further in Section~\ref{subsec:auto}) is notable: despite tuning for lower prediction error (not decision error), \autosklearn{} effectively reduces the regret so that it usually performs better than two-stage linear regression and random forest. However, as the output dimensions increase (e.g. TSP in Figure~\ref{fig:tsp}), \autosklearn{} struggles to remain competitive.

\spo{} and \pfyl{} show their advantage: they perform best, or at least relatively well, in all cases. These methods are comparable to linear regression at low polynomial degrees and depend less on the sample size than a random forest. At high polynomial degrees, \spo{} and \pfyl{}  outperform \autosklearn{}, exposing the limitations of the two-stage approach. Furthermore, the \pfyl{} method offers the benefit of not requiring the presence of true coefficients $\bm{c}$ within the training dataset, setting it apart from \spo{}.

%Although the results of \dbb{} are not ideal, \citet{poganvcic2019differentiation} demonstrated in their experiments that \dbb{} can operate on complex image inputs using convolutional neural networks for shortest path and TSP, showing the ability to extract features from complex raw data; future additions to our package will include such experiments.

\begin{tcolorbox}[colback=white,colframe=gray,title=Finding \#1]
    \spo{} and \pfyl{} can robustly achieve relatively good decisions under different scenarios, often outperforming two-stage baselines.
\end{tcolorbox}

\subsection{Two-stage Method with Automated Hyperparameter Tuning}
\label{subsec:auto}

This method leverages the sophisticated \autosklearn{}~\cite{feurer-neurips15a} tool that uses Bayesian optimization methods for the automated hyperparameter tuning of the~\sklearn{} regression models. The metric of ``2-stage Auto" is the mean squared error of the predicted coefficients, which does not directly reduce the decision error. Due to the limitation of multi-output regression in \autosklearn{} v0.14.6, the choices of the predictor in 2-stage Auto are restricted to five models: k-nearest neighbor (KNN), decision tree, random forest, extra trees, and Gaussian process. Despite these constraints, \autosklearn{} can achieve low regret. Although 2-stage Auto training is time-consuming, it remains a competitive method.

%When we examine the selection of the 2-stage Auto predictor, the results have some similarities. Compared to the other models, decision trees are never the ideal choice. A fine-tuned Gaussian process has impressive performance for datasets with small samples, no noise, or high polynomial degree, while KNN and random forest can also be the best models at times.

\begin{tcolorbox}[colback=white,colframe=gray,title=Finding \#2]
    Even with successful model selection and hyperparameter tuning, the two-stage method performs worse than \spo{} in terms of decision quality, underscoring the value of the end-to-end approach.
\end{tcolorbox}

\subsection{Exact Method and Relaxation}
\label{subsec:rel} 

Training \spo{}, \pfyl{}, and \dbb{} with a linear relaxation instead of solving the integer program improves the computational efficiency. However, the use of a ``weaker" solver theoretically undermines the decision quality. Therefore, an important question arises about the trade-off when using linear relaxation in training. To this end, we compare the performance of end-to-end approaches with their relaxation using the 2D knapsack and TSP as case studies. We ensure consistency by using the same instances, models, and hyperparameters as in our previous experiments.

\begin{figure}[htbp]
    \centering
    \subfloat[]{\includegraphics[width=0.48\textwidth]{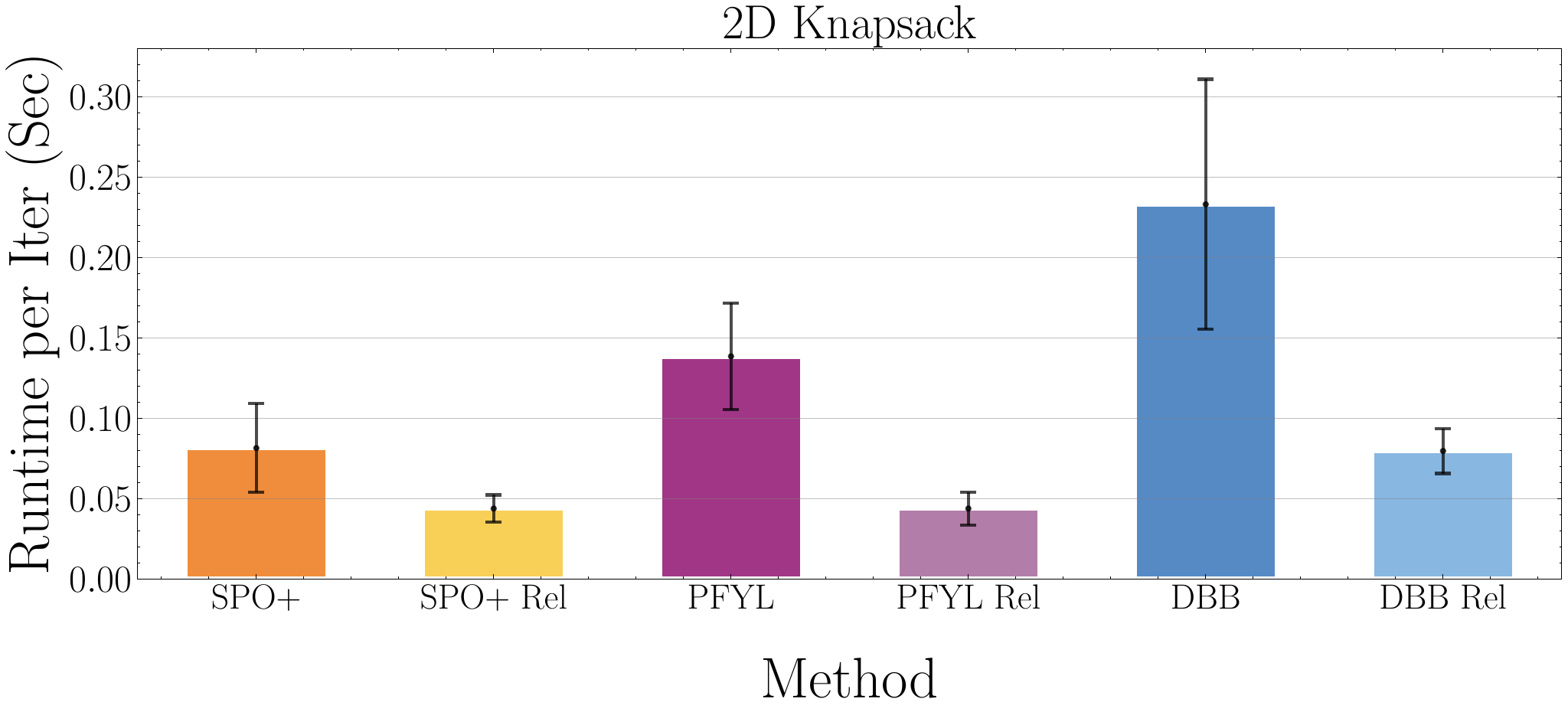}} 
    \subfloat[]{\includegraphics[width=0.48\textwidth]{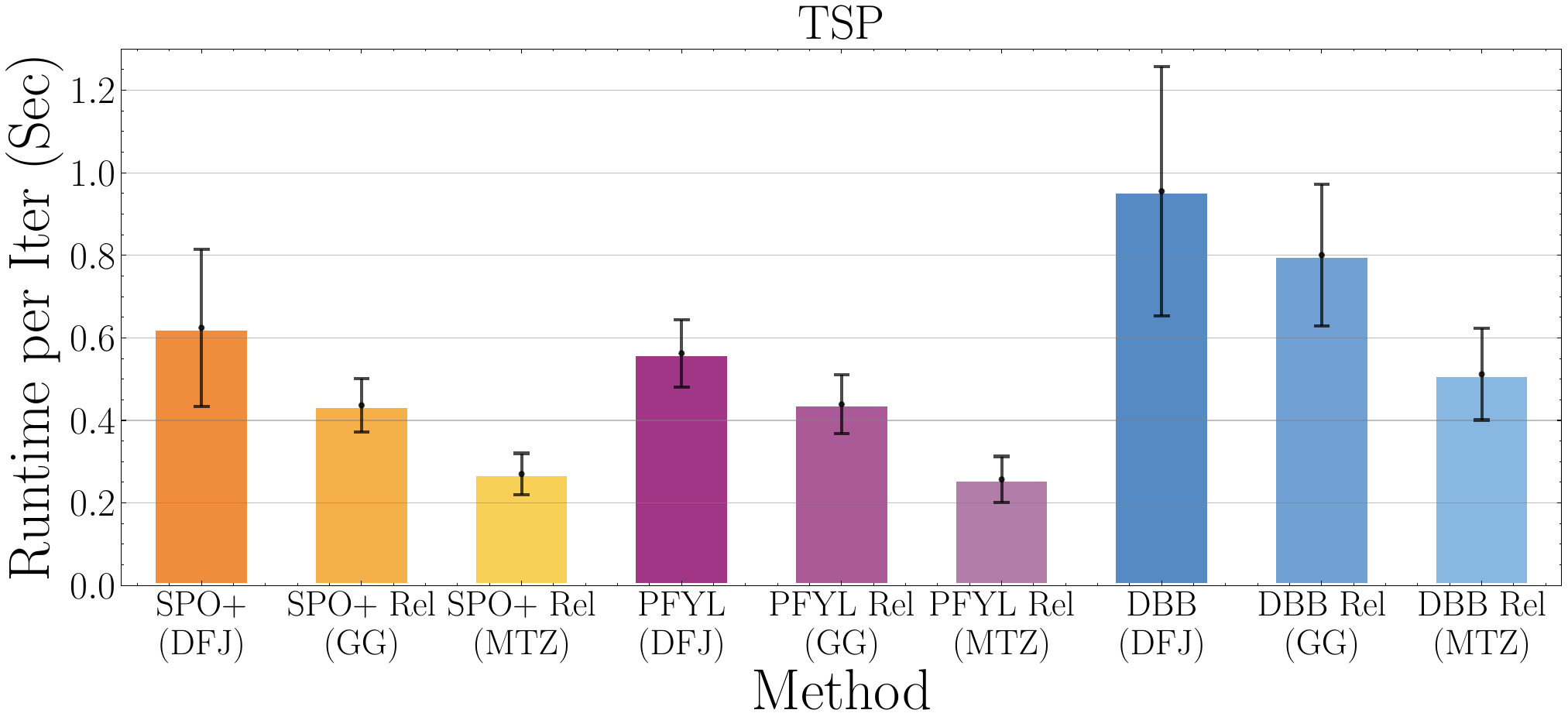}} \\
    \caption{Average training time per iteration for exact and relaxation methods with standard deviation error bars: We visualized the mean training time with standard deviation for the Knapsack and TSP; lower is better.}
    \label{fig:rel-time}
\end{figure}

There are several integer programming formulations for the TSP. In addition to DFJ, we have also implemented the Miller-Tucker-Zemlin (MTZ) formulation \cite{miller1960integer} and the Gavish-Graves (GG) formulation \cite{gavish1978travelling}. Although all formulations yield the same integer solution, they differ in computational efficiency and linear relaxations. Since DFJ requires column generation to handle exponential subtour constraints, its linear relaxation is hard to implement. The GG formulation is shown to have a tighter linear relaxation than MTZ. Thus, we use DFJ for exact \spo{}, \pfyl{}, and \dbb{}, and MTZ and GG for relaxation to investigate the effect of the quality of the solution on regret.

\begin{figure}[htbp]
    \centering
    %\subfloat[]{\includegraphics[width=0.48\textwidth]{img/rel-ks2-n100e0.pdf}} 
    \subfloat[]{\includegraphics[width=0.48\textwidth]{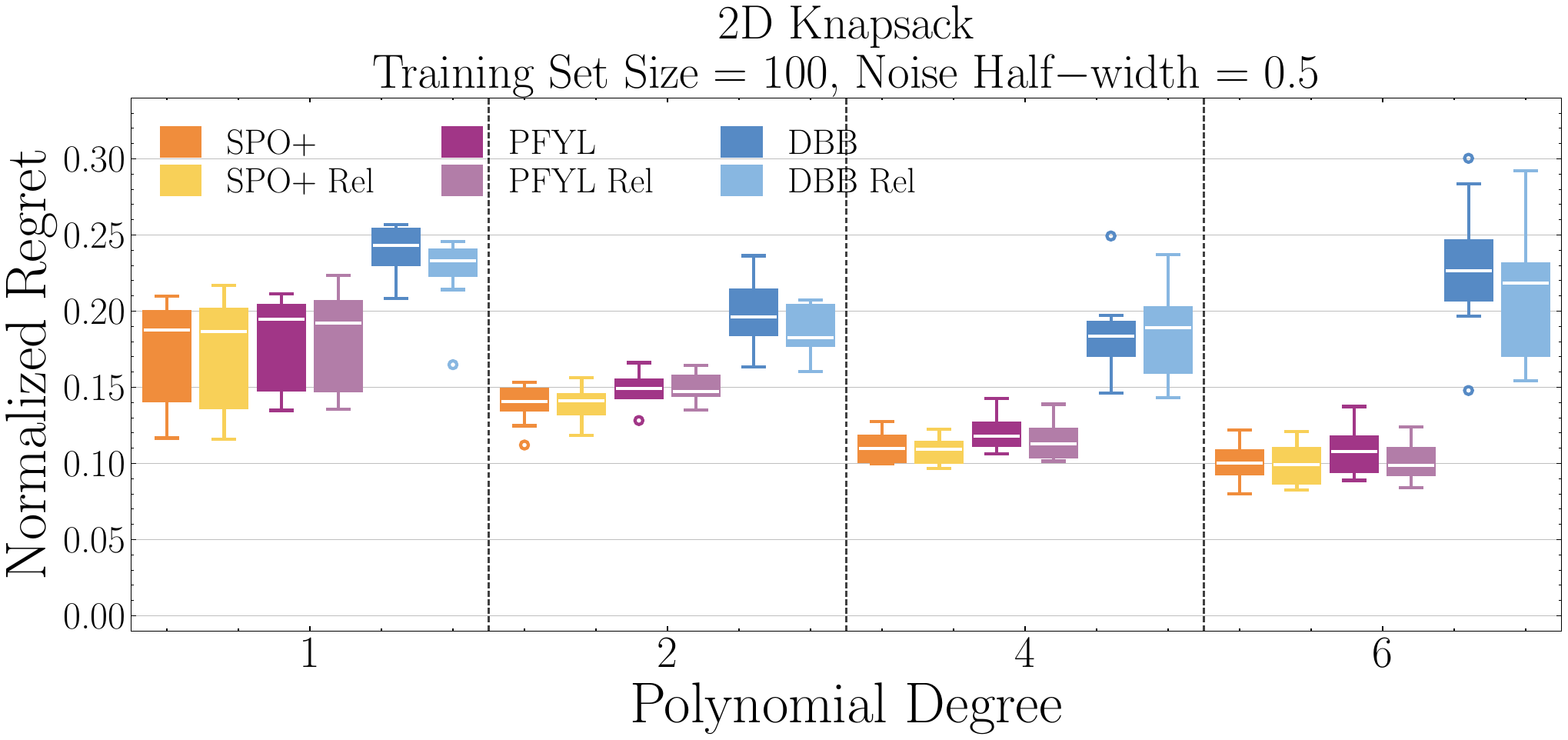}}
    %\subfloat[]{\includegraphics[width=0.48\textwidth]{img/rel-ks2-n1000e0.pdf}} 
    \subfloat[]{\includegraphics[width=0.48\textwidth]{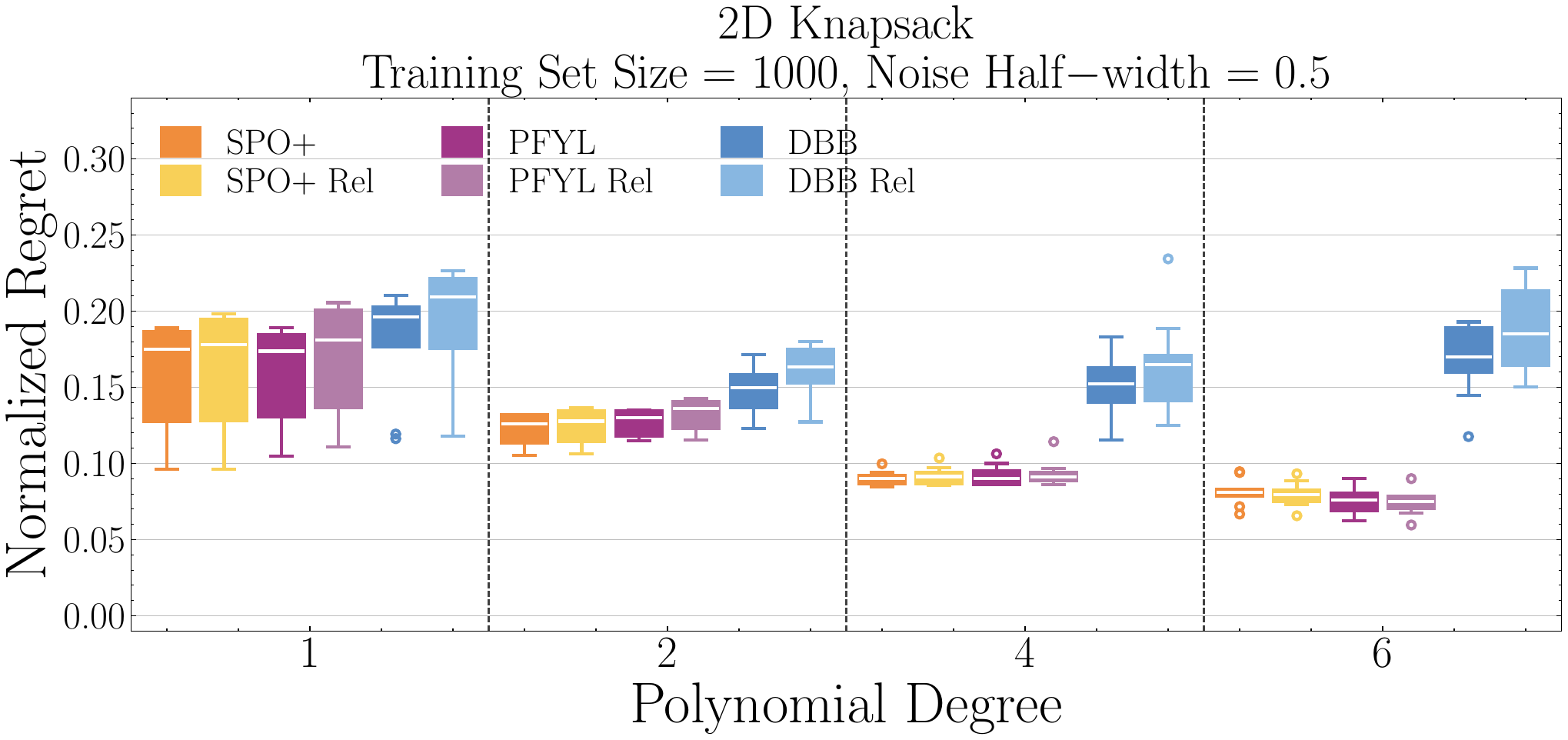}} \vspace{-2\baselineskip}
    %\subfloat[]{\includegraphics[width=0.48\textwidth]{img/rel-tsp-n100e0.pdf}} 
    \subfloat[]{\includegraphics[width=0.48\textwidth]{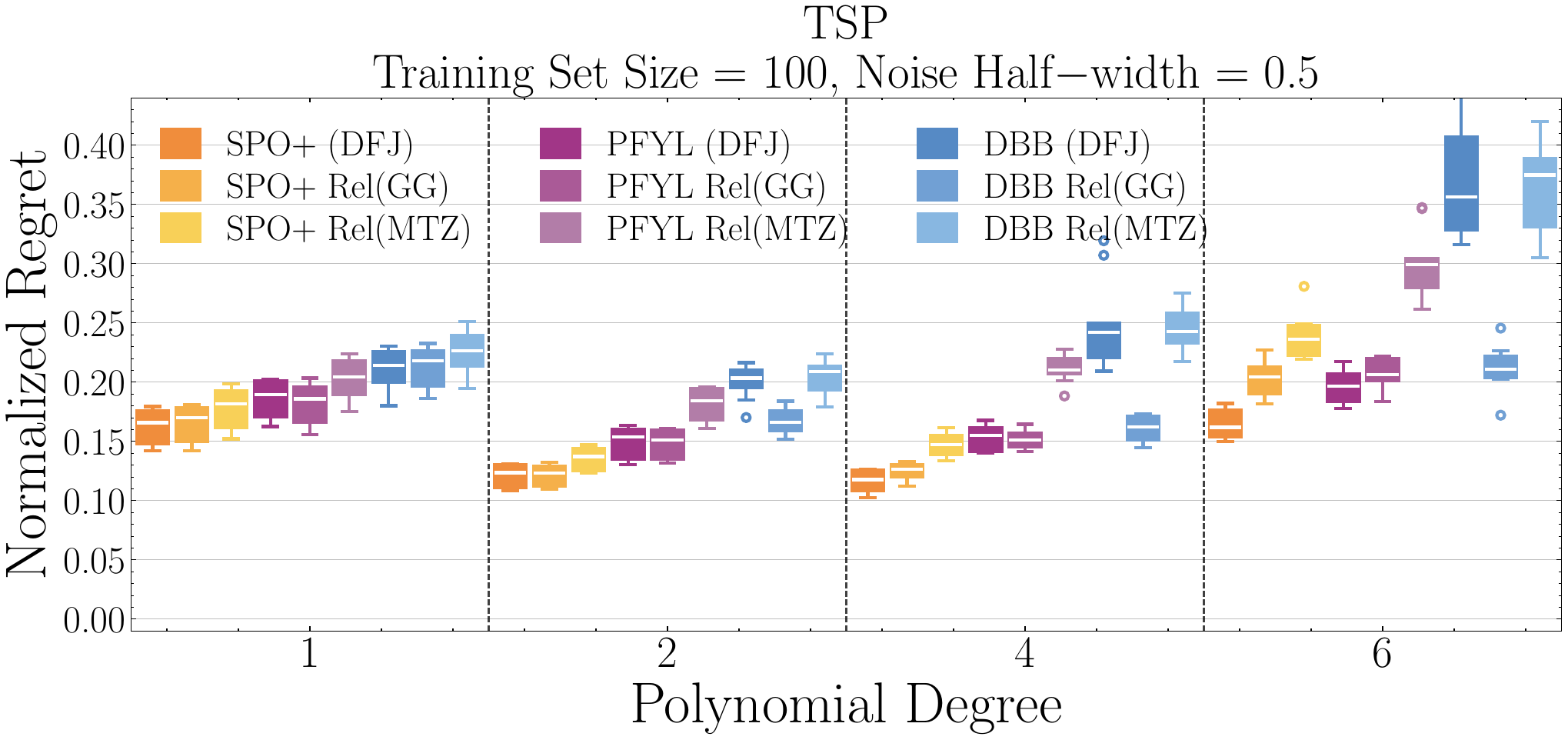}}
    %\subfloat[]{\includegraphics[width=0.48\textwidth]{img/rel-tsp-n1000e0.pdf}} 
    \subfloat[]{\includegraphics[width=0.48\textwidth]{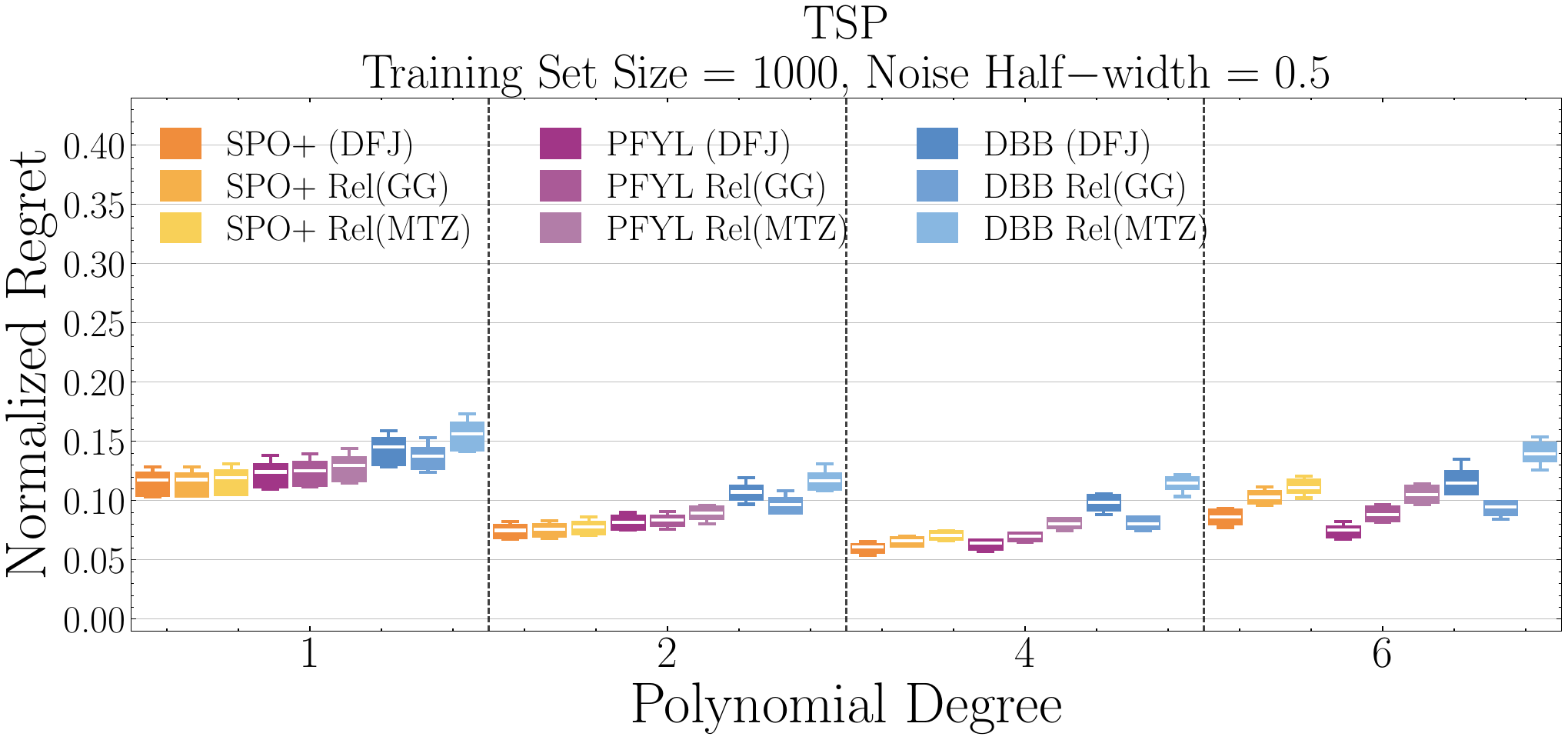}} \\
    \caption{Normalized regret for the 2D knapsack (at the top) and TSP (at the bottom) on the test set:  The methods in the experiment include \spo{}, \pfyl{} and \dbb{} w/o relaxation. Then, we visualize the normalized regret under different sample sizes and polynomial degrees to investigate the impact of the relaxation method. For normalized regret, lower is better.}
    \label{fig:rel}
\end{figure}

According to Figure \ref{fig:rel-time}, the use of linear relaxation significantly reduces training time. As shown in Figure \ref{fig:rel}, the impact on the performance of the knapsack is almost negligible. Interestingly, relaxed methods have the potential to perform better on small data, perhaps because linear relaxation acts as a regularization mechanism that prevents overfitting. For TSP, Figure \ref{fig:rel} demonstrates that a tighter bound does reduce the regret more, and \dbbrel{} with GG formulation shows advantages over \dbb{}. Overall, employing relaxations achieves commendable performance with enhanced computational efficiency. Moreover, formulations with tighter linear relaxation lead to better performance.

\begin{tcolorbox}[colback=white,colframe=gray,title=Finding \#3]
    End-to-end predict-then-optimize with relaxation has excellent potential to improve computation efficiency at a slight degradation in performance, particularly when the true coefficients-generating function is not very nonlinear. 
\end{tcolorbox}

\subsection{Prediction Regularization}
\label{subsec:reg}

As proposed in~\citet{elmachtoub2021smart}, the mean absolute error $\mathcal{L}_{\text{MAE}}(\hat{\bm{c}}, \bm{c}) = \frac{1}{n} \sum_i^n{\norm{\hat{\bm{c}}^i-\bm{c}^i}_1}$ or the mean squared error $\mathcal{L}_{\text{MSE}}(\hat{\bm{c}}, \bm{c}) = \frac{1}{2n} \sum_i^n{\norm{\hat{\bm{c}}_i- \bm{c}^i}_2^2}$ of the predicted coefficient can be added to the decision loss as $l_1$ or $l_2$ regularizers. When using regularization, we set either the $l_1$ regularization parameter $\phi_{1}$ and the $l_2$ regularization parameter $\phi_{2}$ logarithmically from $0.001$ to $10$. For the experiments, we use the same instances, model, and hyperparameters as before, while the number of training samples $n$, the noise half-width $\bar{\epsilon}$, and the polynomial degree $deg$ are fixed at $1000$, $0.5$ and $4$.

\begin{figure}[htbp]
    \centering
    \subfloat[]{\includegraphics[width=0.48\textwidth]{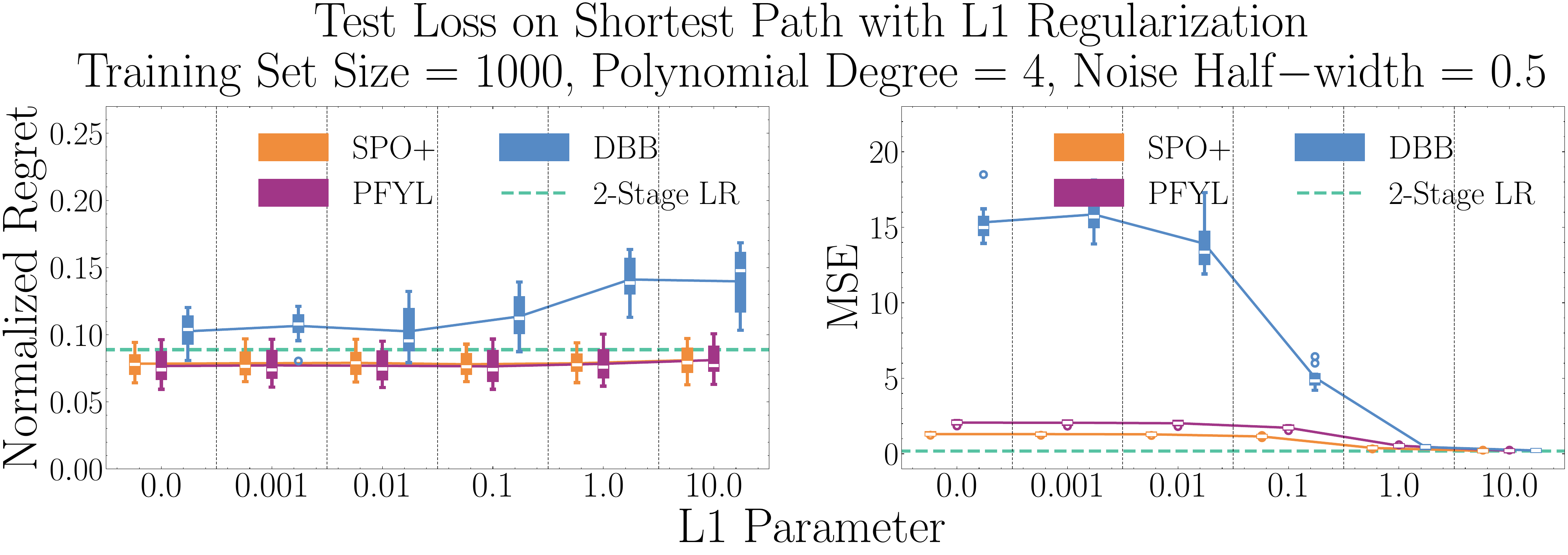}} 
    \subfloat[]{\includegraphics[width=0.48\textwidth]{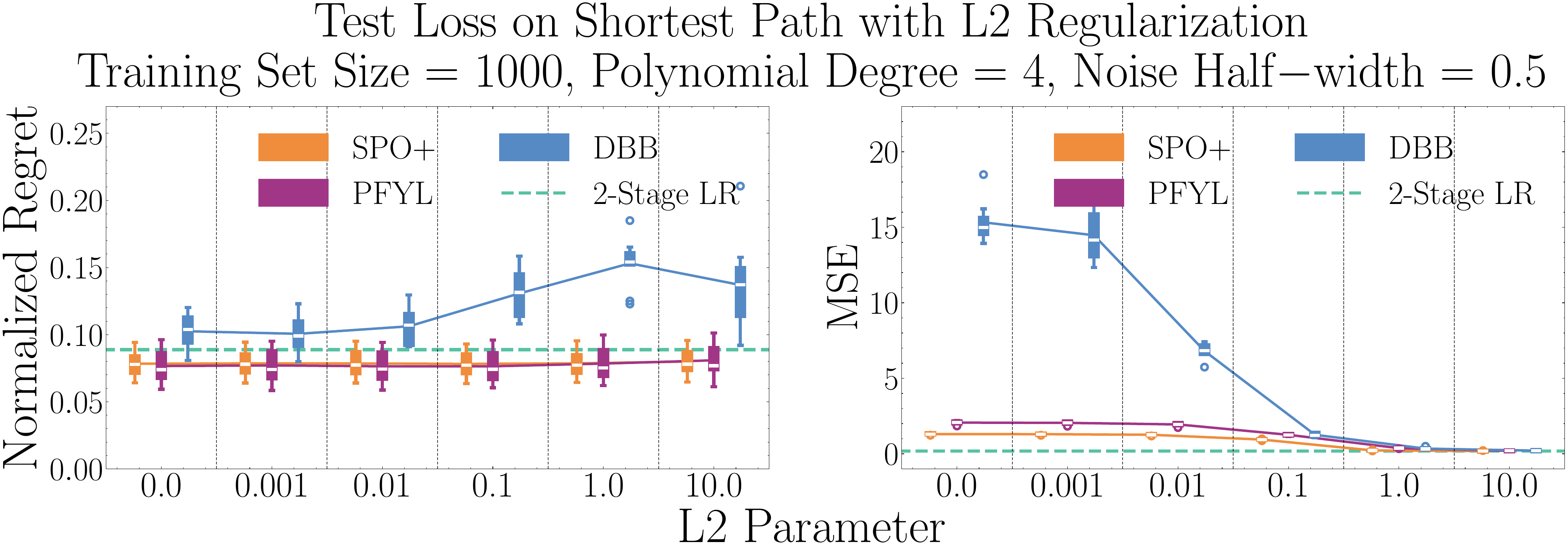}} \vspace{-2\baselineskip}
    \subfloat[]{\includegraphics[width=0.48\textwidth]{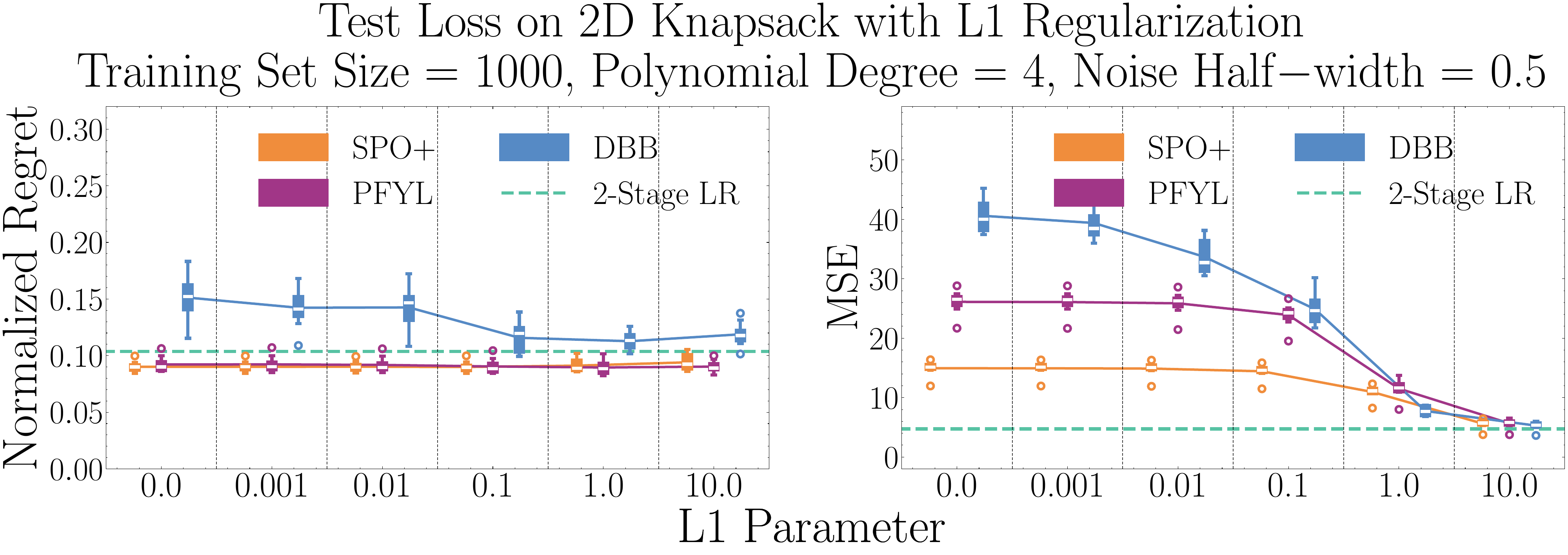}} 
    \subfloat[]{\includegraphics[width=0.48\textwidth]{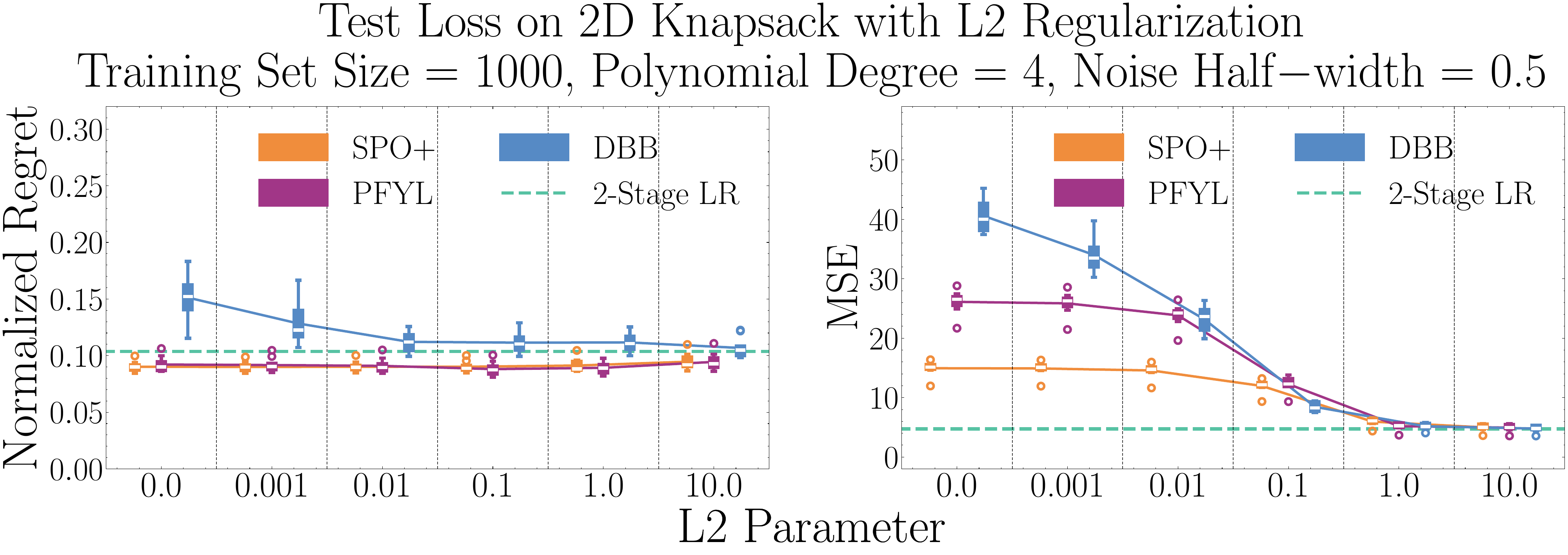}} \vspace{-1.2\baselineskip}
    \caption{Normalized regret and MSE on the test set: In these experiments, we compare the performance on \spo{}, \pfyl{}, and \dbb{} with different levels of $l_1$ or $l_2$ regularization and show the trend line with increasing regularization, while 2-Stage LR (green line) is served as the baseline. The normalized regret and MSE under different problems are shown; lower is better.}
    \label{fig:reg}
\end{figure}

Based on Figure \ref{fig:reg}, regularization appears to be a viable strategy to reduce MSE. With an increase in the regularization parameter, all end-to-end methods achieve an MSE comparable to that of the two-stage linear regression approach. In terms of regret, the impact on \spo{} and \pfyl{} is minor, while more regularization of \dbb{} has the potential to produce improvements.

\begin{tcolorbox}[colback=white,colframe=gray,title=Finding \#4]
    $l_1$ and $l_2$ regularizations are considered appealing techniques for enhancing prediction accuracy while simultaneously preserving decision quality.
\end{tcolorbox}

\subsection{Trade-offs between MSE and Regret}
\label{subsec:trade-off}

\begin{figure}[htbp]
    \centering
    \subfloat[]{\includegraphics[width=0.33\textwidth]{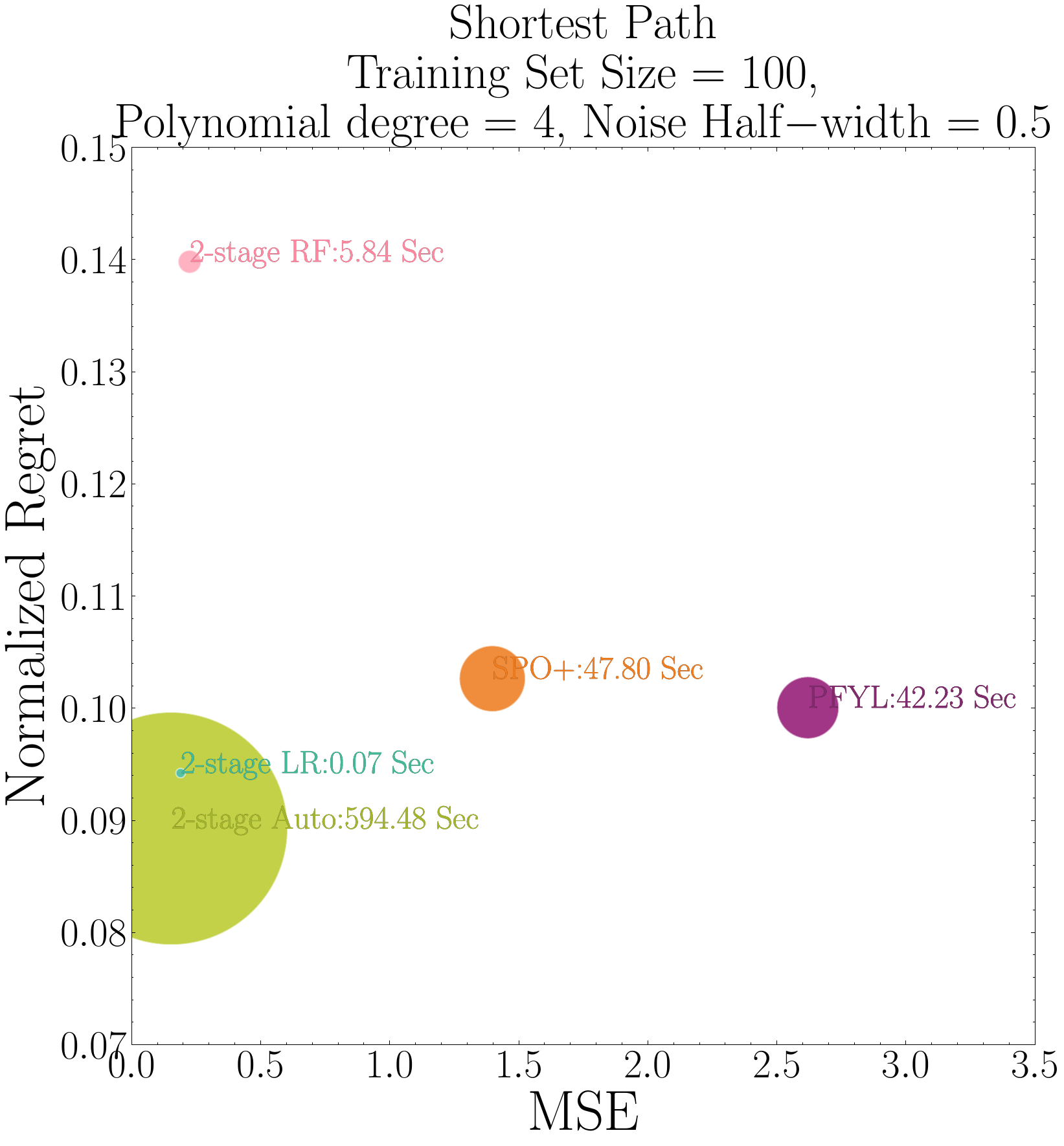}}
    \subfloat[]{\includegraphics[width=0.33\textwidth]{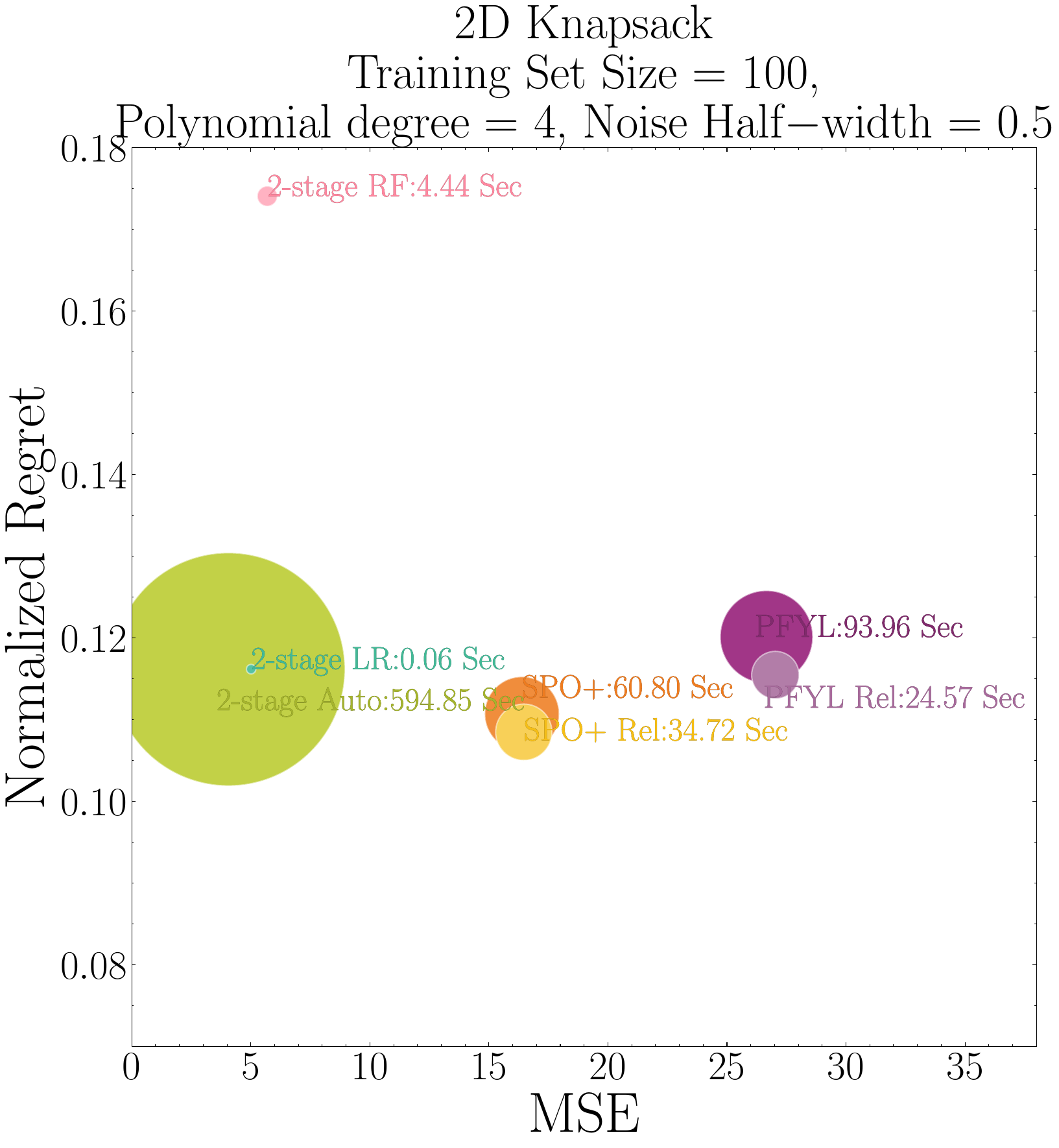}}
    \subfloat[]{\includegraphics[width=0.33\textwidth]{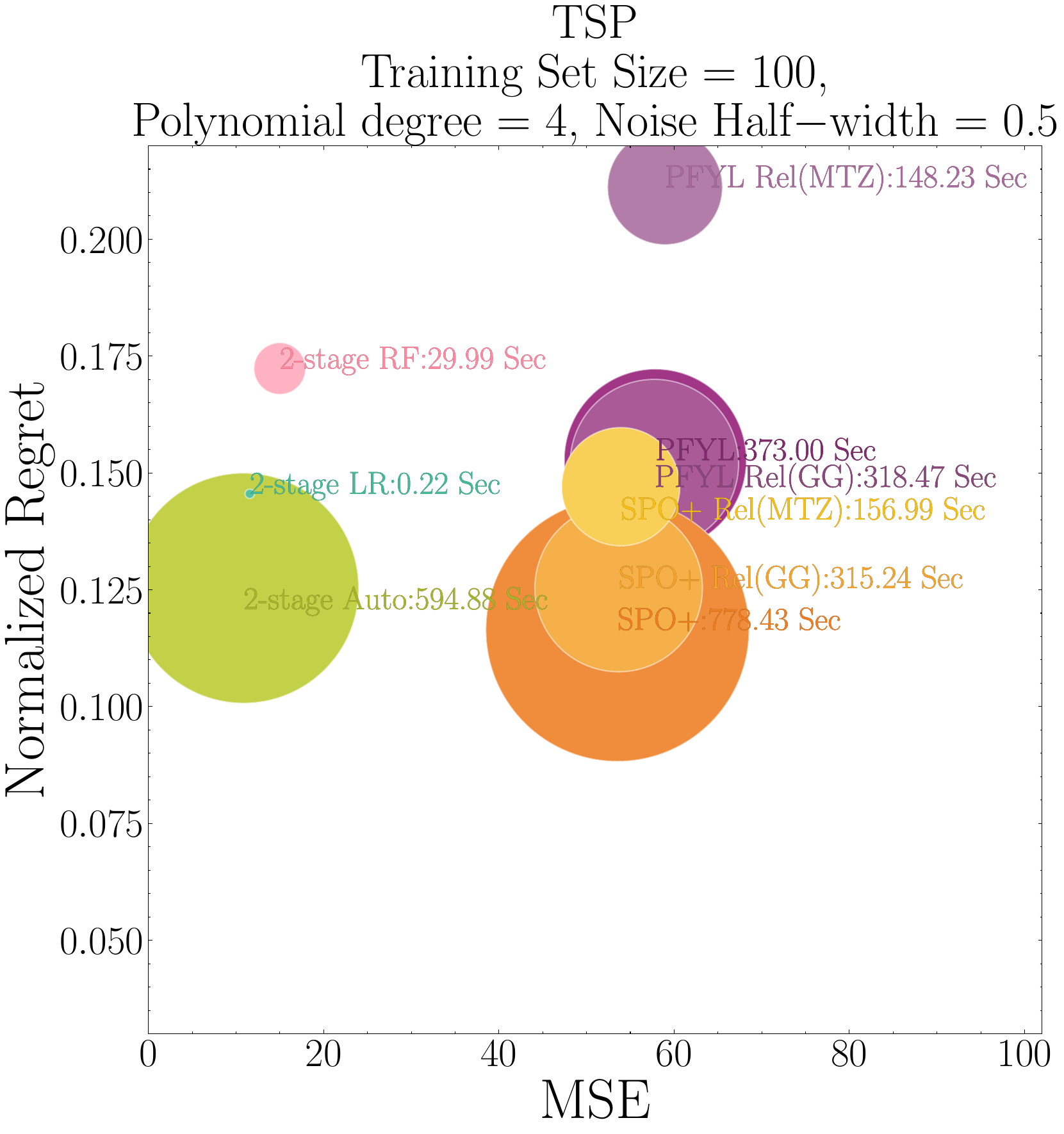}} \vspace{-2\baselineskip}
    \subfloat[]{\includegraphics[width=0.33\textwidth]{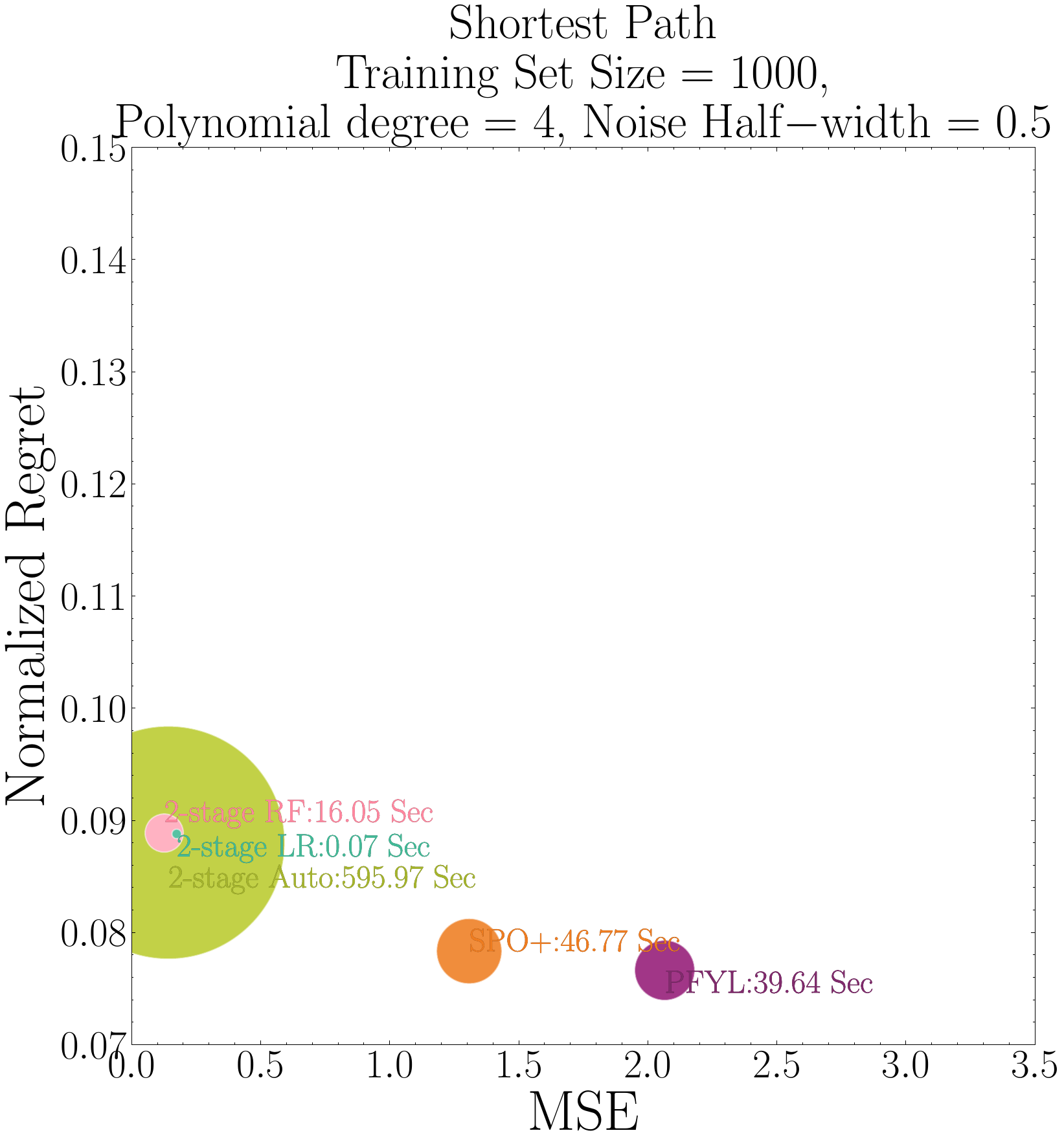}}
    \subfloat[]{\includegraphics[width=0.33\textwidth]{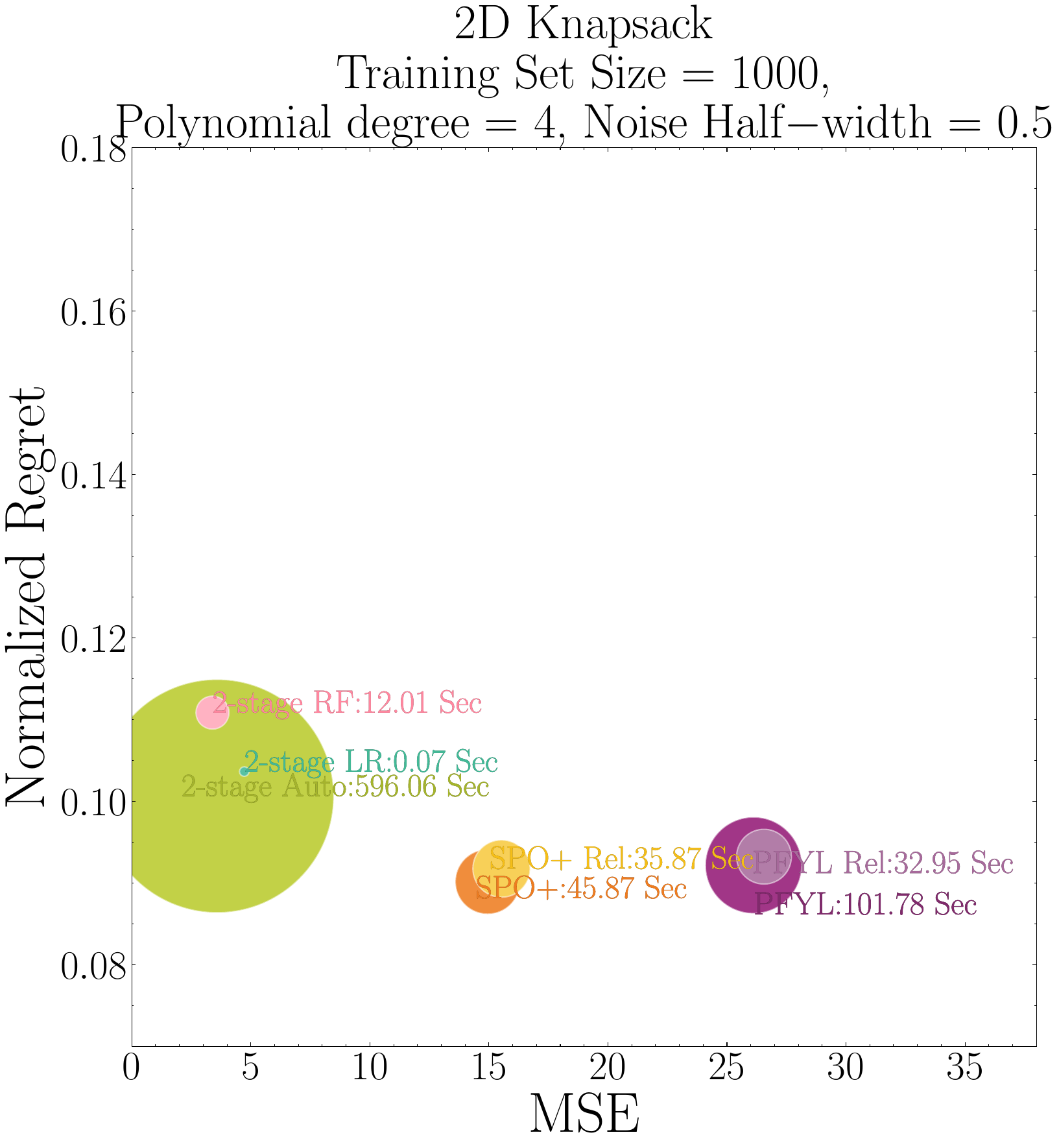}}
    \subfloat[]{\includegraphics[width=0.33\textwidth]{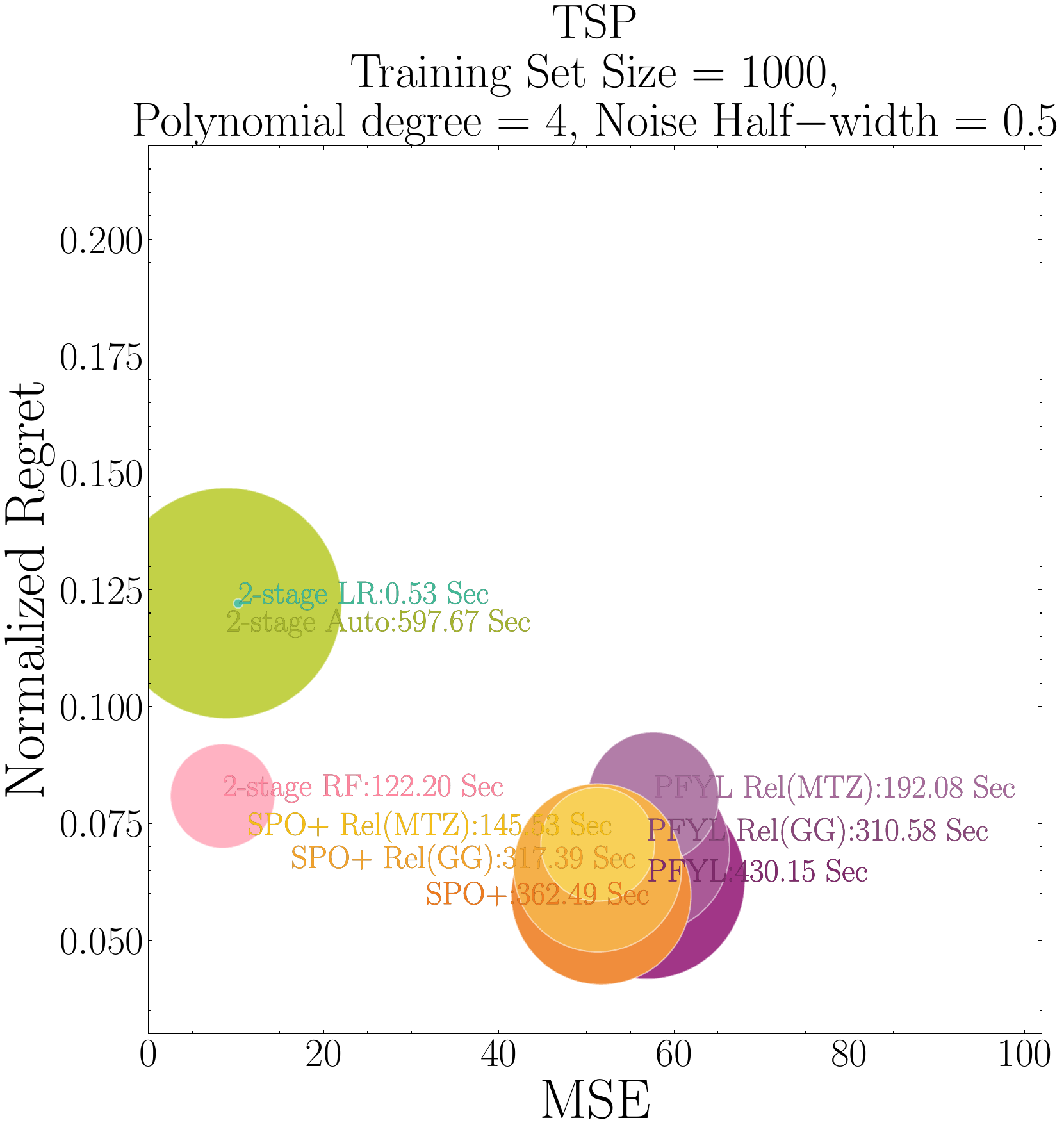}} \\
    \caption{MSE v.s. Regret: The result covers different two-stage methods, \spo{}, \pfyl{}, and their relaxations. \dbb{} is omitted because it is far from others. The size of the circles is proportional to the training time (Sec), so the smaller is better.}
    \label{fig:td}
\end{figure}

%The above experiments have focused only on regret. In addition to decision error, prediction error should also be of concern in many predict-then-optimize tasks. For example, investors seek not only the optimal portfolio but also wonder about the forecasted return of each security.
Figure \ref{fig:td} examines the trade-off between prediction and decision on the shortest path with $100$ and $1000$ training samples, a $0.5$ noise half-width, and polynomial degree $4$. We calculate the average MSE and regret over $10$ repeated random experiments. Apart from this, the mean training time is also annotated, with circle sizes proportional to the time. The circles of \dbb{} have been removed as they are located far in the top right corner and perform significantly worse than other methods.

Figure \ref{fig:td} shows that \spo{} and \pfyl{} can reach low decision errors, particularly with large training sets, at the cost of higher prediction errors. Moreover, the MSE of the predictions of \pfyl{} is significantly higher than that of \spo{}. Further examination reveals that the higher prediction errors of \spo{} and \pfyl{} are mainly due to scale shifts in the predicted coefficient values, which does not alter the optimals of an optimization problem with a linear objective function. In addition, training with \autosklearn{}, which involves automated algorithm selection and hyperparameter tuning, is time-consuming but provides both high-quality prediction error and decision error. However, compared to \spo{} and \pfyl{}, even the competitive 2-stage Auto model does not have an advantage in decision-making with $1000$ training samples.

\begin{figure}[htbp]
    \centering
    {\includegraphics[width=0.45\textwidth]{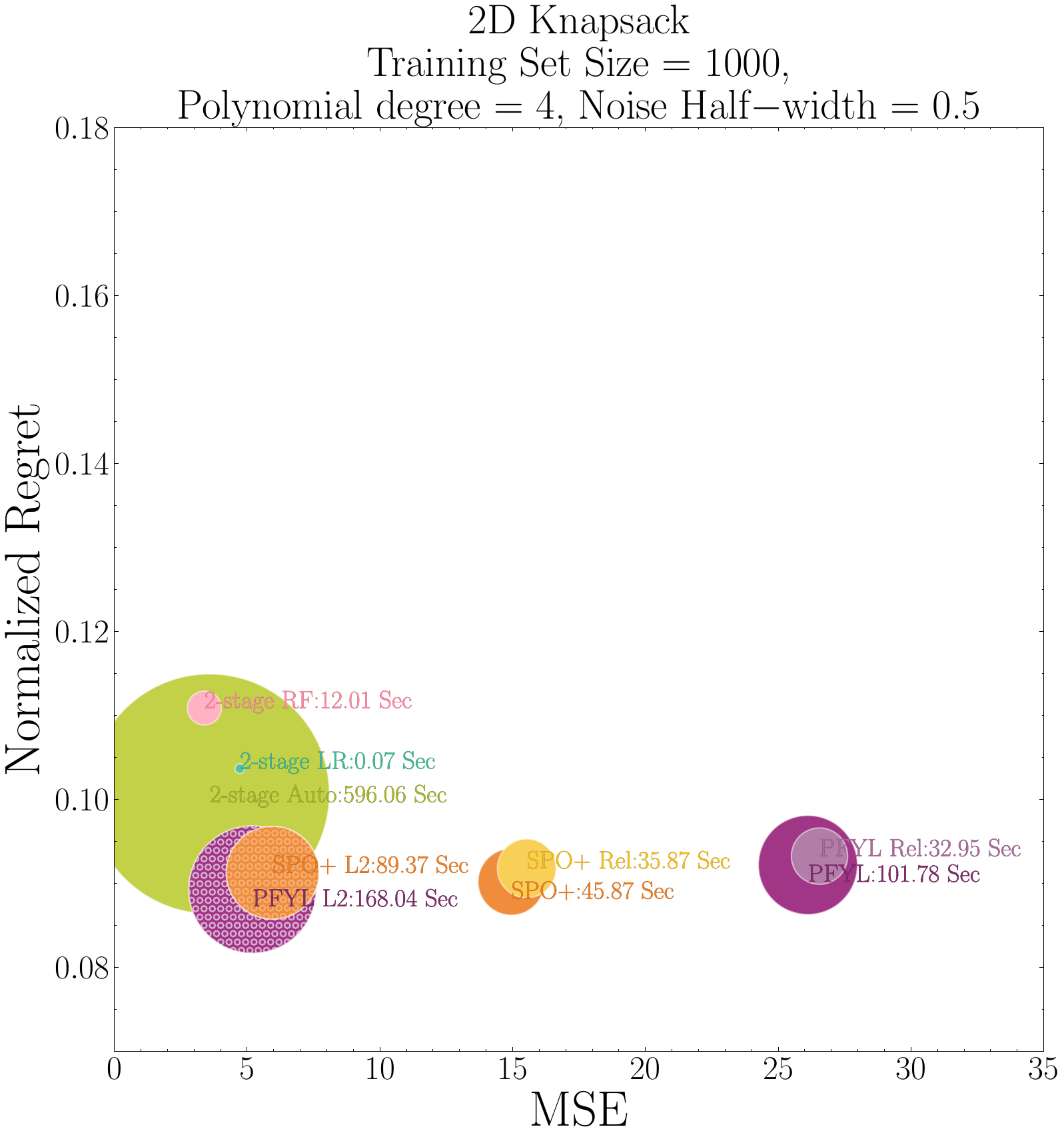}} 
    {\includegraphics[width=0.45\textwidth]{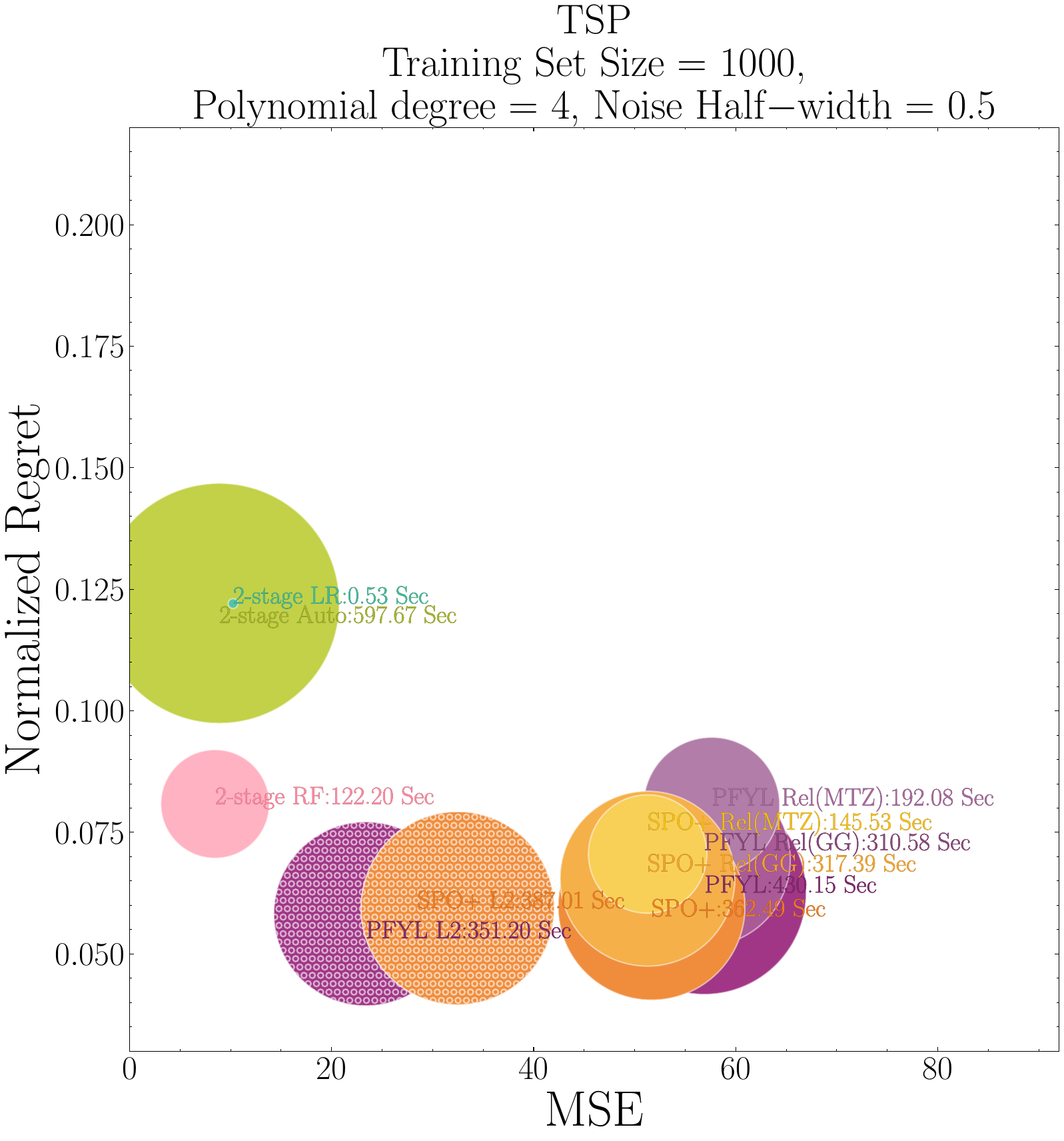}}
    \caption{MSE v.s. Regret: The result covers different two-stage methods, \spo{}, \pfyl{}, and their regularization. \dbb{} is omitted because it is far away from others. The size of the circles is proportional to the training time (Sec), so the smaller is better.}
    \label{fig:tdreg}
\end{figure}

Inspired by Sec~\ref{subsec:reg}, we further incorporated $l_2$ regularization with parameter $\phi_{2} = 1.0$ for \spo{} and \pfyl{}, refering to these as \spolt{} and \pfyllt{}. Figure \ref{fig:tdreg} shows that the appropriate amount of prediction regularization can achieve a favorable balance between MSE and regret, occasionally even reducing regret more.

\begin{tcolorbox}[colback=white,colframe=gray,title=Finding \#5]
    Generally, \spo{} and \pfyl{} can achieve good decisions , albeit with a trade-off in prediction accuracy. For a more balanced approach between decision quality and prediction accuracy, an end-to-end method with prediction regularization might be preferable.
\end{tcolorbox}

\section{Empirical Evaluation for Image-Based Shortest Path}
\label{sec:wcsp}

Following~\citet{poganvcic2019differentiation} and \citet{berthet2020learning}, we employed a truncated ResNet18 convolutional neural network (CNN) architecture consisting of the first five layers for Warcraft terrain images (see Section \ref{subsec:wcdata}). As Table~\ref{tab:cnn} shows, the methods we compare include a two-stage method, \spo{}, \dbb{}, \dpo{}, and \pfyl{}, all using truncated ResNet18.
The CNN is trained over $50$ epochs with batches of size $70$. The learning rate is set to $0.0005$, decaying at the epochs $30$ and $40$, and the hyperparameters $K = 1, \sigma = 1$ for \dpo{} and \pfyl{}, $\lambda = 10$ for \dbb{}. We use the Hamming distance for \dbb{} and the squared error of solutions for \dpo{}, as specified in their original papers.

\begin{table}[htbp]
    \begin{center}
    \resizebox{1.0\textwidth}{!}{
    \begin{tabular}{l|l}
        \hline
        \textbf{Method}     & \textbf{Description} \\ \hline
        2S                  & Two-stage method where the predictor is a truncated ResNet18 \\ \hline
        \spo{}              & Truncated ResNet18 with SPO+ loss  \cite{elmachtoub2021smart} \\ \hline
        \dbb{}              & Truncated ResNet18 with differentiable black-box optimizer and Hamming distance loss
\cite{poganvcic2019differentiation} \\ \hline
        \dpo{}              & Truncated ResNet18 with differentiable perturbed optimizer and squared error loss \cite{berthet2020learning} \\ \hline
        \pfyl{}             & Truncated ResNet18 with perturbed Fenchel-Young loss \cite{berthet2020learning} \\ \hline
    \end{tabular}}
    \end{center}
    \caption{Methods compared in the experiments.}\label{tab:cnn}
\end{table}

The sample size of the test set $n_\text{test}$ is $1000$. To evaluate our methods, we compute the relative regret $\frac{\bm{c}^\intercal \bm{w}^*(\hat{\bm{c}}) - z^*(\bm{c})}{z^*(\bm{c})}$ and path accuracy $\frac{\sum_{j=1}^d \mathbbm{1}(z^*(\bm{c})_j = z^*(\hat{\bm{c}})_j)}{d}$ per instance on the test set; the latter is simply the proportion of edges in the ``predicted" solution that also appear in the optimal solution. 

\begin{figure}[htbp]
    \centering
    \subfloat[]{\includegraphics[width=0.33\textwidth]{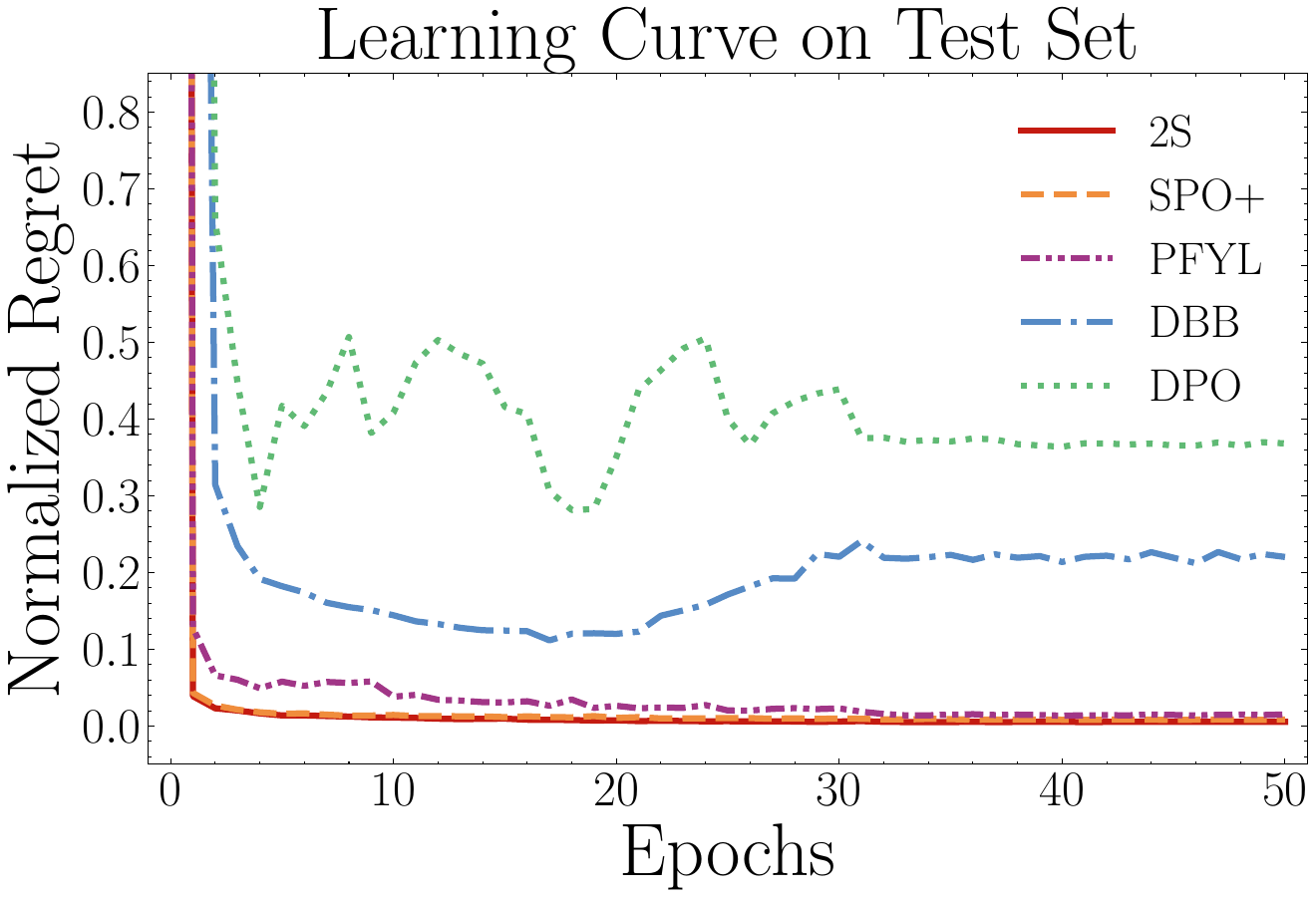}} 
    \subfloat[]{\includegraphics[width=0.33\textwidth]{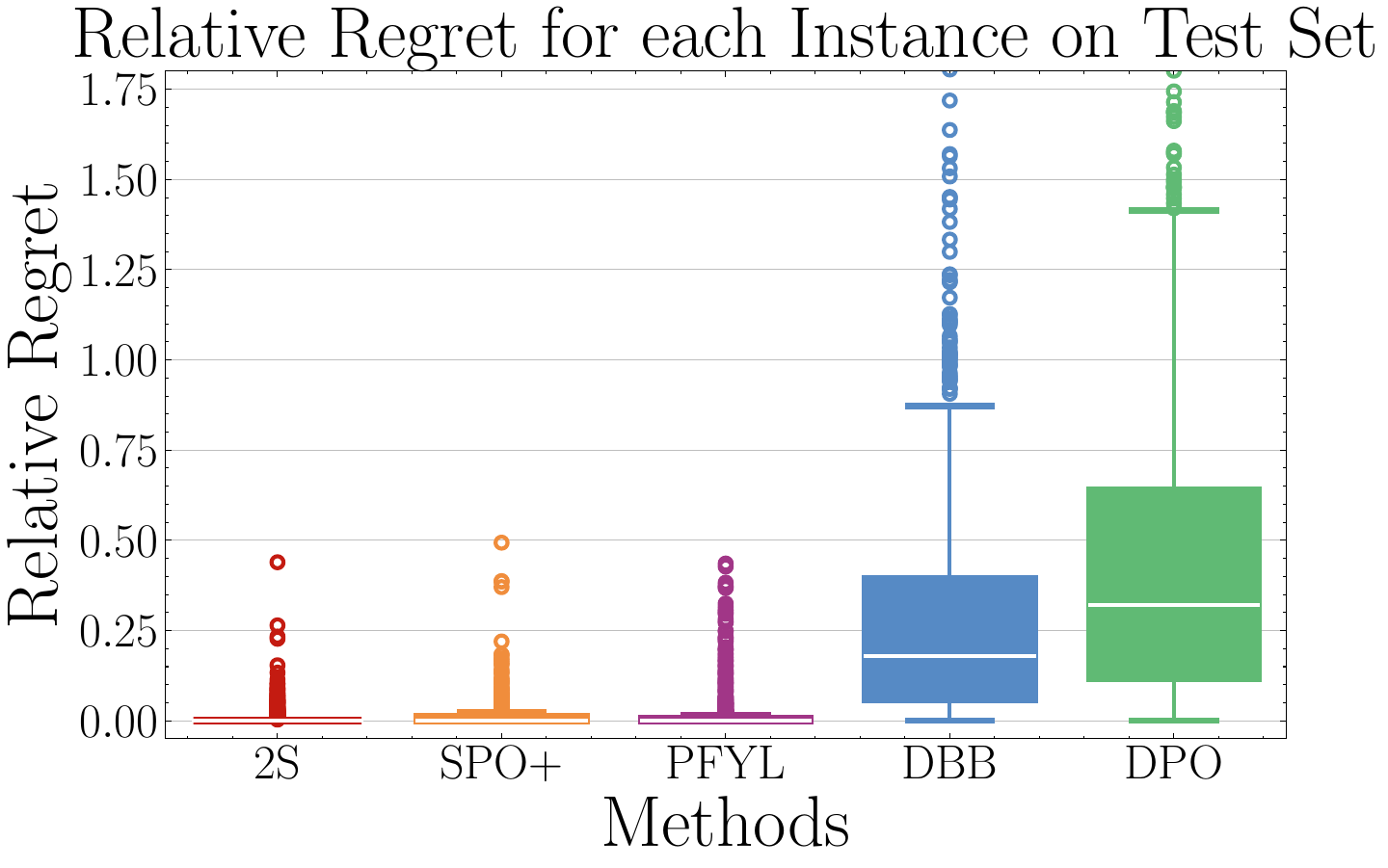}} 
    \subfloat[]{\includegraphics[width=0.33\textwidth]{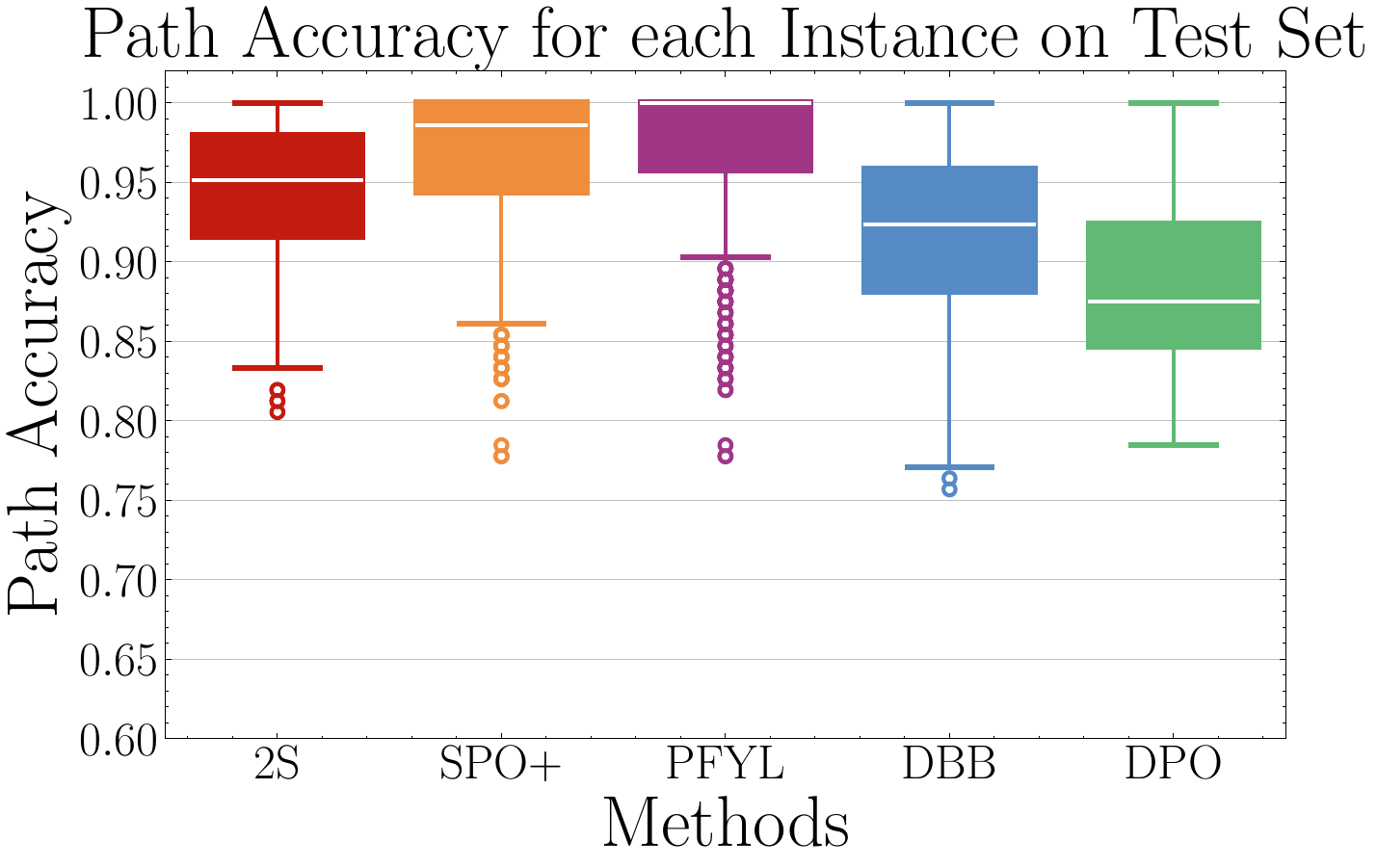}}\\
    \caption{Learning curve, relative regret, and path accuracy for the shortest path problem on the test set: The methods in the experiment include a two-stage neural network, \spo{}, \dbb{}, \dpo{}, and \pfyl{}. The learning curve shows relative regret on the test set, and the box plot demonstrates the distribution of relative regret and path accuracy on the test set. For relative regret, lower is better. For path accuracy, higher is better.}
    \label{fig:war}
\end{figure}

As shown in Figure \ref{fig:war}, the two-stage method, \spo{} and \pfyl{} achieve comparable performance in predicting the shortest path on the Warcraft terrain. Among these methods, \pfyl{} achieves the highest path accuracy, indicating that its solutions align most closely with the given optimal ones. It seems that the Warcraft shortest path problem may not require end-to-end learning. However, it is noteworthy that \pfyl{}, despite the lack of knowledge of the true coefficients, produces a competitive result. This finding encourages researchers to explore a broader range of applications for end-to-end learning.

\begin{tcolorbox}[colback=white,colframe=gray,title=Finding \#6]
    End-to-end learning is effective in rich contextual features such as images. Moreover, the study highlights that \pfyl{} can achieve impressive performance levels even without access to labeled coefficients during training.
\end{tcolorbox}

\section{Conclusion}
\label{sec:conc}

Because of the lack of easy-to-use generic tools, the potential of end-to-end predict-then-optimize has been underestimated or even overlooked in various applications. Our package ~\pyepo{} aims to alleviate the barriers between the theory and practice of the end-to-end approach, making it more accessible and practical for widespread use.

\pyepo{}, the \pytorch-based end-to-end predict-then-optimize tool, is specifically designed for linear objective functions, including linear programming and (mixed) integer programming. The tool is extended from the automatic differentiation function of \pytorch{}, one of the most widespread open-source machine learning frameworks. Hence, \pyepo{} enables users to leverage a vast array of state-of-the-art deep learning models and techniques implemented in \pytorch{}.

While \pyepo{} is currently tailored to linear objective functions and has been primarily tested on linear and linear integer programming, its underlying architecture holds the potential for broader problems. Specifically, it can be extended to accommodate nonlinear constraints and even nonlinear objective functions. For example, \spo{} in \pyepo{} requires convex constraints rather than linear ones.

\pyepo{} allows users to build optimization problems as black boxes or by leveraging interfaces to \gurobipy{}, \coptpy{} and \pyomo{}, thus providing broad compatibility with high-level modeling languages and both commercial and open-source solvers.

With the~\pyepo{} framework, we generate three synthetic datasets and incorporate one image dataset from the literature. Comprehensive experiments and analyses are conducted with these datasets. The results show that the end-to-end methods achieved excellent improvement in decision quality over two-stage methods in many cases. In addition, end-to-end models can benefit from using relaxations and regularization. The code repository includes step-by-step tutorials, enabling new users to replicate the experiments presented in this paper or use them as a starting point for their applications.
Future development efforts of \pyepo{} include: 
%adding new methods and improving current approaches. For example, \qptl{} and its variants are also popular end-to-end predict-then-optimize, but now it is absent in \pyepo{}. In addition, except relaxation, the iterative optimization procedure in both \spo{}, \pfyl{}, and \dbb{} still has the potential to be accelerated, which encourages us to investigate further.

\begin{itemize}
    \item[--] New applications and new optimization problems, including linear objective functions with mixed-integer variables or nonlinear constraints;
    \item[--] Novel training methods that leverage the gradient computation features that \pyepo{} provides;
    \item[--] Innovative approaches that accommodate nonlinear constraints and/or objective function; 
    \item[--] Additional variations and improvements to current approaches such as warm starting and training speed up;
    \item[--] Other existing end-to-end predict-then-optimize approaches, such as \qptl{} and its variants.
\end{itemize}

%\begin{acknowledgements}
%\end{acknowledgements}


\bibliographystyle{spbasic}      % basic style, author-year citations
\bibliography{ref}   % name your BibTeX database

%\begin{thebibliography}{}
% etc
%\end{thebibliography}

\section{Statements and Declarations}
\label{sec:state}

\subsection{Funding}
\label{subsec:fund}
This work was supported by funding from a SCALE AI Research Chair and an NSERC Discovery Grant to Elias B. Khalil.

\subsection{Author Contributions}
\label{subsec:contr}
Tang and Khalil contributed to the conception and design of the project. Software development, data generation, software testing, and benchmarking experiments were performed by Tang. 
Tang wrote the first draft of the submission. Tang and Khalil contributed to finalizing the submission.
Both authors read and approved the final manuscript.

\subsection{Competing Interests}
The authors have no relevant financial or non-financial interests to disclose.

\subsection{Availability of data and materials}
All data analyzed during this study are publicly available.

\subsection{Code availability}
The complete source code used in this study is open source and can be freely accessed and reviewed. We have utilized several open-source packages for our study. Specific references to these packages can be found within the code and article.

\end{document}